\documentclass[11pt,a4paper]{amsart}
\usepackage{amssymb,amscd}
\usepackage{mathrsfs}
\usepackage[mathcal]{eucal}
\usepackage{float}
\usepackage[usenames]{color}
\usepackage{MnSymbol}
\usepackage[table]{xcolor}
\usepackage[normalem]{ulem}
\usepackage{multirow}

\usepackage{ifthen}
\def \showproofsversion {0} 
\def \showproofs {\ifthenelse{\equal{\showproofsversion}{0}}}

\usepackage{tikz}
\usepackage{tikz-3dplot}
\usetikzlibrary{calc}

\definecolor{r}{HTML}{CF1020}

\newcommand{\diag}{\operatorname{diag}}

\theoremstyle{plain}
\newtheorem{theorem}{Theorem}[section]
\newtheorem{lemma}[theorem]{Lemma}

\newtheorem{condition}[theorem]{Condition}

\newtheorem{conjecture}[theorem]{Conjecture}
\newtheorem{corollary}[theorem]{Corollary}
\newtheorem{proposition}[theorem]{Proposition}

\theoremstyle{definition}
\newtheorem{definition}[theorem]{Definition}
\newtheorem{observation}[theorem]{Observation}
\newtheorem*{ackno}{Acknowledgements}

\newtheorem{beispiel}[theorem]{Example}
\newenvironment{example}{\begin{beispiel}%
 \pushQED{\qed}}%
  {\popQED\end{beispiel}}

\theoremstyle{remark}
\newtheorem{remark}[theorem]{Remark}

\numberwithin{equation}{section}

\newcommand{\N}{\mathbb{N}}
\newcommand{\Z}{\mathbb{Z}}
\newcommand{\Q}{\mathbb{Q}}
\newcommand{\R}{\mathbb{R}}
\newcommand{\C}{\mathbb{C}}
\newcommand{\F}{\mathbb{F}}

\newcommand{\GL}{\mathsf{GL}}
\newcommand{\SL}{\mathsf{SL}}
\newcommand{\Mat}{\mathsf{M}}
\newcommand{\Oh}{\mathsf{O}}

\newcommand{\G}{\mathbf{G}}
\newcommand{\Gm}{\mathsf{G}_\mathrm{m}}
\newcommand{\bv}{\mathbf{v}}

\DeclareMathOperator{\End}{End}
\DeclareMathOperator{\Hom}{Hom}
\DeclareMathOperator{\Aut}{Aut}
\DeclareMathOperator{\Span}{span}

\newcommand {\lf} {\operatorname{len}}
\newcommand {\pavbig}[1] {\big\lvert#1\big\rvert_p}

\newcommand{\Eta}{H}
\renewcommand{\rho}{\varrho} 
\renewcommand{\theta}{\vartheta}
\renewcommand{\epsilon}{\varepsilon}

\newcommand{\hlf}{\genfrac{}{}{0.1pt}{1}{1}{2}}

\renewcommand{\tt}{\top}
\newcommand{\tr}{\operatorname{tr}}

\newcommand{\cone}{\mathcal{C}}
\newcommand{\T}{t}
\newcommand{\shift}{0}

\newcommand{\xrightarrowdbl}[2][]{%
  \xrightarrow[#1]{#2}\mathrel{\mkern-14mu}\rightarrow
}

\begin{document}

\title[On pro-isomorphic zeta functions of $D^*$-groups]{On pro-isomorphic zeta functions of \\ $D^*$-groups of even Hirsch
  length}

\author{Mark N.\ Berman}\address{Department of Mathematics, Braude
  College of Engineering, P.O. Box 78, Snunit St., 51, Karmiel
  2161002, Israel} \email{berman@braude.ac.il}

\author{Benjamin Klopsch} \address{Mathematisches Institut,
  Heinrich-Heine-Universit\"at D\"usseldorf, 40225 D\"usseldorf, Germany}
\email{klopsch@math.uni-duesseldorf.de}

\author{Uri Onn} 
\address{Mathematical Sciences Institute, The Australian National University, Canberra, Australia}
\email{uri.onn@anu.edu.au}



\thanks{The research was partially conducted in the framework of the
  DFG-funded research training group GRK 2240: Algebro-Geometric
  Methods in Algebra, Arithmetic and Topology.}


\maketitle

\begin{abstract}
  The pro-isomorphic zeta function of a finitely generated nilpotent
  group is a Dirichlet generating series that enumerates all
  finite-index subgroups whose profinite completion is isomorphic to
  that of the ambient group.  We study the pro-isomorphic zeta
  functions of $\Q$-indecomposable $D^*$-groups of even Hirsch length.
  These groups are building blocks of finitely generated class-two
  nilpotent groups with rank-two centre, up to commensurability.  Due
  to a classification by Grunewald and Segal, they are parameterised
  by primary polynomials whose companion matrices define commutator
  relations for an explicit presentation.  For Grunewald--Segal
  representatives of even Hirsch length of type $f(t)=t^m$, we give a
  complete description of the algebraic automorphism groups of
  associated Lie lattices.  Utilising the automorphism groups, we
  determine the local pro-isomorphic zeta functions of groups
  associated to $t^2$ and $t^3$.  In both cases, the local zeta
  functions are uniform in the prime~$p$ and satisfy functional
  equations.  The functional equations for these groups, not predicted
  by the currently available theory, prompt us to formulate a
  conjecture which prescribes, in particular, information about the
  symmetry factor appearing in local functional equations for
  pro-isomorphic zeta functions of nilpotent groups.  Our description
  of the local zeta functions also yields information about the
  analytic properties of the corresponding global pro-isomorphic zeta
  functions. Some of our results for the $D^*$-groups associated to
  $t^2$ and $t^3$ generalise to two infinite families of class-two
  nilpotent groups that result naturally from the initial groups via
  `base extensions'.
\end{abstract}



\section{Introduction}\label{section:introduction}


\subsection{Setting the scene} Zeta functions of groups and rings were
introduced by Grunewald, Segal and Smith~\cite{GrSeSm88} as an
effective means for studying subgroup growth. Since their inception in
the late 1980s, much progress has been made regarding their analytic
and arithmetic properties; see for instance~\cite{dSGr06, Vo13}. In
this paper we focus on pro-isomorphic zeta functions. Let~$\Gamma$ be
a finitely generated nilpotent group and let $a_n^\wedge(\Gamma)$
denote the number of subgroups $\Delta \le \Gamma$ satisfying
$\lvert \Gamma : \Delta \rvert = n$ and
$\widehat{\Delta} \cong \widehat{\Gamma}$, where $\widehat{\Eta}$
denotes the profinite completion of a group~$\Eta$.  The
\emph{pro-isomorphic zeta function} of $\Gamma$ is the Dirichlet
generating series
\begin{equation*}
\zeta_\Gamma^\wedge(s) =
\sum_{n=1}^\infty a_n^\wedge(\Gamma) \, n^{-s} \qquad (s \in \C).
\end{equation*}

As with subgroup and normal subgroup zeta functions, an immediate
consequence of nilpotency is that the pro-isomorphic zeta function has
an Euler product decomposition over all rational primes:
\begin{equation} \label{equ:euler-decomp}
  \zeta_\Gamma^\wedge(s) = \prod_{p} \zeta_{\Gamma,p}^\wedge(s), \quad
  \text{where} \quad \zeta_{\Gamma,p}^\wedge(s) = \sum_{k=0}^\infty
  a_{p^k}^\wedge(\Gamma) \, p^{-ks}
\end{equation}
is called the \emph{local zeta function} at a prime $p$ and is known to
be a rational function in~$p^{-s}$ over $\Q$; see~\cite{GrSeSm88}.

A special feature of pro-isomorphic zeta functions, in contrast to
other related zeta functions of groups, is that the local zeta
functions can be expressed rather naturally as $p$-adic integrals over
algebraic groups taking the form
\begin{equation}\label{integral}
  Z(\G, p)(s) = \int_{G_p^+}  \lvert \det(g) \rvert_p^{\, s} \,
  \mathrm{d}\mu_p(g). 
\end{equation}
Here $\G \le \GL_d$ is an affine $\Z$-group scheme (the algebraic
automorphism group~$\mathbf{Aut}(L)$ of an associated nilpotent Lie
lattice~$L$), $\mu_p$ denotes a suitably normalised Haar measure on
the locally compact $p$-adic group $G_p = \G(\Q_p)$, and
$G_p^+ = G_p \cap \Mat_d(\Z_p)$ is a compact open subset of~$G_p$; the
precise details are described in Section~\ref{section:background}.

Integrals such as \eqref{integral} have a long history and were
studied for various classical groups by Hey, Weil, Tamagawa, Igusa and
others~\cite{He29, We61, Ta63, Ig89}; for a more detailed account see
~\cite{dSLu96}.  Grunewald, Segal and Smith~\cite{GrSeSm88} discovered
the relevance of such integrals for the study of pro-isomorphic zeta
functions.  Subsequently, du Sautoy and Lubotzky~\cite{dSLu96}
advanced the general theory of integrals of the form \eqref{integral}
by considering non-reductive groups~$\G$; an essential aspect of their
work was to carry out a reduction of the integral, subject to certain
technical assumptions, to an integral over a reductive
subgroup.

It is remarkable that in many cases (e.g., when the algebraic group
$\G$ is irreducibly reductive and split over~$\Q$) the zeta functions
$Z_p(s) = Z(\G,p)(s)$ are given by a single rational function in
$p, p^{-s}$ and display a symmetry upon inversion of the prime, for
almost all primes~$p$:
\[
Z_p(s)_{p \rightarrow p^{-1}} = (-1)^j p^{a-bs}Z_p(s) \qquad
\text{for suitable $a, b, j \in \N_0$.}
\]

Constructions using base extensions lead to slightly more general
situations, where the zeta functions are finitely uniform and a
corresponding finite variation $a = a(p),b = b(p)$ with~$p$ is
observed; compare with \cite[Thm.~4]{GrSeSm88}, \cite[\S3]{dSLu96} and
\cite{BGS22}. In these contexts the functional equation is a
manifestation of the compatibility of the integral with the $p$-adic
Bruhat decomposition and symmetries related to the affine Weyl group
of the reductive group~$\G$; see \cite{Ig89, dSLu96}.  Such a
phenomenon should be compared with the symmetries conjectured by Igusa
and proved by Denef and Meuser~\cite{DeMe91} for integrals over
$\Z_p^{\, d}$ of integral homogeneous polynomials, based on the
principalisation of ideals and the Weil conjectures.  More general
results in this direction, with group-theoretic applications, were
discovered and proved by Voll~\cite{Vo10}.  Since then functional
equations of the kind discussed have been recognised as a widespread,
but not universal feature of zeta functions associated to groups,
rings and modules; for instance, see
\cite{AKOV13,StVo14,Ro17,Vo19,KiKl19,dSWo08, BeKl15}.


\subsection{Main results and a conjecture}
The motivations for the present paper are two-fold. Firstly, we wish
to explore pro-isomorphic zeta functions of nilpotent groups in
situations where a crucial standard assumption, originally introduced
in~\cite{dSLu96} and until now widely used to study integrals of the
form~\eqref{integral}, does \emph{not} hold.  For this purpose, we
consider finitely generated torsion-free class-two nilpotent groups
with rank-two centres; we refer to such groups as \emph{$D^*$-groups}.
An explicit example from this family is studied in depth in this
paper, pertaining to the $D^*$-group $\Gamma_{t^3}$ of Hirsch
length~$8$, associated to the primary polynomial $t^3$; see
Theorem~\ref{x3example} below and the following discussion.  Our
analysis relies in the first place on pinning down the automorphism
group of~$\Gamma_{t^3}$.  More generally, we extend the computation,
initiated in \cite{BeKlOn18}, of the automorphism groups of
Grunewald--Segal representatives of $\Q$-indecomposable $D^*$-groups,
up to commensurability; see Theorem~\ref{thm:structure.of.Aut}.  In
addition to its inherent interest, our description of the automorphism
groups provides a first essential step toward studying the
pro-isomorphic zeta functions of more complicated $D^*$-groups; we
extend our description of the relevant automorphism groups further
in~\cite{BeKlOn21}.  Indeed, after our original work was finished,
Moadim Lesimcha and Schein~\cite{MoSc22} went ahead and studied other
families of $D^*$-groups; they produced a combinatorial description of
local pro-isomorphic zeta functions and derived local functional
equations for the families that they considered.  Secondly, we wish to
establish a conjectural framework for the shape that local functional
equations take in the context of pro-isomorphic zeta functions of
nilpotent groups, when they occur; see Conjecture~\ref{conjecture}.

We now provide more details.  In \cite[\S 6]{GrSe84}, Grunewald and
Segal considered $\mathcal{D}^*$-groups, that is, torsion-free radicable
class-two nilpotent groups of finite rank with rank-two centres.  They
classified the indecomposable constituents of such groups, by giving a
parametrisation in terms of the rank and -- in even rank -- an extra
datum, namely the projective equivalence class of an associated binary
form over~$\Q$.  Each $\mathcal{D}^*$-group is the radicable hull of a
$D^*$-group, determined up to commensurability.  We refer to such
`integral representatives' of indecomposable $\mathcal{D}^*$-groups as
\emph{$\Q$-indecomposable $D^*$-groups}.

In~\cite[Thm.~6.3]{GrSe84}, Grunewald and Segal effectively gave
explicit presentations for certain $\Q$-indecomposable $D^*$-groups,
which cover all such groups up to commensurability.  For convenience,
we refer to these special groups as \emph{Grunewald--Segal
  representatives}.  In passing, we remark that the local normal
subgroup zeta functions of such Grunewald--Segal representatives were
computed in~\cite[\S 3.2]{Vo04}.  The automorphism groups of
Grunewald--Segal representatives for $\Q$-indecomposable $D^*$-groups
of odd Hirsch length were determined in~\cite{BeKlOn18}.  In the
current paper we consider Grunewald--Segal representatives for
$\Q$-indecomposable $D^*$-groups of even Hirsch length; these are
defined explicitly in Section~\ref{section:structure-of-Aut}.  We are
particularly interested in a subfamily of $D^*$-groups $\Gamma_{t^m}$,
$m \in \N$, given by the presentations
\begin{multline} \label{equ:def-Gamma-tm}
  \Gamma_{t^m} = \langle x_1,\dots, x_m, y_1, \dots, y_m, z_1, z_2 \,
  \mid \,  [x_i,y_i]=z_1
  \text{ for $1\leq i\leq m$,} \\ 
  [x_j,y_{j+1}] = z_2 \text{ for $1\leq j<m$}, \quad [x_i,y_j]=1 \text{
    for $1 \le i,j \le m$ with $j-i \not \in \{0,1\}$,} \\ [x_i,x_j] =
  [y_i,y_j] = [x_i,z_1]= [x_i,z_2] = [y_i,z_1] = [y_i,z_2] = 1
  \text{ for $1 \le i,j \le m$} \rangle.
\end{multline}
Observe that $\Gamma_{t^m}$ has Hirsch length $2m+2$ and rank-two
centre $\mathrm{Z}(\Gamma_{t^m}) = \langle z_1, z_2 \rangle$.  For
$m=1$, the presentation yields the decomposable $D^*$-group
$\Gamma_t \cong C_\infty \times \mathrm{Heis}(\Z)$, the direct product
of an infinite cyclic group and the discrete Heisenberg group.  Its
pro-isomorphic zeta function is relatively easy to compute:
$\zeta^{\wedge}_{\Gamma_t}(s) = \zeta(s-2) \zeta(2s-3) \zeta(2s-4)$ is
a product of shifted Riemann zeta functions; this case was already
treated in~\cite[\S 3.3.4]{Be05} and we confirm the result in
Example~\ref{exa:m=1-case}.

For $m \ge 2$, the groups $\Gamma_{t^m}$ constitute one basic family
of Grunewald--Segal representatives for $\Q$-indecomposable
$D^*$-groups.  In Theorem~\ref{thm:structure.of.Aut} below we provide,
for all $m \in \N$, a complete description of the algebraic
automorphism groups of associated Lie lattices.  Based on this result,
we explicitly determine for $m \in \{2,3\}$ the corresponding
pro-isomorphic zeta functions, including all local zeta functions with
no exceptions.

\begin{theorem}\label{x2example}
  For all primes $p$, the $D^*$-group $\Gamma = \Gamma_{t^2}$
  satisfies $\zeta^{\wedge}_{\Gamma,p}(s) = W_{t^2}(p,p^{-s})$,
  where
  \[
    W_{t^2}(X,Y) = \frac{1+X^{10}
      Y^{4}}{(1-X^{8}Y^3)(1-X^{11}Y^4)(1-X^{12}Y^5)}.
  \]
  Thus $\zeta^{\wedge}_{\Gamma,p}(s)$ has abscissa of convergence
  $11/4$ and satisfies the functional equation
  \[
    \zeta^{\wedge}_{\Gamma,p}(s) \vert_{p\to p^{-1}} = (-1)p^{21-8s}\,
    \zeta^{\wedge}_{\Gamma,p}(s).
  \]
\end{theorem}

\begin{corollary} \label{cor:t2} The pro-isomorphic zeta function of
  the $D^*$-group $\Gamma = \Gamma_{t^2}$ is
  \[
    \zeta^{\wedge}_\Gamma(s) = \frac{\zeta(3s-8) \zeta(4s-11)
      \zeta(5s-12) \zeta(4s-10)}{\zeta(8s-20)},
  \]
  where $\zeta(s)$ denotes the Riemann zeta function; in particular,
  it admits meromorphic continuation to the entire complex plane and
  has abscissa of convergence $3$, with a double pole at $s=3$.

  Furthermore, the asymptotic growth of pro-isomorphic subgroups in
  $\Gamma$ is given by
  \begin{equation} \label{equ:tauberian-t2} \sum_{n=1}^N
    a_n^\wedge(\Gamma) \sim c_{t^2} N^3 \log N \qquad \text{as
      $N \to \infty$,}
  \end{equation}
  where $c_{t^2} = \frac{5 \, \zeta(3)}{12 \pi^2} \approx 0.050747$.
\end{corollary}

Theorem~\ref{x2example} and its proof resemble similar results for
other nilpotent groups, for instance the $D^*$-groups studied
in~\cite{BeKlOn18}.  In contrast, the next theorem and its proof open
up several promising new directions for further exploration.

\begin{theorem}\label{x3example}
  For all primes~$p$, the $D^*$-group $\Gamma = \Gamma_{t^3}$
  satisfies $\zeta^{\wedge}_{\Gamma,p}(s) = W_{t^3}(p,p^{-s})$, where
  \begin{multline*}
    W_{t^3}(X,Y) = \\
    \frac{(1-X^{29}Y^{10})(1 +X^{14}Y^5 -X^{15}Y^5 +X^{30}Y^{10}
      -X^{59}Y^{21} +X^{74}Y^{26} -X^{75}Y^{26}
      -X^{89}Y^{31})}{(1-X^{15}Y^5)^2 \, (1-X^{29}Y^9) \,
      (1-X^{30}Y^{11}) \, (1-X^{61}Y^{21})}.
   \end{multline*}
   Thus $\zeta^{\wedge}_{\Gamma,p}(s)$ has abscissa of convergence
   $29/9$ and satisfies the functional equation
  \[
    \zeta^{\wedge}_{\Gamma,p}(s) \vert_{p\to p^{-1}}=(-1)p^{32-10s}\,
    \zeta^{\wedge}_{\Gamma,p}(s).
  \]
\end{theorem}

\begin{corollary} \label{cor:t3} The pro-isomorphic zeta function of
  the $D^*$-group $\Gamma = \Gamma_{t^3}$ has abscissa of convergence
  $10/3$ and admits meromorphic continuation to
  $\{ s \in \C \mid \operatorname{Re}(s) > 3 \}$ via
  \[
    \zeta^{\wedge}_\Gamma(s) = \frac{\zeta(5s-15) \zeta(9s-29)
      \zeta(10s-30) \zeta(11s-30) \zeta(15s-45) \zeta(21s-61)
    }{\zeta(10s-29) \zeta(30s-90)} \; \widetilde{\psi}(s),
  \]
  where $\zeta(s)$ denotes the Riemann zeta function and
  \[
    \widetilde{\psi}(s) = \prod_p
    \frac{\widetilde{W}(p,p^{-s})}{1-p^{15-5s}+p^{30-10s}}
  \]
  for
  $\widetilde{W}(X,Y) = 1 +X^{14}Y^5 -X^{15}Y^5 +X^{30}Y^{10}
  -X^{59}Y^{21} +X^{74}Y^{26} -X^{75}Y^{26} -X^{89}Y^{31}$; moreover,
  the line $\{ s \in \C \mid \operatorname{Re}(s) = 3 \}$ is a natural
  boundary.  In particular, the zeta function
  $\zeta^{\wedge}_\Gamma(s)$ has a simple pole at $s = 10/3$.
\end{corollary}

\begin{remark} \label{rem:asymp-growth-t3} Similar to
  Corollary~\ref{cor:t2}, the asymptotic growth of pro-isomorphic
  subgroups in $\Gamma = \Gamma_{t^3}$ can be described by means of a
  suitable Tauberian theorem:
  \[
    \sum_{n=1}^Na_n^\wedge(\Gamma) \sim c_{t^3}N^{10/3} \qquad
    \text{as $N \to \infty$,}
  \]
  where
  $c_{t^3} = \frac{\zeta(5/3) \, \zeta(10/3) \, \zeta(20/3) \,
    \zeta(5) \, \zeta(9) \, \widetilde{\psi}(10/3)}{30 \, \zeta(13/3)
    \, \zeta(10)} \in \R_{> 0}$ is somewhat unwieldy.
\end{remark}

Following a suggestion of the referee, in Section~\ref{sec:base-ext}
we extend our results for the $\Q$-indecomposable $D^*$-groups
$\Gamma_{t^2}$ and $\Gamma_{t^3}$ to two infinite families,
$\widetilde{\Gamma}_{t^2,k}$ and $\widetilde{\Gamma}_{t^3,k}$ of
class-two nilpotent groups, where $k$ runs through all number fields.
These families of groups result naturally from the initial groups via
`base extensions' of corresponding Lie lattices, and pro-isomorphic
zeta functions of groups constructed in this way were systematically
investigated in~\cite{BGS22}.  For completeness we also
  discuss the family $\widetilde{\Gamma}_{t,k}$ associated to the
  decomposable $D^*$-group $\Gamma_t$.  We state here the
generalisation of Corollary~\ref{cor:t2}; further details about the
set-up and generalisations of some of our other results can be found
in Section~\ref{sec:base-ext}.

\begin{theorem} \label{thm:t2-base-change} Let $k$ be a number field
  of absolute degree $d = [k:\Q]$, with ring of
  integers~$\mathfrak{o}$.  Let
  $\widetilde{\Gamma} = \widetilde{\Gamma}_{t^2,k}$ be the class-two
  nilpotent group of Hirsch length $6d$ and with rank-$2d$ centre,
  corresponding to the class-two nilpotent $\Z$-Lie lattice
  $\widetilde{L} = \widetilde{L}_{t^2,k}$ which results from the Lie
  lattice $L = L_{t^2}$ associated to the group $\Gamma_{t^2}$ by
  extension of scalars from $\Z$ to $\mathfrak{o}$ and subsequent
  restriction of scalars back to $\Z$.

  Then the pro-isomorphic zeta function of the group
  $\widetilde{\Gamma}$ is
  \begin{equation} \label{equ:t2-zeta-Dedekind}
    \zeta^{\wedge}_{\widetilde{\Gamma}}(s) = \frac{\zeta_k(3s-(4d+4))
      \, \zeta_k(4s-(8d+3)) \, \zeta_k(5s-12d) \,
      \zeta_k(4s-(8d+2))}{\zeta_k(8s-(16d+4))},
  \end{equation}
  where $\zeta_k(s)$ denotes the Dedekind zeta function of~$k$; in
  particular, it admits meromorphic continuation to the entire complex
  plane.
\end{theorem}

\begin{remark} \label{rem:t2-base-change} For $k = \Q$, i.e.\ $d=1$,
  we recover Corollary~\ref{cor:t2}.  For quadratic fields $k$, i.e.\
  $d=2$, the abscissa of convergence is~$5$, with a double pole at
  $s=5$.  For number fields $k$ of absolute degree $d \ge 3$, the
  abscissa of convergence is $(12d+1)/5$, with a simple pole at
  $s = (12d+1)/5$.  Similar to Corollary~\ref{cor:t2}, the asymptotic
  growth of pro-isomorphic subgroups in $\widetilde{\Gamma}$ can be
  described by means of a suitable Tauberian theorem.  Via the Euler
  product, the formula~\eqref{equ:t2-zeta-Dedekind} incorporates a
  description of the local pro-isomorphic zeta functions
  $\zeta^{\wedge}_{\widetilde{\Gamma},p}(s)$ for all primes~$p$ and
  thus also yields a generalisation of Theorem~\ref{x2example}.
  Indeed, for $d \ge 2$ the zeta function
  $\zeta^{\wedge}_{\widetilde{\Gamma},p}(s)$ has abscissa of
  convergence $12d/5$ and, if $p$ is unramified in $k$, it satisfies
  the functional equation
  \[
    \zeta^{\wedge}_{\widetilde{\Gamma},p}(s) \vert_{p\to p^{-1}} = \pm
    p^{16d^2+5d-8ds}\, \zeta^{\wedge}_{\widetilde{\Gamma},p}(s).
  \]
\end{remark}
 
Theorem~\ref{x3example} and its proof extend the scope of functional
equations and the complexity of the integrals arising in the context
of pro-isomorphic zeta functions of class-two nilpotent groups.  As
alluded to above, and demonstrated in
Remark~\ref{remark:theta-not-character} below, it is the first
explicitly computed pro-isomorphic zeta function for which a certain
lifting condition \cite[Assumption~2.3]{dSLu96} does not hold.
Furthermore, it involves a technically challenging computation of an
integral with non-multiplicative integrand which requires careful
analysis by certain number-theoretic and combinatorial techniques.  In
particular, one needs to count solutions to congruence equations of
the form ${p^\alpha x^2+p^\beta y z \equiv 0 \mod p^n}$; see
Section~\ref{section:x^3}.  This reveals a new phenomenon in the
setting of pro-isomorphic zeta functions, namely the prominent role
played by counting points on reductions of varieties; previously this
feature was encountered only for other types of zeta functions of
nilpotent groups, such as subgroup and normal subgroup zeta functions;
compare with~\cite{dSGr00, dS01, Vo10}.  Our analysis of the structure
of the automorphism groups of $\Q$-indecomposable $D^*$-groups of even
Hirsch length given in Section~\ref{section:structure-of-Aut} suggests
that this is only the tip of the iceberg, and should be contrasted
with the linearity assumption in~\cite[\S 5]{dSLu96}.

The available theory on integrals of the form \eqref{integral}, which
occupy a central role in our computation, could not be used to predict
\emph{a priori} the resulting form of the local pro-isomorphic zeta
function in any sense.  It is thus somewhat of a surprise that the
zeta functions in Theorem~\ref{x3example} satisfy local functional
equations.  In contrast to the situation for $\Q$-indecomposable
$D^*$-groups of odd Hirsch length~\cite{BeKlOn18}, the values of the
abscissae of convergence -- for the pro-isomorphic zeta functions of
$\Q$-indecomposable $D^*$-groups of even Hirsch length -- remain
elusive.  More work is required, even to produce a promising
conjecture for the family of groups $\Gamma_{t^m}$,
$m \in \N_{\ge 2}$.

\medskip

In order to compare the local functional equations in
Theorems~\ref{x2example}, \ref{x3example} and their generalisations
with data for other groups, we briefly recall further concepts.  To a
finitely generated torsion-free class-$c$ nilpotent group $\Gamma$ of
Hirsch length $d$ one associates, via Lie theory, a class-$c$
nilpotent $\Z$-Lie lattice $L$ of $\Z$-rank~$d$, whose local zeta
functions
$\zeta_{L,p}^\wedge(s) = \zeta_{L_p}^{\mathrm{iso}}(s) =
\sum_{k=0}^\infty a_{p^k}^{\mathrm{iso}}(L_p) p^{-ks}$ satisfy
$\zeta_{\Gamma,p}^\wedge(s) = \zeta_{L,p}^\wedge(s)$ for almost all
primes~$p$; here $L_p = \Z_p \otimes_\Z L$ denotes the $p$-adic
completion of~$L$, and $a_{p^k}^{\mathrm{iso}}(L_p)$ is the number of
Lie sublattices of $L_p$ of index $p^k$ which are isomorphic to~$L_p$.
It was shown in~\cite{GrSeSm88} that each local zeta function
$\zeta_{L,p}^\wedge(s)$ is a rational function in $p^{-s}$ over~$\Q$,
i.e., $\zeta_{L,p}^\wedge(s) = W_p(p^{-s})$ for suitable
$W_p = R_p/Q_p$ with $R_p,Q_p \in \Q[Y]$. We then define the
\emph{degree} of a local pro-isomorphic zeta function, denoted by
$\deg_{p^{-s}} \zeta_{L,p}^\wedge(s)$, to be the degree of the
rational function $W_p$, viz.\ $\deg W_p = \deg_Y R_p - \deg_Y Q_p$.
The family of local zeta functions $\zeta_{L,p}^\wedge(s)$ is said to
be \emph{finitely uniform} if there exist finitely many rational
functions $\mathcal{W}_1, \ldots, \mathcal{W}_r \in \Q(X,Y)$ in two
variables such that, for each prime~$p$, there is an index $i = i(p)$
for which
$\zeta_{L,p}^\wedge(s) = W_p(p^{-s}) = \mathcal{W}_i(p,p^{-s})$.

Another ingredient relates to the nilpotent $\Z$-Lie lattice $L$
itself: recall that $L$ is \emph{$\N$-graded} if it is equipped with
an additive decomposition $L = \bigoplus_{i \in \N} L_{(i)}$ such that
$[L_{(i)}, L_{(j)}] \subseteq L_{(i+j)}$ for all $i, j \in \N$; for
short, we refer to the latter as a grading on~$L$.  Since $L$ has
finite rank as a $\Z$-module, there exists, for a given grading, a
minimal $l \in \N_0$ such that $L_{(j)}=0$ for $j>l$; the grading then
gives rise to a descending filtration
$L = \overline{L}_{(1)} \supseteq \overline{L}_{(2)} \supseteq \cdots
\supseteq \overline{L}_{(l)} \supseteq \{0\}$ of $L$ by Lie
sublattices
$\overline{L}_{(i)} = \sum_{j=i}^l L_{(j)} \supseteq \gamma_i(L)$.  We
call a grading \emph{natural} if its associated filtration is
precisely the lower central series, i.e., if
$\overline{L}_{(i)} = \gamma_i(L)$ for $1 \le i \le l$ and $l=c$ is
the nilpotency class of~$L$. To a grading on $L$ as above we attach a
\emph{weight} given by
$\sum_{i=1}^l i\ \operatorname{rk}_\Z L_{(i)} = \sum_{i=1}^l
\operatorname{rk}_\Z \overline{L}_{(i)}$, and we call a grading
\emph{minimal} if its weight is minimal amongst all weights of
gradings on~$L$.  In passing, we mention that not all nilpotent Lie
lattices admit a grading.  For instance, Dyer~\cite{Dy70} constructed
a $9$-dimensional class-$6$ nilpotent Lie algebra over $\Q$ whose
algebraic automorphism group is unipotent.  This implies that the Lie
algebra does not possess any grading, since every non-zero graded Lie
algebra admits non-trivial semisimple automorphisms; clearly, no Lie
lattice in such a Lie algebra can possess a grading.

\begin{conjecture}\label{conjecture}
  Let $L$ be a nilpotent $\Z$-Lie lattice that admits at least one
  grading. Then, for almost all primes~$p$, the degree of the local
  pro-isomorphic zeta function of $L$ at $p$ is equal to the weight of
  a minimal grading of~$L$.

  In particular, if the family of local pro-isomorphic zeta functions
  $\zeta^\wedge_{L,p}(s)$ is finitely uniform and the local zeta
  functions satisfy, for almost all primes~$p$, functional equations
  of the form
  \[
    \zeta^\wedge_{L,p}(s) \vert_{p\to p^{-1}}=(-1)^j p^{a-bs} \,
    \zeta^\wedge_{L,p}(s) \qquad \text{for suitable
      $a = a(p), b = b(p), j = j(p) \in \N_0$,}
  \]
  then the integer $b$ in the `symmetry factor' is the same for
  almost all $p$ and is given by the weight of a minimal grading
  of~$L$.
\end{conjecture}

\begin{remark}
  Note that natural gradings, when they exist, are minimal. It follows
  that, if a class-$c$ nilpotent Lie lattice $L$ is naturally graded,
  then -- in accordance with the conjecture -- we expect that
  $\deg_{p^{-s}} \zeta_{L,p}^\wedge(s) = \sum_{j=1}^c
  \operatorname{rk}_\Z \gamma_j(L)$ for almost all primes~$p$.  It is
  curious that this expression already has an interpretation in
  asymptotic group theory: it provides the degree of polynomial word
  growth of finitely generated nilpotent groups~$\Gamma$ giving rise
  to~$L$ via Lie theory; see~\cite{Ba72}.  In particular, every
  class-two nilpotent Lie lattice $L$ is naturally graded and thus we
  expect that the degrees satisfy
  $\deg_{p^{-s}} \zeta_{L,p}^\wedge(s) = \operatorname{rk}_\Z L +
  \operatorname{rk}_\Z [L,L]$ for almost all primes~$p$.
\end{remark}

In spirit, Conjecture~\ref{conjecture} is similar to part of a
conjecture of Voll on submodule zeta
functions~\cite[Conj.~1.11]{Vo19}, but the conjectures involve
different types of filtrations (which can be seen already for the
group $\Gamma_t$, arising from \eqref{equ:def-Gamma-tm} for $m=1$) and
as yet there is no direct link between the two.
We have tested Conjecture~\ref{conjecture} comprehensively for all
nilpotent $\Z$-Lie lattices $L$ for which the local pro-isomorphic
zeta functions are known; this list includes many naturally graded Lie
lattices as well as some Lie lattices not possessing a natural
grading; we refer to \cite{GrSeSm88, Be05, BeKlOn18,
    BGS22,MoSc22} for descriptions of relevant nilpotent $\Z$-Lie
lattices and their pro-isomorphic zeta functions.  The current paper
provides two new infinite families of groups confirming the
conjecture: the integers $b$ in the symmetry factors of the local zeta
functions described in Remarks~\ref{rem:t2-base-change} and
\ref{rem:t3-base-change} indeed match the sum of the ranks of terms of
the lower central series: for the `base extensions' defined in
Theorems~\ref{thm:t2-base-change} and \ref{thm:t3-base-change} one has
$8d = 6d + 2d$ and $10d = 8d + 2d$ for all primes unramified in the
extension.

Our conjecture also holds true for a $\Z$-Lie lattice $L$, constructed
by Berman and Klopsch in~\cite{BeKl15}, with the property that its
local pro-isomorphic zeta functions $\zeta_{L,p}^\wedge(s)$ do
\emph{not} satisfy functional equations for $p>3$.  The relevant Lie
lattice $L$ is not naturally graded, but admits a minimal grading of
weight~$102$; and, indeed, the local zeta functions are uniform
in~$p$, for $p>3$, of degree~$102$.  This example can also be
generalised by means of base extensions; see~\cite{BGS22}.

It is well known and easy to see that there is a link between the
existence of gradings of a $\Z$-Lie lattice $L$ and the occurrence of
diagonalisable elements in the algebraic automorphism group
$\mathbf{Aut}(L)$ of~$L$.  Conjecture~\ref{conjecture} suggests that
there is a somewhat more delicate connection (yet to be discovered)
between minimal gradings of a nilpotent $\Z$-Lie lattice $L$ and the
degrees of its local pro-isomorphic zeta functions, which stand in
close relation to $\mathbf{Aut}(L)$ as indicated in~\eqref{integral}.

\medskip

In order to carry out the computations leading to
Theorems~\ref{x2example} and~\ref{x3example} and their generalisations
we require a structural description of the relevant automorphism
groups.  In fact, we determine the algebraic automorphism groups for
the Lie lattices associated to Grunewald--Segal representatives of
$\Q$-indecomposable $D^*$-groups of even Hirsch length associated to
the primary polynomials $\Delta(t)=t^m$, for all $m \in \N$; as in the
case of odd Hirsch length~\cite{BeKlOn18}, this structure theorem for
the algebraic automorphism groups is of independent interest.  The
presentation \eqref{equ:def-Gamma-tm} for the group $\Gamma_{t^m}$
readily translates into a description~\eqref{eq:structure.of.L} of the
corresponding Lie lattice; compare with
Section~\ref{section:class-2-correspondence}.

\begin{theorem}\label{thm:structure.of.Aut}
  For $m \in \N$, let $\G \le \GL_{2 m + 2}$ be the algebraic
  automorphism group of the $\Z$-Lie lattice (scheme) $L$ associated,
  via~\eqref{eq:structure.of.L} below, to the primary polynomial
  $\Delta(t) =t^m$.  Let $\G_0 \trianglelefteq \G$ be the affine
  subgroup consisting of all automorphisms that fix pointwise the
  centre of~$L$. Then $\G$ splits as
  \[
    \G \cong  \mathbf{B}_2 \ltimes \G_0,
  \]
  where, for every field extension $k$ of $\Q$, the group
  $\mathbf{B}_2(k)$ is the group of invertible lower-triangular
  $2 \times 2$ matrices, and
  \[
    \G_0(k) \cong \SL_2 \big( R \big) \ltimes
    \mathsf{V}_\mathrm{st}(R)^{\oplus 2}, \quad \text{for
      $R = k[t]/(t^m)$ and $\mathsf{V}_\mathrm{st}(R) = R^2$,}
  \]
  with respect to the standard left action. In particular, the
  algebraic group $\G$ is connected.
\end{theorem}

\begin{remark}
  In fact, the description of $\G_0$ given in
  Theorem~\ref{thm:structure.of.Aut} holds true more generally, for
  $\Z$-Lie lattices corresponding to arbitrary primary polynomials;
  see Theorem~\ref{explicit.G0} below.  The description of the
  quotient of $\G$ by $\G_0$, however, becomes more involved;
  see~\cite{BeKlOn21}.
\end{remark}

The proof of Theorem~\ref{thm:structure.of.Aut}, along with explicit
forms of the automorphism groups, is given in
Section~\ref{section:structure-of-Aut}.  Our considerations in this
context overlap somewhat with the treatment
in~\cite{BrOB08}. In~\cite{BeKlOn21} we give a complete description of
the algebraic automorphism groups of all $\Z$-Lie lattices associated
to Grunewald--Segal representatives of $\Q$-indecomposable
$D^*$-groups of even Hirsch length, based on a more technical analysis
of the Lie algebras associated to (subgroups of) the algebraic
automorphism groups.


\subsection{Layout of the paper}
In Section~\ref{section:structure-of-Aut} we analyse and describe the
algebraic automorphism groups of $\Z$-Lie lattices associated to
indecomposable $D^*$-groups of even Hirsch length, corresponding to
primary polynomials of the form $\Delta(t) = t^m$.  In
Section~\ref{section:background} we provide technical background
regarding conditions on the algebraic automorphism group of a Lie ring
that is needed for calculating pro-isomorphic zeta functions of
groups.  In Sections~\ref{section:x^2} and~\ref{section:x^3} we
present calculations of the local pro-isomorphic zeta functions of the
groups $\Gamma_{t^2}$ and $\Gamma_{t^3}$.  The former group can be
dealt with in a quite straightforward manner, while the latter group
is considerably more difficult to handle. From the description of the
local zeta functions we draw conclusions about the analytic behaviour
of the global pro-isomorphic zeta functions of $\Gamma_{t^2}$ and
$\Gamma_{t^3}$; again the treatment of the latter group, which forms
Section~\ref{sec:mero-cont}, is more challenging and displays
interesting features.  In Section~\ref{sec:base-ext} we extend our
results for the groups $\Gamma_{t^2}$ and $\Gamma_{t^3}$ to two
infinite families of class-two nilpotent groups that result via `base
extensions' of corresponding Lie lattices.


\subsection{Basic notation} We denote by $\N_0$ and $\N$ the
non-negative and positive integers, respectively.  For $S\subseteq\R$
and $a \in \R$ we write $S_{\geq a} = \{x \in S \mid x\geq a\}$, and
similarly for $S_{>a}$.  For a prime $p$, we write $\Q_p$ for the
field of $p$-adic numbers with $\Z_p$ its ring of integers.  We denote
the $p$-adic valuation of $x\in \Q_p$ by $v_p(x)$ and write
$\lvert x \rvert_p = p^{\, - v_p(x)}$ for the $p$-adic absolute value.
A Lie lattice over a commutative ring $R$ with $1$ is a finitely
generated free $R$-module, equipped with a suitable Lie bracket.

\begin{ackno}
  The first author thanks Braude College of Engineering for travel grants.  We
  thank Moritz Petschick for technical help with implementing
  Figure~\ref{figure1}.  We are grateful for the referee's feedback
  that led to several improvements in the exposition and prompted us
  to work out the generalisations in Section~\ref{sec:base-ext}.
\end{ackno}



\section{Automorphism groups of $\Q$-indecomposable $D^*$-Lie lattices}\label{section:structure-of-Aut}

For any commutative ring $R$ with~$1$ and any free $\Z$-module~$M$, we
use the notation ${_R M} = R \otimes_\Z M$ to denote the free
$R$-module obtained by extension of scalars; if $M$ carries extra
algebraic structure that is compatible with extension of scalars, such
as the structure of a Lie lattice, we employ the same notation.  Thus
a $\Z$-Lie lattice $L$ sets up a Lie lattice scheme
$R \leadsto {_R L}$.  We realise the \emph{algebraic automorphism
  group} $\mathbf{Aut}(L)$ of~$L$, via a $\Z$-basis of~$L$, as an
affine $\Z$-group scheme $\G \le \GL_d$, where $d = \dim_\Z(L)$ is the
$\Z$-rank of~$L$, so that, in particular,
\[
   \Aut({_k L}) \cong \G(k) \le \GL_{d}(k) \qquad \text{for every
    extension field $k$ of~$\Q$,}
\]
and, thinking of $\GL_d$ as a subgroup of $\SL_{d+1}$ to make the
arithmetic structure tangible,
\[
   \Aut(L) \cong \G(\Z) \quad \text{and} \quad 
  \Aut({_{\Z_p} L}) \cong \G(\Z_p) \text{ for each prime $p$,}
\]
with respect to the chosen basis.  The automorphism groups arising in
this paper come from nilpotent $\Z$-Lie lattices with rank-two centres
and, for short, we refer to these as $D^*$-Lie lattices.  Our aim here
is to describe the algebraic automorphism groups of
$\Q$-indecomposable $D^*$-Lie lattices of even $\Z$-rank $d=2m+2$
which admit a presentation suggested by~\cite[Thm.~6.3(b)]{GrSe84} and
associated with the primary polynomial $\Delta(t)=t^m$; compare with
Section~\ref{section:class-2-correspondence}.  The corresponding task
for $D^*$-Lie lattices of odd $\Z$-rank has been carried out
in~\cite{BeKlOn18}; the case of more general $D^*$-Lie lattices of
even $\Z$-rank is considered in~\cite{BeKlOn21} (and turns out to be
more involved).

We now give a detailed description, in coordinates, that is tailored
also to our investigations of pro-isomorphic zeta functions.  Let
$m \in \N$ and consider the companion matrix
\begin{equation} \label{equ:matrix-K}
  K=C(a_1,\dots,a_{m}) =
    \begin{pmatrix}
      0 & 1 & 0 & \cdots & 0\\
      0 & 0 & 1 & \cdots & 0\\
      \vdots & \vdots & \vdots & \ddots & \vdots\\
      0 & 0 & 0 & \cdots & 1\\
      a_m & a_{m-1} & a_{m-2} & \cdots & a_{1}
    \end{pmatrix}
  \in \Mat_m(\Z)
\end{equation}
of a monic polynomial
\[
\Delta_K = t^m-a_{1}t^{m-1}- \dots - a_{m-1} t -a_m \in \Z[t].
\]
We consider the $\Z$-Lie lattice $L$ of $\Z$-rank $2m+2$ with ordered
$\Z$-basis
\[
  \mathcal{S} = (x_1, \dots, x_m,\, y_1, \dots, y_m,\, z_1,z_2)
\]
and the Lie bracket defined by
\begin{multline}\label{eq:structure.of.L} 
  [x_i,y_j] =\delta_{i,j} z_1 + K_{ij}z_2 \quad \text{and} \quad 
  [x_i,x_j] = [y_i,y_j] = [x_i,z_1] = [x_i,z_2] = [y_i,z_1] =
  [y_i,z_2] = 0, \\
\text{for $1 \le i,j \le m$,}  
\end{multline}
where $\delta_{i,j}$ denotes the Kronecker-delta.  We observe that $L$
is a $D^*$-Lie lattice with centre
\[
  Z = \mathrm{Z}(L) = \Z z_1 + \Z z_2.
\]

Let $\G \le \GL_{2m+2}$ be the algebraic automorphism group of $L$
with the embedding defined by the ordered basis~$\mathcal{S}$.  In
particular, for every integral domain $k$ of characteristic~$0$, the
coordinate maps with respect to $\mathcal{S}$ identify ${_k L}$ with
the module $k^{2m+2}$ of row-vectors, and the action of the group
$\G(k) = \Aut({_k L}) \le \GL_{2m+2}(k)$ on ${_k L}$ corresponds to
matrix multiplication from the right.  We write
$\G_0 \trianglelefteq \G$ for the affine subgroup and $\Z$-subscheme
arising as the kernel of the natural restriction homomorphism
\begin{equation} \label{equ:restr-hom}
  \G_0(k) = \operatorname{Ker} \big( \G(k)
  \xrightarrow{\operatorname{Res}^L_Z} \GL({_k Z} ) \big).
\end{equation}

From now on without further reference, let $k$ denote an integral
domain of characteristic~$0$.  Recall that an $n \times n$ matrix over
$k$ is \emph{regular} (or cyclic) over~$k$, if it is similar over $k$
to a companion matrix; such a matrix yields a regular element of the
Lie lattice $\mathsf{gl}_n(k)$, i.e., an element whose centraliser has
the smallest possible rank~$n$.  The fact that the matrix~$K$ is
regular plays a central role in the elucidation of~$\G$, and it is
convenient to note down two elementary facts.

\begin{remark}\label{remark.regular.matrices}
  Let $X, Y \in \Mat_n(k)$ be regular $n \times n$ matrices over~$k$.
  Then
  \begin{enumerate}
  \item The centraliser of $X$ is the polynomial algebra that it
    generates: $\mathrm{C}_{\Mat_n(k)}(X) = k[X]$.
  \item If $X$ and $Y$ have the same characteristic polynomial, then
    $X$ and $Y$ are similar over~$k$.
  \end{enumerate}
\end{remark}


\medskip

The Lie bracket of ${_k L}$ induces an anti-symmetric bilinear map
\begin{equation}\label{eqn:bilinear.form.Z} [\cdot,\cdot] \colon\, {_k
    L}/{_k Z} \times {_k L}/ {_k Z} \longrightarrow {_k Z}
\end{equation}
with values in ${_k Z}$ which, by a slight abuse of notation, we
continue to denote by $[\cdot,\cdot]$. The structure of $\G(k)$ is
tightly connected with the symmetries of two $k$-valued bilinear forms
on the free module~$k^{2m} \cong {_k L} / {_k Z}$ that can be derived
from the map described in~\eqref{eqn:bilinear.form.Z}.  For any matrix
$Q \in \Mat_m(\Z)$, the matrix
\[
  \mathcal{J}_Q =
  \begin{pmatrix}
    0 & Q \\ -Q^\tt & 0
  \end{pmatrix}
  \in \Mat_{2m}(\Z)
\]
can be regarded as the structure matrix of an anti-symmetric bilinear
form $\langle \cdot,\cdot \rangle_{\mathcal{J}_Q}$ on $k^{2m}$.
Let~$\Oh_{\mathcal{J}_Q} \le \GL_{2m}$ be the affine $\Z$-group scheme
such that $\Oh_{\mathcal{J}_Q}(k)$ consists of all elements of
$\GL_{2m}(k)$ that preserve the form
$\langle \cdot,\cdot \rangle_{\mathcal{J}_Q}$, that is,
\begin{align*}
  \Oh_{\mathcal{J}_Q}(k) %
  & = \{ g \in \GL_{2m}(k) \mid g \,\mathcal{J}_Q\, g^\tt=\mathcal{J}_Q\}.
\end{align*}
We remark that, if $Q = \mathrm{I}_m$ is the identity matrix, the
group scheme $\Oh_{\mathcal{J}_{\mathrm{I}_m}}$ is simply the
classical symplectic group $\mathsf{Sp}_{2m}$.
  

\subsection{The structure of the algebraic subgroup
  $\G_0$}\label{sec:structure.of.aut.gp}
We start with the structure of $\G_0 \trianglelefteq \G$, the
algebraic subgroup and $\Z$-subscheme, whose group of $k$-points
$\G(k)$ fixes the centre $_k Z = \mathrm{Z}({_k L})$ pointwise. An
element $g \in \G(k) \le \GL_{2m+2}(k)$ can be written as a block
matrix
\begin{equation}\label{equ:g-in-block-form}
  g=\begin{pmatrix} X & U \\ 0 & Y\end{pmatrix}, \quad \text{with}
  \quad 
  X = \left(\begin{smallmatrix} A & B \\ C &
      D \end{smallmatrix}\right) \in \GL_{2} \big( \Mat_m(k) \big),
  \,\,   Y = \left(\begin{smallmatrix} a & b \\
      c & d \end{smallmatrix}\right) \in \GL_2(k), \,\, U
  \in \Mat_{2m,2}(k),
\end{equation}
where $X$ and $Y$ correspond to the automorphisms that $g$ induces
naturally on ${_k L}/{_k Z}$ and ${_k Z}$.  Each of the following
equivalent conditions characterises elements of $\G(k)$ among
arbitrary elements $g$ of the form~\eqref{equ:g-in-block-form}:
\begin{align*}
  & [u,v]g=[ug,vg] && \text{for all $u,v \in {_k L}$;}  \\
  & [\bar{u},\bar{v}]Y = [\bar{u}X,\bar{v}X] %
                   && \text{for all $\bar{u},\bar{v} \in {_k L}/{_k
                      Z}$;} 
\end{align*}
\begin{equation}\label{eq:Group.action.in.terms.of.J}
  a \mathcal{J}_{\mathrm{I}_m} + c \mathcal{J}_K =
  X\mathcal{J}_{\mathrm{I}_m} X^\tt \quad \text{and} \quad
  b \mathcal{J}_{\mathrm{I}_m} + d \mathcal{J}_K = X\mathcal{J}_K X^\tt.
\end{equation}
From \eqref{eq:Group.action.in.terms.of.J} we directly obtain a
characterisation of $\G(k)$.
\begin{proposition}\label{prop:characterisation.of.G.and.G0} Let
  $g \in \GL_{2m+2}(k)$ be a block matrix of the
  form~\eqref{equ:g-in-block-form}.  
  Then
  \begin{enumerate}
  \item $g \in \G(k)$ if and only if the following four conditions are satisfied:
    \begin{enumerate}
    \item[(i)] $BA^\tt=AB^\tt$ and $BK^\tt A^\tt=AK B^\tt$,
    \item[(ii)] $C D^\tt=D C^\tt$ and $CKD^\tt=D K^\tt C^\tt$,
    \item[(iii)] $a \mathrm{I}_m +cK=A D^\tt -B C^\tt$,
    \item[(iv)] $b \mathrm{I}_m + dK=AKD^\tt-BK^\tt C^\tt$.
    \end{enumerate}
  \item $g \in \G_0(k)$ if and only if $Y = \mathrm{I}_2$ and
    $X \in \Oh_{\mathcal{J}_{\mathrm{I}_m}}(k) \cap
    \Oh_{\mathcal{J}_K}(k)$, or explicitly: $\left(\begin{smallmatrix} a & b \\
        c & d \end{smallmatrix}\right) = \mathrm{I}_2$ and
    \begin{enumerate}
    \item[$\mathrm{(i)\phantom{_0}}$] $BA^\tt=AB^\tt$ and $BK^\tt
      A^\tt=AK B^\tt$, 
    \item[$\mathrm{(ii)\phantom{_0}}$] $C D^\tt=D C^\tt$ and
      $CKD^\tt=D K^\tt C^\tt$, 
    \item[$\mathrm{(iii)_0}$] $\mathrm{I}_m =A D^\tt -B C^\tt$,
    \item[$\mathrm{(iv)_0}$] $K=AKD^\tt-BK^\tt C^\tt$.
    \end{enumerate}
    
\end{enumerate}
\end{proposition}

The proof of the following key theorem was inspired by a more
technical analysis of the Lie algebras associated to subgroups of
$\G$, carried out in~\cite{BeKlOn21}, and by-passes the use of Lie
algebras by means of a computational trick.

\begin{theorem}\label{explicit.G0}  The affine group scheme $\G_0$
  splits as follows:
  $\G_0(k) \cong \SL_2(k[K]) \ltimes
  \mathsf{V}_\mathrm{st}(k[K])^{\oplus 2}$, where
  $\mathsf{V}_\mathrm{st}(\cdot)$ denotes the standard left
  $\SL_2(\cdot)$-module.
\end{theorem}

\begin{proof}
  Recall that every square matrix over a field is similar to its
  transpose and that the conjugating matrix may be taken to be
  symmetric.  In fact, for regular matrices it is always symmetric;
  compare with~\cite{TaZa59}.  Therefore, there exists a symmetric
  matrix $\sigma \in \GL_m(\Q)$ such that
  $K^\tt = \sigma K \sigma^{-1}$.  In our special situation, we can
  even arrange that $\sigma \in \GL_m(\Z)$, because the factor groups
  $\Z^m$ modulo the row-span of $K$ and $\Z^m$ modulo the column-span
  of $K$ are isomorphic (cyclic) groups.  We set
  \begin{equation} \label{equ:def-Sigma}
    \Sigma = \left(\begin{matrix} \mathrm{I}_m & \\ &
        \sigma \end{matrix}\right) \in  \GL_{2m}(\Z), \qquad
    \text{where $K^\tt = \sigma K \sigma^{-1}$,}
  \end{equation}
  and claim that for every $g \in \GL_{2m+2}(k)$ of the
  form~\eqref{equ:g-in-block-form}, with $Y = \left(\begin{smallmatrix} a & b \\
      c & d \end{smallmatrix}\right) = \mathrm{I}_2$, the following
  holds:
  \begin{equation} \label{equ:G0-SL2} g \in \G_0(k) \qquad \text{if
      and only if} \qquad \Sigma^{-1} X \Sigma \in \SL_2(k[K]).
  \end{equation}
  
  First suppose that $g \in \G_0(k)$.  By
  Proposition~\ref{prop:characterisation.of.G.and.G0}~(2), this
  implies that
  $X = \left(\begin{smallmatrix} A & B \\ C &
      D \end{smallmatrix}\right) \in \GL_{2} \big( \Mat_m(k) \big)$
  satisfies conditions $\mathrm{(i)}$--$\mathrm{(iv)}_0$.  From (i),
  (ii), $\mathrm{(iii)_0}$ and $\mathrm{(iii)}_0^\tt$ -- the transpose
  of $\mathrm{(iii)_0}$ -- we obtain
  \begin{equation}
    \left(\begin{matrix} A & B \\ C & D \end{matrix}\right)^{-1} =
    \left(\begin{matrix} D^\tt & -B^\tt \\ -C^\tt &
        A^\tt \end{matrix}\right). 
  \end{equation}
  Now, using the fact that the inverse $g^{-1} \in \G_0(k)$
   satisfies a similar set of equations, we get
  \begin{enumerate}
  \item[$\mathrm{(i)}'$] $B^\tt D=D^\tt B$ and
    $B^\tt K^\tt D=D^\tt K B$,
  \item[$\mathrm{(ii)}'$] $C^\tt A=A^\tt C$ and
    $C^\tt K A=A^\tt K^\tt C$,
  \item[$\mathrm{(iii)_0'}$] $\mathrm{I}_m =D^\tt A-B^\tt C$,
  \item[$\mathrm{(iv)_0'}$] $K=D^\tt K A-B^\tt K^\tt C$.
  \end{enumerate}
  Using these additional conditions we deduce that
  \begin{equation}\label{equ:4-cons}
    AK=KA, \qquad K B = B K ^\tt, \qquad C K = K^\tt C, \qquad K^\tt
    D=D K^\tt.
  \end{equation}
  Indeed, multiplying $\mathrm{(iv)_0'}$ by $A$ on the left gives
  \begin{multline*}
    AK = AD^\tt KA - AB^\tt K^\tt C
    \overset{\mathrm{(iii)_0}}{=}(\mathrm{I}_m+BC^\tt) KA - A B^\tt
    K^\tt C
    \overset{\mathrm{(ii)'}}{=} KA+BA^\tt K^\tt C - A B^\tt K^\tt C \\
    \overset{\mathrm{(i)}}{=} KA;
  \end{multline*}
  multiplying $\mathrm{(iv)_0'}$ by $B^\tt$ on the right gives
  \begin{multline*}
    KB^\tt = D^\tt KA B^\tt - B^\tt K^\tt C B^\tt
    \overset{\mathrm{(iii)_0^\tt}}{=}D^\tt KA B^\tt - B^\tt
    K^\tt(DA^\tt-\mathrm{I}_m) \\ \overset{\mathrm{(i)}}{=} D^\tt K B
    A^\tt -B^\tt K^\tt D A^\tt +B^\tt K^\tt \overset{\mathrm{(i)'}}{=}
    B^\tt K^\tt;
  \end{multline*}
  multiplying $\mathrm{(iv)_0'}$ by $C$ on the left gives
  \begin{multline*}
    CK = CD^\tt KA - CB^\tt K^\tt C
    \overset{\mathrm{(iii)_0^\tt}}{=}CD^\tt KA - (DA^\tt-\mathrm{I}_m)
    K^\tt C \overset{\mathrm{(ii)}}{=} D C^\tt K A - DA^\tt K^\tt C
    +K^\tt C \\ \overset{\mathrm{(ii)'}}{=} K^\tt C;
  \end{multline*}
  and multiplying $\mathrm{(iv)_0'}$ by $D^\tt$ on the right gives
  \begin{multline*}
    K D^\tt = D^\tt KAD^\tt - B^\tt K^\tt C D^\tt
    \overset{\mathrm{(iii)_0}}{=}D^\tt K(\mathrm{I}_m+B C^\tt) -
    B^\tt K^\tt C D^\tt \\
    \overset{\mathrm{(ii)}}{=} D^\tt K +D^\tt K B C^\tt - B^\tt K^\tt
    D C^\tt \overset{\mathrm{(i)'}}{=} D^\tt K.
  \end{multline*}

  Recalling the definition of $\Sigma$ in~\eqref{equ:def-Sigma} and
  rewriting the relations~\eqref{equ:4-cons}, we get
  \[
    AK=KA, \qquad K(B\sigma)=(B\sigma)K, \qquad
    (\sigma^{-1}C)K=K(\sigma^{-1}C), \qquad (\sigma^{-1}D\sigma) K= K
    (\sigma^{-1}D \sigma).
  \]
  By Remark~\ref{remark.regular.matrices}, this implies that
  $A, B\sigma, \sigma^{-1}C, \sigma^{-1}D\sigma \in k[K]$, that is,
  \[
    \Sigma^{-1} X \Sigma = \left(\begin{matrix} A &
        B\sigma \\ \sigma^{-1}C &
        \sigma^{-1}D\sigma \end{matrix}\right) \in \GL_2(k[K]).
  \]
  From $B\sigma, \sigma^{-1}C\in k[K]$ and the symmetry of $\sigma$,
  one readily obtains that $B, C$ are symmetric.
  From~$(\mathrm{iii})_0$ and $\sigma^{-1} D \sigma = D^\tt$ we obtain
  that $\Sigma^{-1} X \Sigma \in \SL_2(k[K])$.

  Conversely, suppose that
  \[
    \left(\begin{matrix} A & B\sigma \\ \sigma^{-1}C &
        \sigma^{-1}D\sigma \end{matrix}\right) = \Sigma^{-1} X \Sigma
    \in \SL_2(k[K]).
  \]
  It suffices to check the conditions
  $\mathrm{(i)}$--$\mathrm{(iv)}_0$ in
  Proposition~\ref{prop:characterisation.of.G.and.G0}~(2).  This can
  be done by routine computations, using
  $K^\tt = \sigma K \sigma^{-1}$ and the fact that $k[K]$ is
  commutative.  For instance, from
  $\sigma^{-1} D \sigma , \sigma^{-1} C \in k[K]$ and
  $\sigma^\tt = \sigma$ we obtain
  $D = \sigma (\sigma^{-1} D \sigma) \sigma^{-1} = (\sigma^{-1} D
  \sigma)^\tt = \sigma D^\tt \sigma^{-1}$, thus
  $\sigma^{-1} D \sigma = D^\tt$, and
  $C \sigma^{-1} = \sigma (\sigma^{-1} C) \sigma^{-1} = (\sigma^{-1}
  C)^\tt = C^\tt \sigma^{-1}$, thus $C = C^\tt$.  This yields
  \[
    \mathrm{I}_m = \det \left(\begin{matrix} A & B\sigma \\
        \sigma^{-1}C & \sigma^{-1}D\sigma \end{matrix}\right) = A
    \cdot \sigma^{-1} D \sigma - B \sigma \cdot \sigma^{-1} C = A
    D^\tt - B C^\tt,
  \]
  and $\mathrm{(iii)_0}$ holds.  This concludes the justification
  of~\eqref{equ:G0-SL2}.
  
  \smallskip
  
  Finally, the block matrix $U \in \Mat_{2m,2}(k)$ in
  \eqref{equ:g-in-block-form} remains unconstrained in
  Proposition~\ref{prop:characterisation.of.G.and.G0} and therefore
  the group is isomorphic to $\SL_2(k[K]) \ltimes \Mat_{2m,2}(k)$. We
  can identify the natural $k[K]$-module $k^m$ with the standard
  $k[K]$-module $k[K]$, by mapping a cyclic generator of $k^m$ to the
  cyclic generator~$K$ of~$k[K]$.  Therefore $\Mat_{2m,2}(k)$ can be
  replaced by a direct sum of two copies of the standard
  $\SL_2(k[K])$-module $\mathsf{V}_\mathrm{st}(k[K])$.
\end{proof}


\subsection{The structure of the algebraic automorphism group $\G$ for
  $\Delta_K = t^m$}\label{subsec:aut.gp.xm}
Now we focus on the special case $\Delta_K = t^m$; that is, the case
\begin{equation} \label{equ:matrix-K-t^m}
  K=
    \begin{pmatrix}
      0 & 1 & 0 & \cdots & 0\\
      0 & 0 & 1 & \cdots & 0\\
      \vdots & \vdots & \vdots & \ddots & \vdots\\
      0 & 0 & 0 & \cdots & 1\\
      0 & 0 & 0 & \cdots & 0
    \end{pmatrix}
  \in \Mat_m(\Z).
\end{equation}
In this situation we can take
\begin{equation} \label{equ:matrix-sigma}
  \sigma =
    \begin{pmatrix}
      0 & 0  & \cdots & 1\\
      \vdots & \vdots &  \udots & \vdots\\
      0 & 1 & \cdots & 0\\
      1 & 0 & \cdots & 0
    \end{pmatrix}
  \in \GL_m(\Z), \quad \text{and}\quad \Sigma = \left(\begin{matrix}
      \mathrm{I}_m &  \\  & \sigma \end{matrix}\right) \in \GL_{2m}(\Z)
\end{equation}
in the analysis carried out in Section~\ref{sec:structure.of.aut.gp}.
We remark that this particular choice of $\sigma$ corresponds to the
longest element in the symmetric group $\mathrm{Sym}(m)$, with respect
to the standard generators.

\begin{proposition} \label{pro:short-split-sequ} Suppose that $K$ has
  characteristic polynomial $\Delta_K=t^m$. Then the natural
  restriction homomorphism~\eqref{equ:restr-hom} sets up, over $\Z$, a
  split short exact sequence
  \[
    \G_0(k) \quad \hookrightarrow \quad \G(k) \quad
    \xrightarrowdbl{\operatorname{Res}}
    \underbrace{\mathbf{B}_2(k)}_{\cong \G(k)/\G_0(k)} \le \GL_2(k),
  \]
  where $\mathbf{B}_2(k)$ is the group of invertible
  lower-triangular $2\times 2$ matrices.
\end{proposition}

\begin{proof}
  We show below that the image of $\G(k)$ in $\GL_2(k)$ under the
  restriction homomorphism
  \begin{enumerate}
  \item[(a)] contains $\mathbf{B}_2(k)$ by exhibiting an explicit
    section over~$\Z$, but
  \item[(b)] does not contain elements of the form
    $\left(\begin{smallmatrix} 1 & b \\ 0 &
        1 \end{smallmatrix}\right)$ with $b \neq 0$. 
  \end{enumerate}

  From this it follows that the image is precisely $\mathbf{B}_2(k)$,
  because, once we replace $k$ by its field of fractions, there are no
  properly intermediate subgroups between $\mathbf{B}_2(k)$ and
  $\GL_2(k)$, as can be seen from the Bruhat decomposition.

  \smallskip
  
  To prove (a), we define for $a,d \in k^\times$ and $c \in k$ the
  following elements of $\GL_{2m+2}(k)$:
  \begin{equation} \label{equ:U-V-W}
    \begin{split}
      U(a) & =\operatorname{diag}\left(a,a^2,\ldots,a^m, \,
        1,a^{-1},\ldots,a^{-m+1}, \, a,1 \right),   \\
      V(d)&=\operatorname{diag}\left(d^{-1},d^{-2},\ldots,d^{-m}, \,
        d,d^{2},\ldots,d^{m}, \, 1,d \right), \\
      W(c)&=\operatorname{diag}\left(\exp(c E_m), \exp(c E_m^\vee) ,
        \begin{pmatrix} 1 & 0 \\ c &1 \end{pmatrix} \right),
    \end{split}
  \end{equation}
  where $\exp(t) = \sum_{n=0}^\infty t^n/(n!)$ denotes the exponential
  series (which, evaluated on nilpotent $m \times m$-matrices, can be
  truncated after the $m$th term and thus produces finite sums) and
    \[
      E_m =
      \begin{pmatrix}
        0 & 1 & 0  & \cdots & 0\\
        0&    0 & 2      &    &     0\\
        0&    0      & 0 & \ddots &   \\
        \vdots  &    \vdots & \ddots & \ddots & (m-1)\\
        0&    0 & \cdots & 0 & 0 \\
      \end{pmatrix}, \qquad E_m^\vee =
      \begin{pmatrix}
        0 & 0  \\
        0 & -{E_{m-1}}^\tt
      \end{pmatrix}
      =
      \begin{pmatrix}
        0 & 0  & 0 &\cdots & 0 & 0\\
        0 & 0 &  0& \ddots &  \vdots & \vdots \\
        0 & -1 & 0 &\ddots &   0 & 0 \\
        0 & 0 & -2 & \ddots &   0 & 0 \\
        \vdots  & \vdots &\ddots & \ddots & 0 & 0 \\
        0& 0& \cdots & 0 & -(m-2) & 0
      \end{pmatrix}.
    \]
    A direct calculation reveals that the elements $U(a),V(d)$ and
    $W(c)$ satisfy (iii) and (iv) of
    Proposition~\ref{prop:characterisation.of.G.and.G0}, while (i) and
    (ii) hold trivially; thus $U(a), V(d), W(c) \in \G(k)$.  Moreover,
    there is an affine subgroup and $\Z$-subscheme
    $\mathbf{B} \le \mathbf{G}$ such that
    \[
      \mathbf{B}(k) = \{ U(a)V(d)W(c) \mid a,d \in k^\times \text{ and
      } c \in k \} \qquad \text{and} \qquad \mathbf{B}(k) \cong \mathbf{B}_2(k)  
    \]
    via the natural restriction homomorphism, which satisfies
    \[
      U(a) \mapsto \begin{pmatrix} a & 0 \\ 0 & 1 \end{pmatrix}, \quad
      V(d) \mapsto \begin{pmatrix} 1 & 0 \\ 0 &
        d \end{pmatrix}, \quad W(c) \mapsto \begin{pmatrix} 1 & 0 \\
        c & 1 \end{pmatrix};
    \]
    the inverse can be built from the morphisms $a \mapsto U(a)$,
    $d \mapsto V(d)$ and $c \mapsto W(c)$ which are defined over $\Z$.
    The latter is clear for $U(\cdot)$ and $V(\cdot)$, and requires a
    routine calculation for $W(\cdot)$: by induction, one sees
      that the factorials in the denominators coming from the
      exponential series duly cancel out with the entries of the relevant
      powers of $c E_m$ and $cE_m^\vee$.

    \smallskip
  
    To prove (b) we observe that $\mathrm{(iii)_0'}$ in the proof of
    Theorem~\ref{explicit.G0} holds also for elements $g \in \G(k)$ of
    the form ~\eqref{equ:g-in-block-form} which satisfy
    $Y = \left(\begin{smallmatrix} 1 & b \\ 0 &
        1 \end{smallmatrix}\right)$.  Taking the trace in equation
    (iv) of Proposition~\ref{prop:characterisation.of.G.and.G0}~(2),
    we obtain
    \begin{align*}
      mb+\tr(K)=\tr(b\mathrm{I}_m+ K) %
      & =\tr(AKD^\tt-BK^\tt C^\tt) && \text{by taking the trace in $\mathrm{(iv)}$}\\
      & =\tr(KD^\tt A-C^\tt BK^\tt ) && \text{by permuting matrices}\\
      & =\tr(KD^\tt A-K B^\tt C) && \text{by transposing the second matrix} \\
      & =\tr(K ) && \text{by applying $\mathrm{(iii)_0'}$},
    \end{align*}
    and this implies $b=0$.

    We remark that, alternatively, one can prove (b) as follows.
    Every $g \in \G(k)$ restricts to an automorphism of the centre
    ${_k Z}$ of ${_k L}$, which is represented by
    $Y = \left(\begin{smallmatrix} a & b \\ c &
        d \end{smallmatrix}\right) \in \GL_2(k)$ with respect to the
    chosen basis $z_1, z_2$, and similarly for $g^{-1}$.  The image
    $(z_1',z_2')$ of the pair $(z_1,z_2)$ under $g^{-1}$ yields two
    antisymmetric bilinear forms which encode the Lie bracket;
    inspection of the form associated to $z_2'$ shows that
    $b \mathrm{I}_m + d K$ should have the same rank as~$K$, namely
    $m-1$; thus $b = 0$.
\end{proof}

\begin{proof}[Proof of Theorem~\ref{thm:structure.of.Aut}]
  In view of Theorem~\ref{explicit.G0} and
  Proposition~\ref{pro:short-split-sequ}, it only remains to show that
  the algebraic group $\G$ is connected.  As
  \[
    \G(k) \cong \mathbf{B}_2(k) \ltimes \left(\SL_2(k[K]) \ltimes
      \mathsf{V}_\mathrm{st}(k[K])^{\oplus 2}\right)
  \]
  by an isomorphism of group schemes over $\Z$, the connectedness of
  $\G$ follows from the fact that $\G$ is generated by one-parameter
  subgroups, which are, in particular, affine irreducible varieties
  containing~$1$; for instance, see \cite[Prop.~1.16]{MaTe11}.
\end{proof}

For our next step we record also the following consequence of
Theorem~\ref{explicit.G0} and Proposition~\ref{pro:short-split-sequ}.

\begin{corollary}\label{Aut.after.Sigma}
  Suppose that $K$ has characteristic polynomial $\Delta_K =t^m$.
  Then the group $\G_0(k)$ is conjugate to the subgroup of
  $\GL_{2m+2}(k)$ consisting of elements of the form
  \[
    \begin{pmatrix} A & B & E \\ C & D & F \\ 0 & 0 &
      \mathrm{I_2}\end{pmatrix},
  \]
  where $A,B,C, D \in \Mat_{m}(k)$ satisfy $AD-BC=\mathrm{I}_m$ and
  are in Toeplitz form, that is,
  \begin{equation}\label{T.form}
    A=\begin{pmatrix} a_1 & a_2 & \cdots & a_{m}  \\
      & a_1 & \ddots & \vdots \\
      & & \ddots & a_2  \\
      & &  & a_1 
    \end{pmatrix}, \;
    B=\begin{pmatrix} b_1 & b_2 & \cdots & b_{m}  \\
      & b_1 & \ddots & \vdots \\
      & & \ddots & b_2  \\
      & & & b_1
    \end{pmatrix}, \;
    C=\begin{pmatrix} c_1 & c_2 & \cdots & c_{m}  \\
      & c_1 & \ddots & \vdots \\
      & & \ddots & c_2  \\
      & & & c_1
    \end{pmatrix}, \;
    D=\begin{pmatrix} d_1 & d_2 & \cdots & d_{m}  \\
      & d_1 & \ddots & \vdots \\
      & & \ddots & d_2  \\
      & & & d_1
    \end{pmatrix}
  \end{equation}
  with suitable entries $a_1, \ldots, d_m \in k$ and entries $0$ in
  white spaces, and $E,F \in \Mat_{m,2}(k)$.
  
  The group $\G(k)$ is generated by $\G_0(k)$ and the elements $U(a)$,
  $V(d)$ and $W(c)$, for $a,d \in k^\times$ and $c \in k$, which are
  defined in the proof of Proposition~\ref{pro:short-split-sequ}.
\end{corollary}


\subsection{Change of coordinates}

For $\Delta_K = t^m$, the Lie lattice $L$ is intimately linked to the
nilpotent group $\Gamma_{t^m}$, defined in~\eqref{equ:def-Gamma-tm},
and the algebraic automorphism group $\G$ plays a central role in the
treatment of the pro-isomorphic zeta function of~$\Gamma_{t^m}$; see
Section~\ref{section:background}.  With a view towards the computation
of the pro-isomorphic zeta function of the group $\Gamma_{t^m}$, we
perform a change of basis
\[
  \text{from} \quad \mathcal{S} = (x_1,x_2, \dots,x_m,\, y_1, y_2,
  \dots, y_m,\, z_1,z_2) \quad \text{to} \quad \mathcal{S}^* =
  (x_1, y_m,\, x_2, y_{m-1},\, \dots,\, x_m, y_1, \, z_2,z_1).
\]
This basis change is achieved by conjugating first with
$\operatorname{diag}(\Sigma,\mathrm{I}_2)$, already built into
Corollary~\ref{Aut.after.Sigma} and reversing the order of
$y_1, \ldots, y_m$, and then with
$\operatorname{diag} \left(\Theta,\left( \begin{smallmatrix} 0&1 \\
      1&0 \end{smallmatrix} \right) \right)$, where $\Theta$
corresponds to the permutation of $\{1,2,\ldots, 2m\}$ given by
\begin{equation}\label{equ:change-basis-2}
  \begin{cases}
    i \mapsto 2i-1& \text{if $1 \le i \le m$,} \\
    i \mapsto 2(i-m)  & \text{if $m < i \le 2m$.} \\
  \end{cases}
\end{equation}
From the results in Section~\ref{subsec:aut.gp.xm} we obtain the
following description of $\G(k)$, with respect to the
basis~$\mathcal{S}^*$.

\begin{proposition}\label{prop:Grp-new-coordinates}
  Suppose that $K$ has characteristic polynomial $\Delta_K=t^m$.
  Then, with respect to the basis~$\mathcal{S}^*$, the elements
  of~$\G_0(k)$ take the form
  \begin{equation}\label{explicit.G0.conjugated}
    \left(
      \begin{array}{ccccc|cc}
		X_1 & X_2 & X_3& \cdots & X_m& *& *\\
		& X_1 &\ddots &\ddots&\vdots &\vdots & \vdots\\
		&& \ddots &\ddots &X_3 &*& *\\
		&&& X_1 &X_2 &* & *\\
		&&&& X_1 &*&*\\
		\hline
		&&&&& 1& 0\\
		&&&&& & 1
      \end{array}
    \right), \quad
    \begin{array}{r} \text{with
      $X_i = \begin{pmatrix} a_i & b_i \\ c_i &
          d_i \end{pmatrix} \in \Mat_2(k)$ for $1 \le i \le m$ and} \\
      \text{arbitrary entries in the positions
      marked $*$,}
    \end{array}
  \end{equation}
  such that the matrices $A,B,C,D$ defined as in~\eqref{T.form}
  satisfy $AD-BC = \mathrm{I}_m$.

  \smallskip

  Furthermore, still with respect to the basis~$\mathcal{S}^*$, the
  group $\G(k)$ is generated by $\G_0(k)$ and 
  \begin{equation} \label{equ:U'-V'-W'}
    \begin{split}
      U'(a) & = T^{-1} \,\big( U(a) V(a) \big)\, T =
      \operatorname{diag}\left(\begin{pmatrix} 1 & \\ &
          a \end{pmatrix}, \ldots, \begin{pmatrix} 1 & \\ &
          a \end{pmatrix}, \begin{pmatrix} a & \\ &
          a \end{pmatrix} \right),   \\
      V'(d) & = T^{-1} \,\big( U(d)^{m-1} V(d)^{m} R(d)^m \big)\, T
      =\operatorname{diag} \left( d^{m-1}\mathrm{I_2},
        d^{m-2}\mathrm{I_2}, \ldots , d\mathrm{I_2},
        \mathrm{I_2}, \begin{pmatrix} d^{m} & \\ &
          d^{m-1} \end{pmatrix} \right), \\
      W'(c) & = T^{-1} \, W(c) \, T      
    \end{split}
  \end{equation}
  for
  $a,d \in k^\times$ and $c \in k$, where $T = \operatorname{diag}
  \left(\Sigma\Theta,\left( \begin{smallmatrix} 0&1 \\ 
        1&0 \end{smallmatrix} \right) \right)$, the one-parameter
  groups $U(\cdot), V(\cdot), W(\cdot)$ are as in~\eqref{equ:U-V-W}
  and 
  $R(d) = \operatorname{diag}(d,d,\ldots,d,d^{-1},d^{-1},\ldots,
  d^{-1}, 1,1) \in \G_0(k)$.
\end{proposition}

\begin{corollary}\label{cor:corollary-reductive-part} Suppose that $K$
  has characteristic polynomial $\Delta_K=t^m$.  Then the quotient of
  $\G$ by its unipotent radical $\mathbf{N}$ is isomorphic to
  $\GL_2 \times \GL_1$, with an explicit section defined over $\Z$
  with respect to the basis~$\mathcal{S}^*$ as follows:
  \[
    \GL_2(k) \times \GL_1(k) \to \mathbf{H}(k), \quad (A,\nu) \mapsto
    \left(
      \begin{array}{ccccc|cc}
        \nu^{m-1}A&&&&0&&\\
                   & \nu^{m-2}A&&&&&\\
                   && \ddots &&& &\\
                   &&& \nu A &&&\\
        0&&&& A &&\\
        \hline
                   &&&&& \nu^m \det A& 0\\
                   &&&&& 0& \nu^{m-1}\det A
      \end{array}
    \right).
  \]
\end{corollary}

\begin{proof}
  We consider the affine subgroup $\mathbf{N}$ of $\G$ such that
  $\mathbf{N}(k)$ is generated by elements of the form
  \eqref{explicit.G0.conjugated} with $X_1=\mathrm{I}_2$ together with
  elements of the subgroup $\{ W'(c) \mid c \in k\}$:
  see~\eqref{equ:U'-V'-W'}.  The group $\mathbf{N}$ is a connected
  unipotent normal subgroup of~$\G$.

  Moreover, the quotient $\G(k) / \mathbf{N}(k)$ is generated by the
  block-diagonal matrices
  $\operatorname{diag}(X_1, \dots, X_1, \mathrm{I}_2)$ with
  $X_1 \in \SL_2(k)$ and by the one-parameter subgroups
  $\{ U'(a) \mid a \in k^\times \}$ and
  $\{ V'(\nu) \mid \nu \in k^\times \}$; this analysis also provides a
  section for $\G \to \G/\mathbf{N}$ over~$\Z$.  Finally
  $\G/\mathbf{N} \cong \GL_2 \times \GL_1$ is reductive, and thus
  $\mathbf{N}$ is the unipotent radical of~$\G$.
\end{proof}

\begin{remark}\label{U'V'} For computational purposes we replaced the
  generators $U(\cdot)$ and $V(\cdot)$ by the generators $U'(\cdot)$
  and $V'(\cdot)$.  They generate the same torus, modulo $\G_0$ and up
  to coordinate change; see~\eqref{equ:U'-V'-W'}.

  For similar reasons, a further simplification of the computation of
  the pro-isomorphic zeta function can be achieved by replacing the
  one-parameter subgroup $W'(\cdot)$ in \eqref{equ:U'-V'-W'} by
  \[
    c \mapsto W''(c) = T^{-1} \, \operatorname{diag}\left(\exp \left(
        cE_m+\hlf(1-m)cK \right), \exp \left(c E_m^\vee-\hlf(1-m) c
        K^\tt \right),
      \begin{pmatrix} 1 & 0 \\ c &1 \end{pmatrix} \right) \, T.
  \]  
  This switch is inspired by Lie algebra considerations and works for
  an arbitrary $\Z$-algebra $k$ if $m$ is odd; for even $m$ the switch
  requires that $k$ is a $\Z[\hlf]$-algebra.  For the applications in
  the present paper, we do not need the variant for $m=2$ and we use
  it only for $m=3$.  Hence no primes need to be excluded when we
  compute the local pro-isomorphic zeta functions for
  Theorems~\ref{x2example} and \ref{x3example}.
\end{remark}

\begin{example} \label{exa:unipotent-radical-m23} In order to compute
  later on the pro-isomorphic zeta functions of the groups
  $\Gamma_{t^m}$ for $m\in \{2,3\}$, we record in these cases explicit
  descriptions of the unipotent radical~$\mathbf{N}$ of $\G$, with
  respect to the basis~$\mathcal{S}^*$.  For completeness we also
  provide a description for $m=1$ which is straightforward; compare
  with~\cite[\S 3.3.4]{Be05}.  We have
  \begin{align*}
    \mathbf{N}(k) %
    & = \left\{\left(
      \begin{array}{c|cc}
        \mathrm{I}_2 & *&*\\
        \hline
                      &1&\lambda\\
                     & &1
      \end{array}
                          \right)
                          \,\Bigg|\,
                          \begin{array}{l}
                            \lambda \in k,\,  \text{and arbitrary entries}
                            \\ \text{in the positions marked $*$}
                          \end{array}
    \right\}  && \text{if $m=1$,} \\    
    \mathbf{N}(k) %
    & = \left\{\left(
      \begin{array}{cc|cc}
        \mathrm{I}_2 & X_2  & *&*\\
        0 & \mathrm{I}_2 &*&*\\
        \hline
                     & &1&\tr(X_2)\\
                     && &1
      \end{array}
                          \right)
                          \,\Bigg|\,
                          \begin{array}{l}
                            X_2 \in \Mat_2(k),\,  \text{and arbitrary}
                            \\ \text{entries in the
                            positions marked $*$}
                          \end{array}
    \right\}  && \text{if $m=2$,} \\
    \mathbf{N}(k) %
    &= \left\{\left(
      \begin{array}{ccc|cc}
        \mathrm{I}_2 & X_2&X_3 &*&*\\
        0& \mathrm{I}_2  & X_2+\lambda  \mathrm{I}_2 & * & *\\
        0&0& \mathrm{I}_2 & * & *\\
        \hline
                     & & &1&\lambda\\
                     & && &1
      \end{array}
                            \right)
                            \,\Bigg|\,
                            \begin{array}{l}
                            	X_2, X_3\in \Mat_2(k) \text{ with } \tr(X_2)=0,\\
                               \tr(X_3)+\det(X_2)=0,\,
                              \lambda \in k, \text{ and}\ \\
                              \text{arbitrary entries in the 
                              positions}\\ \text{marked $*$}
                            \end{array}
    \right\} && \text{if $m=3$.}
  \end{align*}

  \smallskip

  Indeed, for $m=2$ we substitute $X_1=\mathrm{I}_2$ in
  \eqref{explicit.G0.conjugated} and use the determinant equation in
  Proposition~\ref{prop:Grp-new-coordinates} to obtain
  \[
    \begin{pmatrix} a_1 & a_2 \\ & a_1 \end{pmatrix} \begin{pmatrix}
      d_1 & d_2 \\ & d_1 \end{pmatrix} -\begin{pmatrix} b_1 & b_2 \\ &
      b_1 \end{pmatrix} \begin{pmatrix} c_1 & c_2 \\ &
      c_1 \end{pmatrix} = \begin{pmatrix} 1 & \\ & 1 \end{pmatrix}
  \]
  with $a_1-1=d_1-1=b_1=c_1=0$, and therefore $a_2+d_2=0$, namely
  $\tr(X_2)=0$; this accounts for the contribution of~$\G_0(k)$.  The
  explicit form of $W'(c)$ for $c \in k$ is
  \[
    W'(c) = \left(
      \begin{array}{cc|c}
        \begin{matrix} 1 & 0 \\ 0 & 1 \end{matrix}  %
                                  & \begin{matrix} c
                                    & 0 \\ 0 & 0 \end{matrix}  &   \begin{matrix} 0 & 0 \\ 0 &
                                    0 \end{matrix} \\ 
                         &   \begin{matrix} 1 & 0 \\ 0 &
                           1 \end{matrix}  &   \begin{matrix} 0 & 0 \\
                           0 & 0 \end{matrix} \\ 
        \hline
                         & &  \begin{matrix} 1 & c \\ 0 & 1 \end{matrix} 
      \end{array}
    \right).
  \]
  Combining the contributions, the result for $m=2$ follows.

  \smallskip

  For $m=3$, we start again by substituting $X_1=\mathrm{I}_2$ in
  \eqref{explicit.G0.conjugated}, which together with the determinant
  equation
  \[
    \begin{pmatrix} a_1 & a_2 & a_3 \\ & a_1 & a_2 \\ & &
      a_3 \end{pmatrix} \begin{pmatrix} d_1 & d_2 & d_3 \\ & d_1 & d_2
      \\ & & d_3 \end{pmatrix} -\begin{pmatrix} b_1 & b_2 & b_3 \\ &
      b_1 & b_2 \\ & & b_3 \end{pmatrix} \begin{pmatrix} c_1 & c_2 &
      c_3 \\ & c_1 & c_2 \\ & & c_3 \end{pmatrix} = \begin{pmatrix} 1
      & & \\ & 1& \\ & & 1 \end{pmatrix}
  \]
  gives $a_1-1=d_1-1=b_1=c_1=0$ and
  \[
    \begin{split}
      \tr(X_2)&=a_2+b_2=0, \\
      \tr(X_3)+\det(X_3)&=a_3+d_3+a_2d_2-b_2c_2=0;
    \end{split}
  \]
  this yields the intersection of the unipotent radical
  $\mathbf{N}(k)$ with $\G_0(k)$.

  Offsetting $W'(\cdot)$ in accordance with Remark~\ref{U'V'}, we get the
  one-parameter subgroup $W''(\cdot)$ which takes the form
  \[
    W''(c)= \left(
      \begin{array}{ccc|cc}
        \mathrm{I}_2 & 0 &0 & &\\
        0& \mathrm{I}_2  & c  \mathrm{I}_2 &  & \\
        0&0& \mathrm{I}_2 &  & \\
        \hline
                     & & &1&c\\
                     & && &1
      \end{array}
    \right), \qquad \text{for $c \in k$.}
  \]
  Combining the contributions, we arrive at the result for $m=3$.
\end{example}


\section{Machinery for computing $p$-adic integrals over algebraic
  groups} \label{section:background}

In this section we collect various facts and notation in order to use
the technology developed in~\cite{GrSeSm88, Ig89, dSLu96, Be11}.  The
general treatment produces a finite, but typically unspecified set of
`exceptional' primes; we take care to verify that, for the
applications in this paper, there is no need to exclude any primes.

\subsection{Lie correspondence for class-two nilpotent groups}\label{section:class-2-correspondence}

Let $\Gamma$ be a finitely generated torsion-free nilpotent group.
Grunewald, Segal and Smith~\cite[Thm.~4.1]{GrSeSm88} showed that the
local pro-isomorphic zeta functions of~$\Gamma$ are closely linked to
the local pro-isomorphic zeta functions of a nilpotent $\Z$-Lie
lattice~$L$ that can be constructed from~$\Gamma$; indeed,
$\zeta_{\Gamma,p}^\wedge(s) = \zeta_{L,p}^\wedge(s)$ for almost all
primes~$p$.  Furthermore, they remark that, if $\Gamma$ has nilpotency
class two, a suitable Lie correspondence can be implemented more
directly, and they highlight consequences for other types of zeta
functions.  The direct correspondence has been reinterpreted and put
to use, for instance, in~\cite[\S 2.4.1]{StVo14}.  For the record, we
state and explain the implications of the special construction in
nilpotency class two for pro-isomorphic zeta functions, where it is
applied not merely to a group, but also to its lattice of subgroups;
compare with~\cite[Rem.~2.2]{BeKlOn18}.

Let $\Gamma$ be a finitely generated torsion-free class-two nilpotent
group of Hirsch length~$d$, and let $Z = \mathrm{Z}(\Gamma)$ denote
its centre.  Then the isomorphism type of $\Gamma$ is uniquely
determined by
$\Gamma/Z = \langle x_1Z, \ldots, x_aZ \rangle \cong \Z^a$,
$Z = \langle y_1, \ldots, y_{d-a} \rangle \cong \Z^{d-a}$ and the map
$\gamma \colon \Gamma/Z \times \Gamma/Z \to Z$,
$(gZ,hZ) \mapsto [g,h]$.  In fact, this data yields a $\Z$-Lie lattice
\begin{equation} \label{equ:Lie-lattice-in-class-2}
  L = \Z \dot x_1 \oplus \ldots \oplus \Z \dot x_a
  \oplus \Z \dot y_1 \oplus \ldots \oplus \Z
  \dot y_{d-a} \cong \Gamma/Z \oplus Z,
\end{equation}
where the Lie bracket is induced by the anti-symmetric bi-additive map
$\gamma$ and the stipulation that
$\Z \dot y_1 \oplus \ldots \oplus \Z \dot y_{d-a}$ be central in~$L$:
\begin{align*}
  & [\dot x_i, \dot x_j]_\mathrm{Lie} = \sum_{k=1}^{d-a} c_{i,j,k} \, \dot
  y_k %
  & & \text{for $1 \le i,j \le a$, where $\gamma(x_iZ,x_jZ) = [x_i,x_j] =
    \prod_{k=1}^{d-a} y_k^{\, c_{i,j,k}}$,} \\
  & [\dot x_i, \dot y_j]_\mathrm{Lie} = [\dot y_j, \dot
    y_k]_\mathrm{Lie} = 0 %
  & & \text{for $1 \le i \le a$ and $1 \le j \le k \le d-a$.}
\end{align*}
Conversely, given such a Lie lattice one can define a class-two
nilpotent group, essentially by factoring out from the free class-two
nilpotent group on $d$ generators
$\hat x_1, \ldots, \hat x_a, \hat y_1, \ldots, \hat y_{d-a}$ the
relations
\begin{align*}
  & [\hat x_i, \hat x_j]_\mathrm{Lie} = \prod_{k=1}^{d-a} \hat y_k^{\,
    c_{i,j,k}} %
  & & \text{for $1 \le i,j \le a$, where $[\dot x_i, \dot x_j]_\mathrm{Lie} =
    \sum_{k=1}^{d-a} c_{i,j,k}\, y_k$,} \\
  & [\hat x_i, \hat y_j] =  [\hat y_j, \hat y_k] = 1 %
  & & \text{for $1 \le i \le a$ and $1 \le j \le k \le d-a$.}
\end{align*}
Moreover, the two constructions set up a $1$-to-$1$ correspondence, up
to isomorphism, between finitely generated torsion-free class-two
nilpotent groups of Hirsch length~$d$ and class-two nilpotent $\Z$-Lie
lattices of dimension~$d$.  For short, we call this the
  \emph{class-two Lie correspondence}.

We observe that, for any prime~$p$, essentially the same constructions
yield a `local' class-two Lie correspondence, up to isomorphism,
between torsion-free class-two nilpotent pro-$p$ groups of rank $d$
and class-two nilpotent $\Z_p$-Lie lattices of dimension~$d$;
compare with~\cite[\S 1]{GrSe84} and \cite[\S 2.4.1]{StVo14}.

\begin{proposition} \label{pro:class-2-correspondence} Let $\Gamma$ be
  a finitely generated torsion-free class-two nilpotent group of
  Hirsch length~$d$, with centre $Z = \mathrm{Z}(\Gamma)$, such that
  $\Gamma/Z = \langle x_1Z, \ldots, x_aZ \rangle \cong \Z^a$ and
  $Z = \langle y_1, \ldots, y_{d-a} \rangle \cong \Z^{d-a}$.  Let $L$
  be the $\Z$-Lie lattice associated to $\Gamma$ under the class-two
  Lie correspondence as in~\eqref{equ:Lie-lattice-in-class-2}.  Then
  the map
  \begin{equation} \label{equ:corresp-wrt-basis}
    \Gamma \to L, \quad \prod_{i=1}^a x_i^{\, m_i} \prod_{j=1}^{d-a}
    y_j^{\, n_j} \mapsto \sum_{i=1}^a m_i \, \dot x_i + \sum_{j=1}^{d-a}
    n_j \, \dot y_j
  \end{equation}
  induces an index-preserving $1$-to-$1$ correspondence between
  finite-index subgroups $\Delta \le \Gamma$ and finite-index Lie
  sublattices $M \le L$.

  Furthermore, subgroups $\Delta$ satisfying
  $\widehat{\Delta} \cong \widehat{\Gamma}$ are bijectively matched
  with Lie sublattices $M$ such that the $\Z_p$-Lie lattices
  $\Z_p \otimes_\Z M$ and $L_p = \Z_p \otimes_\Z L$ are isomorphic for
  all primes~$p$.  In particular, this implies that
  \[
    \zeta_{\Gamma,p}^\wedge(s) = \zeta_{L,p}^\wedge(s) =
    \zeta_{L_p}^\mathrm{iso}(s) \qquad \text{for all primes~$p$.}
    \]
\end{proposition}

\begin{proof}
  It is elementary to check that the (non-canonical)
  map~\eqref{equ:corresp-wrt-basis} sets up an index-preserving
  $1$-to-$1$ correspondence between finite-index subgroups of $\Gamma$
  and finite-index Lie sublattices $L$; this was already remarked
  in~\cite{GrSeSm88}, just after the proof of Theorem~4.1 in that
  paper.

  We fix a prime~$p$, a finite-index subgroup $\Delta \le \Gamma$ and
  its image $M \le L$ under the map~\eqref{equ:corresp-wrt-basis}.  It
  remains to justify that
  $\widehat{\Delta}_p \cong \widehat{\Gamma}_p$ if and only if
  $\Z_p \otimes M \cong L_p$.  First we observe that
  $\mathrm{C}_\Gamma(\Delta) = Z$ and thus
  $\Delta \cap Z = \mathrm{Z}(\Delta)$.  This implies that $M$ is
  isomorphic to the $\Z$-Lie lattice associated canonically to
  $\Delta$ via the class-two Lie correspondence.  Since the
  constructions that lead to the class-two Lie correspondences for
  discrete nilpotent groups and for nilpotent pro-$p$ groups are
  essentially the same, we see that the $\Z_p$-Lie lattice associated
  canonically to the pro-$p$ completion $\widehat{\Delta}_p$ can be
  obtained from $M$ by extension of scalars, i.e., it is isomorphic to
  $\Z_p \otimes_\Z M$.  The same analysis applies, of course, also to
  $\Gamma$ in place of~$\Delta$.  Applying the local class-two Lie
  correspondence, we deduce that
  $\widehat{\Delta}_p \cong \widehat{\Gamma}_p$ if and only if
  $\Z_p \otimes M \cong L_p$.
\end{proof}

\begin{remark}
  The map \eqref{equ:corresp-wrt-basis} used in
  Proposition~\ref{pro:class-2-correspondence} depends on the implicit
  choice of coset representatives $x_1, \ldots, x_a$ for $\Gamma$
  modulo~$Z$.  However, if $\Delta \le \Gamma$ is a finite-index
  subgroup, then it admits a generating $d$-tuple of the form
  \begin{multline*}
    x_1^{\, e_{11}} x_2^{\, e_{12}} \cdots x_a^{\, e_{1a}} y_1^{\,
      f_{11}} y_2^{\, f_{12}} \cdots y_{d-a}^{\, f_{1(d-a)}}, \quad
    x_2^{\, e_{22}} \cdots x_a^{\, e_{2a}} y_1^{\, f_{21}} y_2^{\,
      f_{22}} \cdots y_{d-a}^{\, f_{2(d-a)}}, \quad \ldots, \quad
    x_a^{\, e_{aa}} y_1^{\,
      f_{a1}} y_2^{\, f_{a2}} \cdots y_{d-a}^{\, f_{a(d-a)}}, \\
    y_1^{\, g_{11}} y_2^{\, g_{12}} \cdots y_{d-a}^{\, g_{1(d-a)}},
    \quad y_2^{\, g_{22}} \cdots y_{d-a}^{\, g_{2(d-a)}}, \quad
    \ldots, \quad y_{d-a}^{\, g_{(d-a)(d-a)}}
  \end{multline*}
  with integer exponents.  Moreover, $\Delta$ is isomorphic to the
  subgroup $\Delta_1 \le \Gamma$ generated by
  \begin{equation*}
    x_1^{\, e_{11}} \cdots x_a^{\, e_{1a}}, \quad
    x_2^{\, e_{22}} \cdots x_a^{\, e_{2a}}, \quad \ldots, \quad
    x_a^{\, e_{aa}},
    y_1^{\, g_{11}}  \cdots y_{d-a}^{\, g_{1(d-a)}},
    \quad y_2^{\, g_{22}} \cdots y_{d-a}^{\, g_{2(d-a)}}, \quad
    \ldots, \quad y_{d-a}^{\, g_{(d-a)(d-a)}}.
  \end{equation*}
  Similarly, the image $M$ of $\Delta$
  under~\eqref{equ:corresp-wrt-basis} is isomorphic to the image $M_1$
  of $\Delta_1$ under~\eqref{equ:corresp-wrt-basis}, which has $\Z$-basis
  \[
    e_{11} \, \dot x_1 + \ldots + e_{1a} \, \dot x_a, \quad \ldots, \quad
    e_{aa} \, \dot x_a, \quad g_{11}\, \dot y_1 + \ldots + g_{1(d-a)} \,
    \dot y_{d-a}, \quad  \ldots, \quad g_{(d-a)(d-a)} \, \dot y_{d-a}.
  \]
  In this way, we see that there is a canonical map from finite-index
  subgroups of~$\Gamma$ to finite-index \emph{graded} Lie sublattices
  of~$L$, with finite fibers, where $L$ is regarded as a graded
  $\Z$-Lie lattice with respect to the decomposition
  $L = L_{(1)} \oplus L_{(2)}$ with $L_{(1)} = \Gamma/Z$ and
  $L_{(2)} = Z$.  For any graded Lie sublattice
  $M = M_{(1)} \oplus M_{(2)} \le L$, the fiber above $M$ has size
  $\lvert L_{(2)} : M_{(2)} \rvert^a$.
\end{remark}

\subsection{Local pro-isomorphic zeta functions as integrals over reductive groups}\label{section:reduction-integral}

Recall from Section~\ref{section:structure-of-Aut} the notion of the
algebraic automorphism group $\mathbf{Aut}(L)$ of a $\Z$-Lie
lattice~$L$; via a $\Z$-basis of $L$, the group $\mathbf{Aut}(L)$ is
realised as an affine $\Z$-group scheme $\G \le \GL_d$, where $d$ is
the $\Z$-rank of~$L$.  As before, for any commutative ring $R$
with~$1$ we write ${_R L} = R \otimes_\Z L$ and, for short,
we set
\[
  L_p = {_{\Z_p} L} \qquad \text{for every prime~$p$.}
\]

\begin{proposition}[{Grunewald, Segal,
    Smith~\cite[Prop.~3.4]{GrSeSm88}}] \label{pro:integral-formula}
  Let $L$ be a nilpotent $\Z$-Lie lattice of $\Z$-rank~$d$, and let
  $\G = \mathbf{Aut}(L) \le \GL_d$ denote the algebraic automorphism
  group of~$L$ with respect to some $\Z$-basis.  For each prime $p$,
  let
  \[
    G_p = \G(\Q_p) \qquad \text{and} \qquad G^+_p = G_p \cap
    \Mat_d(\Z_p) \cong \Aut({_{\Q_p} L}) \cap \End({_{\Z_p} L}),
  \]
  equipped with the right Haar measure $\mu_{G_p}$ on the locally
  compact group $G_p$ such that $\mu_p(\G(\Z_p)) =1$.  Then for all
  primes $p$,
  \begin{equation}\label{equ:Lp-integral}
    \zeta_{L_p}^{\mathrm{iso}}(s) = \int_{G^+_p} \lvert \det g \rvert_p^{\, s} \,
    \mathrm{d}\mu_{G_p}(g)
  \end{equation}
  where $\zeta_{L_p}^{\mathrm{iso}}(s)$ enumerates Lie sublattices that
  are isomorphic to~$L_p$.
\end{proposition}

We may decompose the $1$-component $\G^\circ$ into a semidirect
product $\G^\circ = \mathbf{N} \rtimes \mathbf{H}$ of its unipotent
radical $\mathbf{N}$ and a reductive group~$\mathbf{H}$; compare
with~\cite[\S VIII.4]{Ho81}.  Fix a prime $p$ and
write~$G = \G(\Q_p)$, $N = \mathbf{N}(\Q_p)$, $H = \mathbf{H}(\Q_p)$.
Let $V = {_{\Q_p}L} \cong \Q_p^{\, d}$ be the $\Q_p^{\, d}$-vector
space on which~$G$ acts from the right. In~\cite[\S 2]{dSLu96}, du
Sautoy and Lubotzky provide a general framework for reducing an
integral of the form~\eqref{equ:Lp-integral} to an integral over a
suitable subset $H^+ \subseteq H$.  Their reduction depends, in
general, on several technical assumptions (some of which can be
realised by excluding finitely many primes):
\begin{itemize}
\item[(a)] $\G = \G^\circ$ is connected.
\item[(b)] There exists a vector space decomposition
  $V=\bigoplus_{i=1}^c U_i$, with associated flag
  $V_j=\bigoplus_{i=j}^c U_i$, $1 \le j \le c+1$, such that each $U_i$
  is $H$-invariant, each $V_j$ is $N$-invariant and the induced action
  of $N$ on each quotient $V_j/V_{j+1}$, $1 \le j \le c$, is trivial.
\item[(c)] A certain lifting condition holds with respect to this
  decomposition; see~\cite[Assumption~2.3]{dSLu96} for a complete
  description and Condition~\ref{lifting-condition} below for a
  specific instance.
\end{itemize}
The requirement that the action of $N$ on the quotients $V_j/V_{j+1}$
be trivial is not actually needed for the reduction. However, it is
usually desirable -- both for theoretical and practical
applications. We will shortly see that in our applications we need to
drop this requirement.

We now specialise to the case where $L$ is a $D^*$-Lie lattice
associated, via~\eqref{eq:structure.of.L} above, to the polynomial
$\Delta(t) =t^m$ for some integer $m\geq 2$.  Note that $L$ is a
class-two nilpotent $\Z$-Lie lattice of rank $d=2m+2$ with rank-two
centre and $\mathrm{Z}(L) = [L,L]$.  Our aim is to identify modified
versions of the above technical assumptions in order to carry out a
reduction of the integral in the spirit of du Sautoy and Lubotzky,
without excluding any primes.  In our setting, $\G$ is connected and
the splitting $\G = \mathbf{N} \rtimes \mathbf{H}$ is very explicit;
see Corollary~\ref{cor:corollary-reductive-part}.  Thus we are not
worried about~(a).  We write $V = U_1\oplus U_2$, where
$U_2 =[_{\Q_p} L, _{\Q_p} L]$ and $U_1$ is an $H$-stable complement
to~$U_2$ in~$V$, corresponding to the abelianisation of~${_{\Q_p} L}$;
in the case of interest to us, $U_1$ is the $\Q_p$-span of a natural
set of generators for the Lie lattice~$L_p$.  Note that $U_2$ is
automatically invariant under the action of~$G$, while $U_1$ is
$H$-invariant by construction; however, our decomposition is `coarse'
in the sense that the actions of $N$ on $V/U_2$ and on $U_2$ are not
trivial as stipulated in~(b).

We now go about describing a weak version of~(c) that suffices for our
purposes.  Remarkably, \cite[Assumption~2.3]{dSLu96} does not apply to
the $D^*$-Lie lattice associated to $t^3$; compare with
Remark~\ref{remark:theta-not-character} below.  Let
$N_{1} = N \cap \ker(\psi'_2)$, where
$\psi'_2 \colon G \rightarrow \Aut(V/U_2)$ denotes the natural action.
Since $U_2$ is $N$-invariant, we may define the induced map
$\psi_2 \colon G/N_{1} \rightarrow \Aut(V/U_{2}) \le \GL_{2m}(\Q_p)$,
and the set
\[
  (G/N_{1})^+ = \psi_{2}^{-1}\big(\psi_{2}(G/N_{1})  \cap \Mat_{2m}(\Z_p) \big)  ,
\]
where $2m=\dim V/U_2$ is the dimension of the abelianisation
of~$L_p$.  

\begin{condition}\label{lifting-condition}
  For every $g_0 N_{1} \in (G/N_{1})^+$ there exists $g \in G^+$ such
  that $g_0 N_{1} = g N_{1}$.
\end{condition}

\begin{remark}\label{remark-lifting}
  The effect of Condition~\ref{lifting-condition} is weaker than that
  of \cite[Assumption~2.3]{dSLu96}, because in our situation $N$ does
  not act trivially on~$V/U_2$. Condition~\ref{lifting-condition} is
  trivially satisfied due to the freedom to replace $g_0$ by
  $g \in g_0 N_{1}$ such that $vg$ has zero component in $U_2$ for all
  $v\in U_1$. In matrix terms, this amounts to replacing the top-right
  block `above the centre' by zeros.  The action of $g_0$ and $g$ on
  $U_2$ is the same and induced by the action on $V/U_2$; as the
  action on $V/U_2$ is `integral', it is also integral on~$U_2$.
\end{remark}

Define $\theta_{0} \colon H \rightarrow \R_{\geq 0}$ by setting
\[
  \theta_{0}(h) = \mu_{N/N_1} \big( \{u N_1 \in N/N_1 \mid uh N_1 \in
  (G/N_{1})^+ \} \big),
\]  
where $\mu_{N/N_1}$ denotes the right Haar measure on~$N/N_1$,
normalised such that the set
$\psi_{2}^{-1} \big(\psi_{2}(N/N_1) \cap \Mat_{2m}(\Z_p) \big)$ has
measure~$1$.  Similarly, define
$\theta_{1} \colon H \rightarrow \R_{\geq 0}$ by setting
\[
  \theta_{1}(h) = \mu_{N_{1}} \big( \{u \in N_1 \mid nh \in G^+ \}
  \big),
\]  
where $\mu_{N_1}$ denotes the right Haar measure on~$N_1$, normalised
such that the set ${N_1^+=N_1(\Z_p)}$ has measure~$1$.

Write $\mu_G$, respectively $\mu_H$, for the right Haar measure
on~$G$, respectively $H$, normalised such that
$\mu_G \big( \G(\Z_p) \big) = 1$, respectively
$\mu_H \big( \mathbf{H}(\Z_p) \big) = 1$.  From $G = N \rtimes H$ one
deduces (using Condition~\ref{lifting-condition} and
Remark~\ref{remark-lifting}) that
$\mu_G = \mu_{N/N_1}\cdot \mu_{N_{1}} \cdot \mu_H$.  Setting
$G^+ = G \cap \Mat_{2m+2}(\Z_p)$ and $H^+ = H \cap \Mat_{2m+2}(\Z_p)$,
one obtains the following by a mild adaptation of the proof of
\cite[Thm.~2.2]{dSLu96} to the coarse decomposition
$V = U_1 \oplus U_2$.

\begin{theorem} \label{thm:dS-Lu} In the set-up described above, we
  have
  \[
    \int_{G^+} \lvert \det g \rvert_p^{\, s} \, \mathrm{d}\mu_G (g) =
    \int_{H^+} \lvert \det h \rvert_p^{\, s} \, \theta_0(h) \, \theta_1(h) \,
    \mathrm{d}\mu_H(h).
  \]
\end{theorem}

In our applications we will see that $\theta_1(h)$ is straightforward
to calculate, while $\theta_0(h)$ appears to be rather complicated to
track down for large~$m$.  For short, we set
$\theta(h) = \theta_0(h) \theta_1(h)$ for $h \in H$.  In view
of~\cite[Thm.~2.3]{dSLu96}, one could suspect the function
$\theta \colon H \to \R_{>0}$ to be a character on~$H$, but it was
demonstrated in~\cite{BeKlOn18} that, for general class-two nilpotent
groups, one cannot expect this to be the case.  Indeed, in
Sections~\ref{section:x^2} and~\ref{section:x^3} we will see that
$\theta$ is a character for the group $\Gamma_{t^2}$, but that it is
not a character for the group~$\Gamma_{t^3}$; see
Remark~\ref{remark:theta-not-character}.  Subject to the modifications
detailed above, the three technical assumptions (a), (b), (c) of
\cite[\S 2]{dSLu96} are, indeed, satisfied in our setting for every
prime~$p$.  For a general class-two nilpotent Lie lattice, our methods
leading to Theorem~\ref{thm:dS-Lu} work for almost all primes $p$ and
may prove to be useful in other contexts, where
\cite[Assumption~2.3]{dSLu96} does not hold.


\subsection{Utilising a $p$-adic Bruhat decomposition}\label{section:p-adic}
We recall the machinery developed by Igusa~\cite{Ig89}, du Sautoy and
Lubotzky~\cite{dSLu96} and the first author~\cite{Be11} for utilising
a $p$-adic Bruhat decomposition in order to compute integrals over
reductive groups; the reference \cite{Be11} is useful for practical
purposes, where the notation (and some further choices) are
well-suited to the current paper.  We apply this theory in
Sections~\ref{section:x^2} and~\ref{section:x^3}.

Suppose that the group $\mathbf{H}$ is isomorphic to an affine
$\Z$-group scheme $\dot{\mathbf{H}} \le \GL_{\dot{d}}$ and denote by
$\rho \colon \dot{\mathbf{H}} \to \mathbf{H}$ a corresponding
isomorphism.  In our applications, we have
$\dot{\mathbf{H}} = \GL_2 \times \GL_1 \le \GL_3$ and $\rho$ is the
isomorphism described in Corollary~\ref{cor:corollary-reductive-part}.
It is useful to keep this special situation in mind for a concrete
interpretation of the following general approach.  We write
$\dot{H} = \dot{\mathbf{H}}(\Q_p)$, equipped with the right Haar
measure $\mu_{\dot{H}}$ normalised such that
$\mu_{\dot{H}}(\dot{\mathbf{H}}(\Z_p)) = 1$.  We take interest in the
$p$-adic integral
\[
  \mathcal{Z}_{\dot{\mathbf{H}},\rho,\theta,p}(s) = \int_{H^+ \rho^{-1}}
  \lvert \det h^\rho \rvert_p^{\, s} \, \theta(h^\rho) \,
  \mathrm{d}\mu_{\dot{H}}(h),
\]
where $H^+\rho^{-1}$ denotes the full pre-image of $H^+$ under $\rho$
(in the literature this pre-image is usually denoted by~$\dot{H}^+$,
for short, but we prefer the more descriptive form to avoid
misunderstandings).  In our applications, $\rho$ induces a
measure-preserving map from $\dot{H}$ to $H$, as
$\dot{\mathbf{H}}(\Z_p) \rho = \mathbf{H}(\Z_p)$; in this situation,
one could even get away with `identifying' $\mathbf{H}$ and
$\dot{\mathbf{H}}$. 

We fix a maximal torus $\mathbf{T}$ in~$\dot{\mathbf{H}}$ and assume
that $\mathbf{T}$ splits over~$\Q$; this can be arranged in our
applications.  Under an assumption of `good reduction', elements of
$\mathbf{T}$ act by conjugation on minimal closed unipotent subgroups
of~$\dot{\mathbf{H}}$; this action gives rise to a root system
$\Phi \subseteq \Hom(\mathbf{T}, \Gm)$.  The (finite) Weyl group $W$
of $\dot{\mathbf{H}}$ corresponds to
$\mathrm{N}_{\dot{\mathbf{H}}}(\mathbf{T})/\mathbf{T}$, where
$\mathrm{N}_{\dot{\mathbf{H}}}(\mathbf{T})$ is the normaliser of
$\mathbf{T}$ in~$\dot{\mathbf{H}}$.  We suppress here some necessary
requirements of good reduction since these will all trivially hold in
our applications; the technical requirements are detailed in
\cite{Be11}.  We choose a set of simple roots
$\alpha_1,\dots, \alpha_\ell$ which define the positive
roots~$\Phi^+$.  Let $\Xi = \Hom(\Gm, \mathbf{T})$ denote the set of
co-characters of~$\mathbf{T}$.  We refer to~\cite{dSLu96} for a
description of the Iwahori subgroup
$\mathcal{B} \le \dot{\mathbf{H}}(\Z_p)$ with respect to the simple
roots~$\alpha_1, \dots, \alpha_\ell$.  Let $\pi$ denote a fixed
uniformising parameter for~$\Z_p$, e.g., $\pi = p$.  The $p$-adic
Bruhat decomposition theorem of Iwahori and Matsumoto~\cite{IwMa65}
gives
\[
  \dot{H} = \dot{\mathbf{H}}(\Q_p) = \coprod_{\substack{w\in W\\ \xi
      \in \Xi}} \mathcal{B} \, w \, \xi(\pi) \, \mathcal{B} \qquad
  \text{and} \qquad \dot{\mathbf{H}}(\Z_p) = \coprod_{w\in W}
  \mathcal{B} \, w \, \mathcal{B},
\] 
where elements $w \in W$ in this context are to be read as coset
representatives
$g_w \in \mathrm{N}_{\dot{\mathbf{H}}}(\mathbf{T})(\Z_p)$.  One
defines $\Xi^+=\{ \xi \in \Xi \mid \xi(\pi) \in H^+ \rho^{-1} \}$ and
considers, for $w \in W$,
\[
  w \Xi_w^+= \big\{ \xi \in \Xi^+ \mid
  \text{$\alpha_i(\xi(\pi)) \in \Z_p$ for $1 \le i \le \ell$, and
    $\alpha_i(\xi(\pi))\in p\Z_p$ whenever $\alpha_i\in w(\Phi^-)$}
  \big\},
\]
where $\Phi^-$ denotes the set of negative roots.  Utilising
symmetries in the affine Weyl group and the fact that
$\lvert \det \cdot\,^\rho \rvert_p$, $\theta(\cdot\,^\rho)$ are
constant on double cosets of $\mathcal{B} \le \dot{\mathbf{H}}(\Z_p)$, 
(compare with~\cite[Lem.~3.10]{Be11})
the following generalisation of \cite[(5.4)]{dSLu96} holds.

\begin{proposition}[du Sautoy, Lubotzky; Berman
  \protect{\cite[Prop.~4.2]{Be11}}] \label{proposition:formula-integral}
  If $\mathbf{T}$ splits over $\Q$ then, assuming good reduction,
  \[
    \mathcal{Z}_{\dot{\mathbf{H}},\rho, \theta, p}(s) = \sum_{w \in W}
    p^{-\lf(w)} \sum_{\xi\in w\Xi_w^+} \big| \big( \prod_{\beta \in
      \Phi^+}\beta \big) \big((\xi(\pi) \big) \big|_p^{\, -1} \,\,
    \pavbig{ \det \xi(\pi)^\rho}^{\, s} \,\, \theta \big( \xi(\pi)^\rho
    \big),
  \]
  where $\lf(\cdot)$ is the standard length function on~$W$.
\end{proposition}

Finally we recall a natural pairing between
$\Xi = \Hom(\Gm, \mathbf{T})$ and $\Hom(\mathbf{T}, \Gm)$: this is the
map $(\beta, \xi)\mapsto \langle \beta, \xi\rangle$, where
$\beta(\xi(\tau))=\tau^{\langle \beta, \xi\rangle}$ for all
$\tau\in \Gm$. As in \cite[\S 5.2]{Be11}, it will turn out to be
convenient to judiciously choose a basis for~$\Hom(\mathbf{T}, \Gm)$,
consisting of simple roots and dominant weights for the contragredient
representations of irreducible components of~$\rho$, and then to
determine a dual basis for $\Hom(\Gm,\mathbf{T})$.  This will enable
an explicit description of the set $w\Xi_w^+$.


\begin{beispiel}\label{exa:m=1-case}
  To illustrate the general set-up, we indicate how it can be used to compute
  the pro-isomorphic zeta function of the $D^*$-group
  $\Gamma = \Gamma_t$ of Hirsch length~$4$, defined
  in~\eqref{equ:def-Gamma-tm}.
  Proposition~\ref{pro:class-2-correspondence} shows that
  $\zeta_{\Gamma,p}^\wedge(s) = \zeta_{L,p}^\wedge(s)$ for all
  primes~$p$; here $L$ is the $\Z$-Lie lattice of $\Z$-rank~$4$,
  defined by~\eqref{eq:structure.of.L} with respect to the $\Z$-basis
  $\mathcal{S}$, where $K = (0)$ is the companion matrix of the prime
  polynomial $\Delta_K = t$.  We consider the algebraic automorphism
  group $\G = \mathbf{Aut}(L)$, with respect to the $\Z$-basis
  $\mathcal{S}^* = (x_1,\, y_1,\, z_2,z_1)$ as in
  Corollary~\ref{cor:corollary-reductive-part} and
  Example~\ref{exa:unipotent-radical-m23}.

  Let~$p$ be a prime; our aim is to calculate the local pro-isomorphic
  zeta function $\zeta^\wedge_{L,p}(s)$. The coarse decomposition of
  $V={_{\Q_p} L}$ described in
  Section~\ref{section:reduction-integral} is not suitable, due to the
  fact that here the centre does not coincide with the derived
  sublattice of $L$. Instead we require a refined
  decomposition. Setting $U_1=\Span_{\Q_p} \{ x_1,y_1 \}$,
  $U_2=\Span_{\Q_p} \{ z_2 \}$ and $U_3=\Span_{\Q_p} \{z_1 \}$, we
  write $G=\G(\Q_p)$, $H=\mathbf{H}(\Q_p)$, $N=\mathbf{N}(\Q_p)$;
  these groups act on $V={_{\Q_p} L} = U_1\oplus U_2\oplus U_3$ in a
  suitable way. We now require the following subgroups of the
  unipotent radical: $N_{1} = N \cap \ker(\psi'_2)$, where
  $\psi'_2 \colon G \rightarrow \Aut(V/(U_2+U_3))$ denotes the natural
  action, and $N_{2} = N \cap \ker(\psi'_3)$, where
  $\psi'_3 \colon G \rightarrow \Aut(V/U_3)$ denotes the natural
  action. By Corollary~\ref{cor:corollary-reductive-part}, the
  elements of the reductive subgroup $H$ are of the form
  \begin{equation} \label{equ:form-of-h-m=1} \operatorname{diag}(A,
    \nu \det A, \det A), \qquad \text{where
      $(A,\nu) \in \GL_2(\Q_p) \times \GL_1(\Q_p)$,}
  \end{equation}
  and, according to Example~\ref{exa:unipotent-radical-m23}, elements
  of $N$ take the form
  \[
    \left(
      \begin{array}{c|cc}
        \mathrm{I}_2&*&*\\
        \hline
        0&1& \lambda\\
        0&0&1
      \end{array}
    \right), \qquad \text{with $\lambda \in \Q_p$ and arbitrary
      entries in the positions marked~$\ast$.}
  \]
  As explained above, we can utilize
  Proposition~\ref{pro:integral-formula} and Theorem~\ref{thm:dS-Lu}
  to compute $\zeta_{L,p}^\wedge(s)$ via an integral over~$H^+$.  A
  short calculation (using a slightly different analysis of $\theta$,
  based on \cite[\S2]{dSLu96} with respect to the decomposition
  $U_1\oplus U_2\oplus U_3$) shows that, for $h \in H^+$ of the
  form~\eqref{equ:form-of-h-m=1},
  \[
    \theta(h) = \lvert \det A \rvert_p^{\, -5} \lvert
    \nu \rvert_p^{\, -2}.
  \]
  From here on a direct calculations could be carried out; but we
  prefer to illustrate the use of the Bruhat decomposition.  We
  observe that the morphism
  \[
    \varrho \colon \dot{\mathbf{H}} = \GL_2\times \GL_1 \to
    \mathbf{H}, \quad (A,\nu)\mapsto \operatorname{diag}(A, \nu \det
    A, \det A)
  \]
  induces a measure-preserving isomorphism
  $\dot{H} = \dot{\mathbf{H}}(\Q_p) \to H$ such that
  \[
    H^+ \rho^{-1} = \{ (A,\nu) \mid v_p(A) \ge 0 \text{ and } v_p(\det
    A) + v_p(\nu) \ge 0 \},
  \]
  where $v_p \colon \Q_p \to \Z \cup \{\infty\}$ denotes in the first
  place the standard $p$-adic valuation map and also the map
  $\Mat_2(\Q_p) \to \Z \cup \{\infty\}$,
  $(a_{ij}) \mapsto \min \{ v_p(a_{ij}) \mid 1 \le i,j \le 2 \}$.
  Thus we obtain
  \begin{equation*} 
    \zeta^\wedge_{L, p}(s) 
    =\int_{\substack{(A,\nu)\in \dot{H} \text{ with}\\
        v_p(A) \ge 0 \text{ and} \\
        v_p(\det A) + v_p(\nu) \geq 0}}
    \, \lvert \det A \rvert_p^{\, 3s-5} \, \lvert \nu \rvert_p^{\, s-2}
    \, \mathrm{d}\mu_p(A,\nu). 
  \end{equation*}

  For convenience, we consider $\dot{\mathbf{H}}=\GL_2\times \GL_1$ as
  a subgroup of $\GL_3$, embedded as block matrices via
  $(A,\nu) \mapsto \operatorname{diag}(A,\nu)$.  In particular,
  $T = \mathbf{T}(\Q_p) = \{
  \operatorname{diag}(\lambda_1,\lambda_2,\nu) \mid
  \lambda_1,\lambda_2, \nu\in \Q_p^\times \}$ is a maximal torus
  in~$\dot{H}$.  By Proposition~\ref{proposition:formula-integral} we
  obtain
  \[
    \zeta^\wedge_{L, p}(s) = \sum_{w\in W} p^{-\lf(w)} \sum_{\xi\in
      w\Xi_w^+} \lvert \alpha(\xi(\pi)) \rvert_p^{\, -1} \,\, \lvert
    \det(\xi(\pi)^\rho) \rvert_p^{\, s} \,\, \theta(\xi(\pi)^\rho),
  \]
  where we choose
  $\alpha \in \Hom(\mathbf{T},\Gm),
  \alpha(\operatorname{diag}(\lambda_1, \lambda_2,
  \nu))=\lambda_1\lambda_2^{-1}$ as the single positive root, and we
  have
  \[
    w\Xi_w^+=\{\xi\in \Xi^+ \mid \text{$\alpha(\xi(\pi))\in \Z_p$, and
      $\alpha(\xi(\pi))\in p\Z_p$ if $w=w_0$} \},
  \]
  where the Weyl group is $W=\{1, w_0\}$. In order to describe the set
  $w\Xi_w^+$ we consider dominant weights for the contragredient
  representation, following~\cite{dSLu96}. These are given by
  \[
    \omega_1^{-1}(h)=\lambda_2, \quad
    \omega_2^{-1}(h)=\lambda_1\lambda_2\nu\ \qquad \text{for
      $h=\operatorname{diag}(\lambda_1,\lambda_2,\nu) \in T$.}
  \]
  It follows that $\alpha, \omega_1^{-1}, \omega_2^{-1}$ form a
  $\Z$-basis for $\Hom(\mathbf{T}, \Gm)$ whose $\N_0$-span contains
  all the weights of~$\rho$.  Thus to detect whether an element
  $h \in T$ is integral it is sufficient to check whether
  $\alpha(h), \omega_1^{-1}(h), \omega_2^{-1}(h)$ all lie in
  $\Z_p$. We rewrite $\alpha_1=\alpha$, $\alpha_2=\omega_1^{-1}$,
  $\alpha_3=\omega_2^{-1}$ and find that $\xi_1, \xi_2, \xi_3 \in \Xi$
  defined by
  \[
    \xi_1(\tau) = (\tau,1,\tau^{-1}), \quad \xi_2(\tau) =
    (\tau,\tau,\tau^{-2}), \quad \xi_3(\tau) = (1,1, \tau) \qquad
    \text{for $\tau \in \Q_p^\times$.}
  \]
  form a dual basis so that
  \[
    \langle \alpha_i, \xi_j\rangle=
    \begin{cases}
      1& \text{if $i=j$,} \\
      0& \text{if $i\neq j$.}
    \end{cases}
  \]
  A general element of $\Xi$ has the form
  $\xi_\mathbf{e}=\xi_1^{\, e_1} \xi_2^{\, e_2} \xi_3^{\, e_3}$ with
  $\mathbf{e} = (e_1, e_2, e_3) \in \Z^3$ and satisfies
  $\xi_\mathbf{e}(\pi) = \operatorname{diag}(\pi^{e_1+e_2},
  \pi^{e_2},\pi^{-e_1-2e_2+e_3})$.  Hence
  \[
    \xi_\mathbf{e}(\pi)^\rho= \operatorname{diag}( \pi^{e_1+e_2}, \pi^{e_2},
    \pi^{e_3}, \pi^{e_1+2e_2} )
  \]
  and we read off
  \[
    \lvert \det \xi_\mathbf{e}(\pi)^\rho \rvert_p^{\, s}
    =p^{-(2e_1+4e_2+e_3)s}, \qquad \theta \big( \xi_\mathbf{e}(\pi)^\rho \big)
    = p^{3e_1+6e_2+2e_3}.
  \]
  Note that
  $\lvert \alpha(\xi_\mathbf{e}(\pi)) \rvert^{-1}=p^{\langle \alpha,
    \xi_\mathbf{e}\rangle} = p^{\langle \alpha_1, \xi_\mathbf{e} \rangle}=p^{e_1}$ and
  we can rewrite
  \[
    w\Xi^+_w=\{ \xi \in \Xi \mid
    \text{$\langle \alpha_i, \xi \rangle\geq 0$ for $i \in \{1,2,3\}$,
      and $\langle \alpha_1, \xi\rangle > 0$ if $w=w_0$} \},
  \]
  since $\alpha_1\in w(\Phi^-)$ if and only if $w\neq 1$.  Thus we
  obtain
  \begin{align*}
    \mathcal{Z}_{\dot{\mathbf{H}},\rho, \theta, p}(s) = %
    & \sum_{w\in W} p^{-\lf(w)} \sum_{\xi\in
      w\Xi^+_w} p^{\langle \alpha, \xi \rangle} \,
      \pavbig{\det \xi(\pi)^\rho}^{\, s} \, \theta \big( \xi(\pi)^\rho
      \big) \\
    = & \sum_{w\in W}p^{-\lf(w)} \sum_{\substack{\mathbf{e} \in \N_0^{\,3}
        \text{ with} \\ e_1>0\ \text{if}\ w\neq 1}}
    p^{(4-2s)e_1+(6-4s)e_2+(2-s)e_3} \\ 
    = &\frac{1}{(1-p^{6-4s})(1-p^{2-s})}\left(
        p^0\cdot\frac{1}{1-p^{4-2s}}+p^{-1}\cdot\frac{p^{4-2s}}{1-p^{4-2s}}
        \right)\\ 
    = &\frac{1}{(1-p^{3-2s})(1-p^{4-2s})(1-p^{2-s})}, 
  \end{align*}
  confirming the formula that we reported in the introduction, based
  on~\cite[\S 3.3.4]{Be05}.
\end{beispiel}


\section{The local pro-isomorphic zeta functions of
  the group $\Gamma_{t^2}$}\label{section:x^2}
In this section we consider the pro-isomorphic zeta function of the
$D^*$-group $\Gamma = \Gamma_{t^2}$ of Hirsch length~$6$, defined
in~\eqref{equ:def-Gamma-tm}.  We prove Theorem~\ref{x2example} and
obtain Corollary~\ref{cor:t2}; it turns out that we can proceed as
  in Example~\ref{exa:m=1-case}, taking care of a little extra
  complexity along the way.

  Proposition~\ref{pro:class-2-correspondence} shows that
  $\zeta_{\Gamma,p}^\wedge(s) = \zeta_{L,p}^\wedge(s)$ for all
  primes~$p$, where $L$ is the $\Z$-Lie lattice associated
  to~$\Gamma$.  In our setting, $L$ is the $\Q$-indecomposable
  $D^*$-Lie lattice $L$ of $\Z$-rank~$6$, defined
  by~\eqref{eq:structure.of.L} with respect to the $\Z$-basis
  $\mathcal{S}$, where
  $K = \left( \begin{smallmatrix} 0 & 1 \\ 0 &
      0 \end{smallmatrix}\right)$ is the companion matrix of the
  primary polynomial $\Delta_K = t^2$.  We consider the algebraic
  automorphism group $\G = \mathbf{Aut}(L)$, with respect to the
  $\Z$-basis $\mathcal{S}^* = (x_1,y_2,\, x_2,y_1,\, z_2,z_1)$ as in
  Corollary~\ref{cor:corollary-reductive-part} and
  Example~\ref{exa:unipotent-radical-m23}.

Let~$p$ be a prime; we will set about calculating the local
pro-isomorphic zeta function $\zeta^\wedge_{L,p}(s)$.  In the notation
of Section~\ref{section:background}, we set
$U_1=\Span_{\Q_p} \{ x_1,y_2,x_2,y_1 \}$ and
$U_2=\Span_{\Q_p} \{ z_2,z_1 \}$.  We write $G=\G(\Q_p)$,
$H=\mathbf{H}(\Q_p)$, $N=\mathbf{N}(\Q_p)$; these groups act on
$V={_{\Q_p} L} = U_1\oplus U_2$. By
Corollary~\ref{cor:corollary-reductive-part}, the elements of the
reductive subgroup $H$ are of the form
\begin{equation} \label{equ:form-of-h} \left(
    \begin{array}{cc|cc}
      \nu A&0&0&0\\
      0&A&0&0\\
      \hline
      0&0&\nu^2\det A&0\\
      0&0&0&\nu \det A
    \end{array}
  \right),
  \qquad
  \text{where $(A,\nu) \in \GL_2(\Q_p) \times \GL_1(\Q_p)$.}
\end{equation}
The description of the unipotent radical given in
Example~\ref{exa:unipotent-radical-m23} shows that elements of $N$ are
of the form
\[
  \left(
    \begin{array}{cc|cc}
      \mathrm{I}_2&B&*&*\\
      0&\mathrm{I}_2&*&*\\
      \hline
      0&0&1&\tr B\\
      0&0&0&1
    \end{array}
  \right), \qquad \text{where $B \in \Mat_2(\Q_p)$}
\]
and there are arbitrary entries in the positions marked~$\ast$.
As explained in Section~\ref{section:background}, we can utilize
Proposition~\ref{pro:integral-formula} and Theorem~\ref{thm:dS-Lu} to
compute $\zeta_{L,p}^\wedge(s)$ via an integral over~$H^+$.

We now set about calculating the functions $\theta_0, \theta_1$
defined in Section~\ref{section:background}; we refer to
Section~\ref{section:reduction-integral} for definitions of $N_1$,
$\mu_{N/N_1}$ and $\mu_{N_1}$.  Noting that $N/N_1\cong \Q_p^{\, 4}$
and $N_1\cong \Q_p^{\, 8}$, we obtain for $h \in H^+$ of the
form~\eqref{equ:form-of-h},
\[
  \theta_0(h) = \lvert \det A \rvert_p^{\, -2} \qquad \text{and}
  \qquad \theta_1(h) = \lvert \nu^3\det A^2 \rvert_p^{\, -4},
\]
hence
$\theta(h) = \theta_0(h)\theta_1(h) = \lvert \det A \rvert_p^{\, -10}
\lvert \nu \rvert_p^{\, -12}$; in particular,
$\theta \colon H \to \R_{>0}$ is a character.

We observe that the morphism
\[
  \varrho \colon \dot{\mathbf{H}} = \GL_2\times \GL_1 \to \mathbf{H},
  \quad (A,\nu)\mapsto \operatorname{diag}(\nu A, A, \nu^2 \det A, \nu
  \det A)
\]
induces a measure-preserving isomorphism
$\dot{H} = \dot{\mathbf{H}}(\Q_p) \to H$ such that
\[
 H^+ \rho^{-1} = \{ (A,\nu) \mid v_p(A) \ge 0 \text{ and } v_p(A) + v_p(\nu)
\ge 0 \},
\]
where (as in Example~\ref{exa:m=1-case})
$v_p \colon \Q_p \to \Z \cup \{\infty\}$ denotes the standard $p$-adic
valuation map as well as the map
$\Mat_2(\Q_p) \to \Z \cup \{\infty\}$,
$(a_{ij}) \mapsto \min \{ v_p(a_{ij}) \mid 1 \le i,j \le 2 \}$.  Thus
we obtain
\begin{equation} \label{equ:intermediate-integral-t2} \zeta^\wedge_{L,
    p}(s)
  =\int_{\substack{(A,\nu)\in \dot{H} \text{ with}\\
      v_p(A) \ge 0 \text{ and} \\
      v_p(A) + v_p(\nu) \geq 0}} \, \lvert \det A \rvert_p^{\, 4s-10} \,
  \lvert \nu \rvert_p^{\, 5s-12} \, \mathrm{d}\mu_p(A,\nu).
\end{equation}

For convenience, we consider $\dot{\mathbf{H}}=\GL_2\times \GL_1$ as a
subgroup of $\GL_3$, embedded as block matrices via
$(A,\nu) \mapsto \operatorname{diag}(A,\nu)$.  In particular,
$T = \mathbf{T}(\Q_p) = \{
\operatorname{diag}(\lambda_1,\lambda_2,\nu) \mid \lambda_1,\lambda_2,
\nu\in \Q_p^\times \}$ is a maximal torus in~$\dot{H}$.

By Proposition~\ref{proposition:formula-integral} we obtain
\[
  \zeta^\wedge_{L, p}(s) = \sum_{w\in W} p^{-\lf(w)} \sum_{\xi\in
    w\Xi_w^+} \lvert \alpha(\xi(\pi)) \rvert_p^{\, -1} \,\, \lvert
  \det(\xi(\pi)^\rho) \rvert_p^{\, s} \,\, \theta(\xi(\pi)^\rho),
\]
where we choose
$\alpha \in \Hom(\mathbf{T},\Gm),
\alpha(\operatorname{diag}(\lambda_1, \lambda_2,
\nu))=\lambda_1\lambda_2^{-1}$ as the single positive root, and we
have
\[
  w\Xi_w^+=\{\xi\in \Xi^+ \mid \text{$\alpha(\xi(\pi))\in \Z_p$,
  and $\alpha(\xi(\pi))\in p\Z_p$ if $w=w_0$} \},
\]
where the Weyl group is $W=\{1, w_0\}$. In order to describe the set
$w\Xi_w^+$ we will need to consider dominant weights for the
contragredient representation, following \cite{dSLu96}. These are
given by
\[
  \omega_1^{-1}(h)=\lambda_2 \nu, \quad \omega_2^{-1}(h)=\lambda_2,
  \quad \omega_3^{-1}(h)=\lambda_1\lambda_2\nu^2, \quad
  \omega_4^{-1}(h)=\lambda_1\lambda_2\nu\ \qquad \text{for
    $h=\operatorname{diag}(\lambda_1,\lambda_2,\nu) \in T$.}
\]
It follows that $\alpha, \omega_1^{-1}, \omega_2^{-1}$ form a
$\Z$-basis for $\Hom(\mathbf{T}, \Gm)$ whose $\N_0$-span contains all
the weights of~$\rho$.  Thus to detect whether an element $h \in T$ is
integral it is sufficient to check whether
$\alpha(h), \omega_1^{-1}(h), \omega_2^{-1}(h)$ all lie in $\Z_p$. We
rewrite $\alpha_1=\alpha$, $\alpha_2=\omega_1^{-1}$,
$\alpha_3=\omega_2^{-1}$ and seek a dual basis, namely elements
$\xi_1, \xi_2, \xi_3\in\Xi$ such that
\[
  \langle \alpha_i, \xi_j\rangle=
  \begin{cases}
    1& \text{if $i=j$,} \\
    0& \text{if $i\neq j$.}
  \end{cases}
\]
A routine calculation shows that the following elements suffice:
\[
  \xi_1(\tau) = (\tau,1,1), \quad \xi_2(\tau) = (1,1,\tau), \quad
  \xi_3(\tau) = (\tau,\tau, \tau^{-1}) \qquad \text{for
    $\tau \in \Q_p^\times$.}
\]
A general element of $\Xi$ has the form
$\xi_\mathbf{e} =\xi_1^{\, e_1} \xi_2^{\, e_2} \xi_3^{\, e_3}$ with
$\mathbf{e} = (e_1, e_2, e_3) \in \Z^3$ and satisfies
$\xi_\mathbf{e}(\pi) = \operatorname{diag}(\pi^{e_1+e_3},
\pi^{e_3},\pi^{e_2-e_3})$.  Hence
\[
  \xi_\mathbf{e}(\pi)^\rho= \operatorname{diag}( \pi^{e_1+e_2}, \pi^{e_2},
  \pi^{e_1+e_3}, \pi^{e_3}, \pi^{e_1+2e_2}, \pi^{e_1+e_2+e_3} )
\]
and we read off 
\[
  \lvert \det \xi_\mathbf{e}(\pi)^\rho \rvert_p^{\, s}
  =p^{-(4e_1+5e_2+3e_3)s}, \qquad \theta \big( \xi_\mathbf{e}(\pi)^\rho \big)
  = p^{10e_1+12e_2+8e_3}.
\] 
Note that
$\lvert \alpha(\xi_\mathbf{e}(\pi)) \rvert^{-1}=p^{\langle \alpha,
  \xi_\mathbf{e}\rangle} = p^{\langle \alpha_1, \xi_\mathbf{e}
  \rangle}=p^{e_1}$ and we can rewrite
\[
  w\Xi^+_w=\{ \xi \in \Xi \mid \text{$\langle \alpha_i, \xi \rangle\geq 0$ for
    $i \in \{1,2,3\}$, and $\langle \alpha_1, \xi\rangle > 0$ if
    $w=w_0$} \},
\]
since $\alpha_1\in w(\Phi^-)$ if and only if $w\neq 1$.  Thus we
obtain
\begin{align*}
  \mathcal{Z}_{\dot{\mathbf{H}},\rho, \theta, p}(s) = %
  & \sum_{w\in W} p^{-\lf(w)} \sum_{\xi\in
    w\Xi^+_w} p^{\langle \alpha, \xi \rangle} \,
    \pavbig{\det \xi(\pi)^\rho}^{\, s} \, \theta \big( \xi(\pi)^\rho
    \big) \\
  = & \sum_{w\in W}p^{-\lf(w)} \sum_{\substack{\mathbf{e} \in \N_0^{\,3}
      \text{ with} \\ e_1>0\ \text{if}\ w\neq 1}}
      p^{(11-4s)e_1+(12-5s)e_2+(8-3s)e_3} \\ 
  = &\frac{1}{(1-p^{12-5s})(1-p^{8-3s})}\left(
     p^0\cdot\frac{1}{1-p^{11-4s}}+p^{-1}\cdot\frac{p^{11-4s}}{1-p^{11-4s}}
     \right)\\ 
  = &\frac{1+p^{10-4s}}{(1-p^{8-3s})(1-p^{11-4s})(1-p^{12-5s})}, 
\end{align*}
proving Theorem~\ref{x2example}.  The first part of
Corollary~\ref{cor:t2} follows directly from well-known properties of
the Riemann zeta function.  For the assertion about the asymptotic
growth of pro-isomorphic subgroups in~$\Gamma$, we use a Tauberian
theorem as recorded in~\cite[Thm.~4.20]{dSGr00}.  In the notation
employed there, we take $a = 3$,
$g(s) = \big( \frac{1}{12} + g_1(s) \big) \zeta(5s-12) \zeta(4s -10) /
\zeta(8s-20)$ with $g_1(s)$ holomorphic such that $g_1(3)=0$, and
$w=2$ so that~\eqref{equ:tauberian-t2} holds for
\[
  c_{t^2} = \frac{g(a)}{a \mathsf{\Gamma}(w)} = \frac{\tfrac{1}{12}}{3
    \mathsf{\Gamma}(2)} \frac{\zeta(3) \zeta(2)}{\zeta(4)};
\]
using the precise values $\mathsf{\Gamma}(2) = 1$,
$\zeta(2) = \frac{\pi^2}{6}$ and $\zeta(4) = \frac{\pi^4}{90}$ and the
estimate $\zeta(3) \approx 1.202057$ we arrive at the claimed
description of the invariant~$c_{t^2}$.


\section{The local pro-isomorphic zeta functions of the group $\Gamma_{t^3}$}\label{section:x^3}

\subsection{Counting points on a quadratic
  surface} \label{sec:points-on-surface} In preparation for computing
the pro-isomorphic zeta function of the group $\Gamma_{t^3}$, we study
a certain arithmetic function.  In order to make the analysis
transferable to a more general setting, considered in
Section~\ref{sec:base-ext}, we work over a compact discrete valuation
ring ${\scriptstyle \mathcal{O}}$ with maximal
ideal~$\wp = \pi {\scriptstyle \mathcal{O}}$ and residue field
${\scriptstyle \mathcal{O}} / \wp \cong \F_q$ of size $q$ and
characteristic~$p$.  Our primary interest is in the basic set-up:
${\scriptstyle \mathcal{O}} = \Z_p$, $\wp = p \Z_p$ and
$\Z_p / p\Z_p \cong \F_p$.

\begin{definition}\label{def-f} For $\alpha,\beta, m \in \N_0$ and
  an indeterminate $t$, let
  \begin{equation*}
    \begin{split}
      f(\alpha, \beta, m)&= \# \big\{ (x, y, z) \in ({\scriptstyle
        \mathcal{O}} / \pi^m{\scriptstyle \mathcal{O}})^3 \mid
      \pi^\alpha x^2+\pi^\beta y z = 0 \big\}, \\
      F_{\alpha,\beta}(\T)& =\sum_{m=0}^\infty
      f(\alpha,\beta,m)\T^m,
    \end{split}
  \end{equation*}
  and, for $\beta \in \N_0$, let
  \[
    F^\star_{0,\beta}(\T) = \sum_{m=\beta}^{\infty}
    f(0,\beta,m)\T^m.
  \]
\end{definition}

Observing that for $\alpha, \beta, m\in \N_0$ one trivially has
\begin{equation}\label{eqn:rec111}
  f(\alpha+1,\beta+1, m+1) = q^3 f(\alpha,\beta,m),
\end{equation}
we focus on the cases where either $\alpha$ or $\beta$ is zero.

\begin{proposition}\label{prop:F-fns}
  For $\alpha \in \N_0$, we have
  \[
    \begin{split}
      \mathrm{(i)} \quad F_{\alpha,0}(\T)&= q^{2\alpha}t^\alpha
      F_{0,0}(\T)+\frac{(1-q\T)(1-q^{2\alpha}\T^\alpha)}{(1-q^2\T)^2},\\
      \mathrm{(ii)} \quad F_{0,\alpha}(\T)&= (q^{5}t^2)^{ \lfloor
        \alpha/2 \rfloor} F_{0,\overline{\alpha}}(\T)+
      (1+q^2\T)\frac{1-q^{5 \lfloor \frac{\alpha}{2} \rfloor} \T^{2
          \lfloor \frac{\alpha}{2} \rfloor}}{1-q^5\T^2},
    \end{split}
  \]
  where $\overline{\alpha} = 0$ for $\alpha$ even and
  $\overline{\alpha} = 1$ for $\alpha$ odd.  In particular,
  \[
      F_{0,0}(\T) = \frac{1-q^2\T^2}{(1-q^2\T)(1-q^3\T^2)} \qquad
      \text{and} \qquad 
      F_{0,1}(\T) =\frac{1-2q^3\T^2+q^4\T^2}{(1-q^2\T)(1-q^3\T^2)}.
  \]
\end{proposition}

To prove Proposition~\ref{prop:F-fns} we use the following recurrence
relations.  Parts (1) and (2) of Lemma~\ref{lem:rec.f} below form the
basis for the recursion in $\alpha,\beta$ given in (3) and~(4).
Together with \eqref{eqn:rec111} they determine $f(\alpha,\beta,m)$
completely.

\begin{lemma}\label{lem:rec.f} For $\alpha, \beta, m \in \N_0$ the
  following hold:
  \begin{enumerate}
  \item $f(0,0,m+2) = q^2(q^2-1)q^{2m} + q^3f(0,0,m)$,
  \item $f(0,1,m+2) = 2q^3(q-1)q^{2m} + q^3f(0,1,m)$,
  \item $f(0,\beta+2,m+2) = q^5f(0,\beta, m)$,
  \item $f(\alpha+1,0,m+1) = q(q-1)q^{2m} + q^2f(\alpha,0,m)$.
  \end{enumerate}
\end{lemma}

\begin{proof}  
  To prove (1), we observe that for the finite field~$\F_q$, the set
  of $\F_q$-rational points of the affine variety defined by $x^2+yz$,
  viz.\
  \[
    \{(x,y,-y^{-1}x^2) \mid x \in \F_q, y \in \F_q^\times\} \cupdot
    \{(0,0,z) \mid z \in \F_q\},
  \]
  has $q^2$ points and is smooth away from the origin.  By Hensel's
  lemma each of the $(q^2-1)$ smooth points lifts to $q^{2(m+1)}$
  solutions of $x^2+yz = 0$ over
  ${\scriptstyle \mathcal{O}} / \pi^{m+2}{\scriptstyle
    \mathcal{O}}$.  All the other solutions over
  ${\scriptstyle \mathcal{O}} / \pi^{m+2} {\scriptstyle
    \mathcal{O}}$ are of the form $(\pi x, \pi y, \pi z)$, thus
  $x, y, z$ are perturbations in
  ${\scriptstyle \mathcal{O}} / \pi^{m+1} {\scriptstyle
    \mathcal{O}}$ of solutions modulo~$\pi^m$ and the claim
  follows.

  The argument for part (2) is similar, but as the $\F_q$-points of
  the variety defined by $x^2 + \pi yz \equiv x^2$ are all non-smooth,
  we consider higher levels.  The set of solutions of
  $x^2 + \pi yz = 0$ in
  $({\scriptstyle \mathcal{O}} / \pi^{m+2} {\scriptstyle
    \mathcal{O}})^3$ is a subset of the set
  \begin{multline*}
    \{(\pi x,y,z) \in ({\scriptstyle \mathcal{O}} / \pi^{m+2}
    {\scriptstyle \mathcal{O}})^3 \mid \text{exactly one of $y$ or $z$
      is a unit} \} \\ \cupdot \{ (\pi x,\pi y,\pi z) \in
    ({\scriptstyle \mathcal{O}} / \pi^{m+2} {\scriptstyle
      \mathcal{O}})^3 \mid x^2+\pi yz \equiv 0 \mod \pi^m \}.
  \end{multline*}
  The number of solutions of the second type is $q^3f(0,1,m)$.  For
  the first type, assuming that $z$ is a unit and $\pi \mid y$, we are
  left to solve $\pi \tilde{x}^2+\tilde{y}z\equiv 0 \mod \pi^m$,
  where $x = \pi\tilde{x}$ and $y = \pi\tilde{y}$.  Note that
  $\tilde{y}$ is completely determined by $\tilde{x}, z$.  Counting in
  redundancy from the reduction, we find $(q-1)q^{2m+3}$ solutions.
  By symmetry, the total number of solutions for this type is
  $2(q-1)q^{2m+3}$.

  To prove part (3) consider the equation $x^2 + \pi^{\beta+2}yz = 0$
  over
  ${\scriptstyle \mathcal{O}} / \pi^{m+2}{\scriptstyle
    \mathcal{O}}$.  Note that a triple $(x,y,z)$ is a solution if and
  only if $x = \pi\tilde{x}$, and the triple $(\tilde{x},y,z)$ is a
  solution of the equation $\pi^2(\tilde{x}^2 + \pi^{\beta}yz) = 0$
  over
  ${\scriptstyle \mathcal{O}} / \pi^{m+2}{\scriptstyle
    \mathcal{O}}$.  Thus $f(0,\beta+2,m+2) = q^5 f(0,\beta,m)$,
  where the factor $q^5$ comes from the redundancy arising from the
  reduction to $\tilde{x}^2 + \pi^{\beta}yz \equiv 0 \mod \pi^{m}$.

  For part (4), to solve the equation $\pi^{\alpha+1}x^2 + yz = 0$
  over
  ${\scriptstyle \mathcal{O}} / p^{m+1}{\scriptstyle \mathcal{O}}$
  we consider two cases: that $y$ is divisible by $\pi$ or that $y$ is a
  unit.  Using arguments similar to those above, we find in the first
  case $q^2 f(\alpha, 0, m)$ solutions and in the second
  $q(q-1)q^{2m}$ solutions.
\end{proof}

\begin{proof}[Proof of Proposition~\ref{prop:F-fns}]
  We first compute $F_{0,0}(\T)$.  We multiply both sides of
  equation~(1) in Lemma~\ref{lem:rec.f} by~$\T^{m+2}$ and sum over the
  non-negative integers to obtain
  \[
    \sum_{m=0}^{\infty}
    f(0,0,m+2)\T^{m+2} = q^2(q^2-1) \sum_{m=0}^{\infty}
    q^{2m}\T^{m+2} + q^3 \sum_{m=0}^\infty f(0,0,m)\T^{m+2}.
  \]
  Using the fact that $f(0,0,0) = 1$ and $f(0,0,1) = q^2$ we get
  \[
    F_{0,0}(\T)-1 - q^2\T = \frac{q^2(q^2-1)\T^2}{1-q^2\T} + q^3\T^2
    F_{0,0}(\T), 
  \]
  which implies the formula for $F_{0,0}(\T)$.  The derivation of
  $F_{0,1}(\T)$ is similar.

  To prove part (i) we multiply both sides of equation~(4) in
  Lemma~\ref{lem:rec.f} by $\T^{m+1}$ and sum over the non-negative
  integers.  This gives
  \[
    \underbrace{\sum_{m=0}^{\infty} f(\alpha + 1,0,m + 1)
      \T^{m+1}}_{F_{\alpha+1,0}(\T)-1} = q(q-1)
    \sum_{m=0}^\infty q^{2m}\T^{m+1} +
    \underbrace{q^2\sum_{m=0}^\infty
      f(\alpha,0,m)\T^{m+1}}_{q^2\T F_{\alpha,0}(\T)},
  \]
  and thus yields the recurrence
  \[
    F_{\alpha+1,0}(\T) = \frac{1-q\T}{1-q^2\T} + q^2 \T F_{\alpha,0}(\T).
  \]
  A recurrence of this form, namely, $A_{\alpha+1} = d + c A_\alpha$
  ($\alpha \in \N_0$), has the following solution
  \begin{equation}\label{eqn:rec-gen}
    A_\alpha= d \, \frac{1-c^\alpha}{1-c}+c^\alpha A_0, \quad \alpha \in \N_0,
  \end{equation}
  which implies part (i) of the proposition.

  Similarly, to prove part (ii) we multiply both sides of equation~(3)
  in Lemma~\ref{lem:rec.f} by $\T^{m+2}$ and sum over the
  non-negative integers:
  \[
    {F_{0,\alpha+2}(\T)-1-q^2\T} = {\sum_{m=0}^{\infty}
      f(0,\alpha+2,m+2) \T^{m+2}} = {q^5\sum_{m=0}^\infty
      f(0,\alpha,m)\T^{m+2}} = {q^5\T^2 F_{0,\alpha}(\T)}.
  \]
  We get the recurrence relation 
  \[
    F_{0,\alpha+2}(\T) = 1+q^2\T+{q^5\T^2 F_{0,\alpha}(\T)}.
  \]
  This is solved separately for even and odd $\alpha$, via
  \eqref{eqn:rec-gen}, giving
  \[
    F_{0,\alpha}(\T) = (1+q^2\T) \frac{1-q^{5 \lfloor \frac{\alpha}{2}
        \rfloor} \T^{2 \lfloor \frac{\alpha}{2}
        \rfloor}}{1-q^5\T^2} + (q^5\T^2)^{\lfloor \alpha/2
      \rfloor}F_{0,\overline{\alpha}}(\T). \qedhere
  \]
\end{proof} 

We need to pin down the variant $F^\star_{0,\alpha}(\T)$ of
$F_{0,\alpha}(t)$, which was introduced in Definition~\ref{def-f}.

\begin{lemma}\label{lemma:Fstar-fns-new}
  For $\alpha\in \N_0$, set $\overline{\alpha} = 0$ for $\alpha$ even
  and $\overline{\alpha} = 1$ for $\alpha$ odd.  Then
  \[
    \begin{split}
      F^\star_{0,\alpha}(\T) %
      & =\left(q^5\T^2\right)^{\lfloor \alpha/2\rfloor}
      F^\star_{0,\overline{\alpha}}(\T) \\
      &=
      	\begin{cases}
          q^{\frac{5\alpha}{2}}\T^\alpha F_{0,0}(\T) =
          q^{\frac{5\alpha}{2}}\T^\alpha
          \frac{1-q^2\T^2}{(1-q^2\T)(1-q^3\T^2)}
          & \text{for $\alpha$ even,} \\
          q^{\frac{5(\alpha-1)}{2}}\T^{\alpha-1}
          \left(F_{0,1}(\T)-1\right)=
          q^{\frac{5\alpha-1}{2}}\T^{\alpha}
          \frac{(1-q\T)(1+q^2\T)}{(1-q^2\T)(1-q^3\T^2)} & \text{for
            $\alpha$ odd.}
	\end{cases}
    \end{split}
  \]

  Furthermore, employing another indeterminate $Y$, we have
  \[
    \sum_{\alpha=0}^\infty Y^\alpha F^\star_{0,\alpha}(\T)=\frac{
      (1-q\T) \big( 1+q\T +Yq^2\T(1+q^2\T)
      \big)}{(1-q^5\T^2Y^2)(1-p^2\T)(1-q^3\T^2)}.
  \]
\end{lemma}

\begin{proof}
  Multiplying both sides of equation~(4) in Lemma~\ref{lem:rec.f}
  by~$\T^{m+2}$ and summing over $m \ge \alpha$, we obtain
  \[
    F_{0,\alpha+2}^\star(\T) = \sum_{m=\alpha}^\infty f(0,\alpha
    +2,m+2)\T^{m+2} = q^5\T^2 \sum_{m=\alpha}^\infty
    f(0,\alpha,m)\T^{m} = q^5\T^2 F_{0,\alpha}^\star(\T).
  \]
  Writing $\alpha=2j+\epsilon$ with $\epsilon \in \{0,1\}$, we deduce that
  \[
    F_{0,2j+\epsilon}^\star(\T) = q^{5j}\T^{2j}F_{0,\epsilon}^\star(\T).
  \]
  By substituting $F^\star_{0,0}(\T)=F_{0,0}(\T)$ and
  $F^\star_{0,1}(\T)=F_{0,1}(\T)-1$ we arrive at the desired formula.

  The last part follows by substituting the formulae obtained into
  \[
    \sum_{\alpha=0}^\infty Y^\alpha F^\star_{0,\alpha}(\T)=
    \sum_{j=0}^\infty Y^{2j}F^\star_{0,2j}(t)+\sum_{j=0}^\infty
    Y^{2j+1}F^\star_{0,2j+1}(t). \qedhere
  \]
\end{proof}


\subsection{Applying a $p$-adic Bruhat decomposition}\label{section:t3-Bruhat}

We now turn our attention to the pro-isomorphic zeta function of the
$D^*$-group $\Gamma=\Gamma_{t^3}$ of Hirsch length $8$, defined
in~\eqref{equ:def-Gamma-tm}, and we prove Theorem~\ref{x3example}.

Proposition~\ref{pro:class-2-correspondence} shows that
$\zeta_{\Gamma,p}^\wedge(s) = \zeta_{L,p}^\wedge(s)$ for all
primes~$p$, where $L$ is the $\Z$-Lie lattice associated
to~$\Gamma$.  In our setting, $L$ is the $\Q$-indecomposable $D^*$-Lie
lattice $L$ of $\Z$-rank~$8$, defined by~\eqref{eq:structure.of.L}
with respect to the $\Z$-basis $\mathcal{S}$, where
$K = \left( \begin{smallmatrix} 0 & 1&0 \\ 0 & 0&1\\
    0&0&0 \end{smallmatrix}\right)$ is the companion matrix of the
primary polynomial $\Delta_K = t^2$.  We consider the algebraic
automorphism group $\G = \mathbf{Aut}(L)$, with respect to the
$\Z$-basis
$\mathcal{S}^* = (x_1,y_3,\, x_2,y_2,\, x_3,y_1,\, z_2,z_1)$ as in
Corollary~\ref{cor:corollary-reductive-part} and
Example~\ref{exa:unipotent-radical-m23}.

Let~$p$ be a prime; we will set about calculating the local
pro-isomorphic zeta function $\zeta^\wedge_{L,p}(s)$.  In the notation
of Section~\ref{section:background}, we set
$U_1=\Span_{\Q_p} \{ x_1,y_3,\, x_2,y_2,\, x_3,y_1\}$ and
$U_2=\Span_{\Q_p} \{ z_2,z_1 \}$.  We write $G=\G(\Q_p)$,
$H=\mathbf{H}(\Q_p)$, $N=\mathbf{N}(\Q_p)$; these groups act on
$V={_{\Q_p} L} = U_1\oplus U_2$.  In accordance with
Corollary~\ref{cor:corollary-reductive-part}, the elements of the
reductive subgroup $H$ can be written in the form
\begin{equation} \label{equ:form-of-h-for-t3}
  \left(
    \begin{array}{ccc|cc}
      \nu^{-1} A&0&0&0&0\\
      0&A&0&0&0\\
      0&0&\nu A&0&0\\		
      \hline
      0&0&0&\nu^{-1}\det A&0\\
      0&0&0&0&\det A
    \end{array}
  \right),
  \qquad
  \text{where $(A,\nu) \in \GL_2(\Q_p) \times \GL_1(\Q_p)$;}
\end{equation}
observe that we have performed a routine reparametrisation
$\nu \mapsto \nu^{-1}$ and $A \mapsto \nu A$: for our computation
of~$\theta$ we prefer to have the powers of $\nu$ appearing along the
diagonal to be `small'.

The description of the unipotent radical given in
Example~\ref{exa:unipotent-radical-m23} shows that elements of $N$ are
of the form
\begin{equation}\label{eqn:unipotent element}
  u(B, C)=\left(
    \begin{array}{ccc|cc}
      \mathrm{I}_2&B&C&*&*\\
      0&\mathrm{I}_2&B+\lambda \mathrm{I}_2&*&*\\
      0&0&\mathrm{I}_2&*&*\\
      \hline
      0&0&0&1&\lambda\\
      0&0&0&0&1
    \end{array}
  \right), \qquad \begin{array}{l}\text{where $B,C \in \Mat_2(\Q_p)$
                    and $\lambda \in \Q_p$}\\ 
                    \text{with $\tr(B)=0$, $\tr(C)+\det(B)=0$}\end{array}
\end{equation}
and there are arbitrary entries in the positions marked~$\ast$.
As explained in Section~\ref{section:background}, we can utilize
Proposition~\ref{pro:integral-formula} and Theorem~\ref{thm:dS-Lu} to
compute $\zeta_{L,p}^\wedge(s)$ via an integral over~$H^+$.

We now return to our coarse decomposition and set about calculating
the functions $\theta_0, \theta_1$ defined in
Section~\ref{section:background}; we refer to
Section~\ref{section:reduction-integral} for definitions of $N_1$,
$\mu_{N/N_1}$ and $\mu_{N_1}$.  Noting that $N_1\cong \Q_p^{\, 8}$, we
obtain for $h \in H^+$ of the form~\eqref{equ:form-of-h-for-t3},
\[
  \theta_1(h) = \lvert \nu^{-1}\det A^2 \rvert_p^{\, -6} = \lvert \det
  A \rvert_p^{\, -12} \lvert \nu \rvert_p^{\, 6},
\]
hence
$\theta(h) = \theta_0(h)\theta_1(h) = \theta_0(h) \lvert \det A
\rvert_p^{\, -12} \lvert \nu \rvert_p^{\, 6}$.  We defer until the
next section a calculation of $\theta_0$, since this is the most
involved and lengthy aspect of the analysis.

We observe that the morphism
\[
  \varrho \colon \dot{\mathbf{H}} = \GL_2\times \GL_1 \to \mathbf{H},
  \quad (A,\nu)\mapsto \operatorname{diag}(\nu^{-1} A, A, \nu A, \,
  \nu^{-1} \det A, \det A)
\]
induces a measure-preserving isomorphism
$\dot{H} = \dot{\mathbf{H}}(\Q_p) \to H$ such that
\[
  H^+ \rho^{-1} = \{ (A,\nu) \mid v_p(A) \ge 0 \text{ and } v_p(A) -
  \lvert v_p(\nu) \rvert \ge 0 \},
\]
where $v_p$ is defined as in Example~\ref{exa:m=1-case} and in
Section~\ref{section:x^2}.  Thus we obtain
\begin{equation} \label{equ:intermediate-integral-t3}
  \zeta^\wedge_{L, p}(s)
  = \int_{\substack{(A,\nu)\in \dot{H} \text{ with}\\
  v_p(A) \ge 0 \text{ and} \\
  v_p(A) + |v_p(\nu)| \geq 0}}
  \, \lvert \det A \rvert_p^{\, 5s-12} \, \lvert \nu \rvert_p^{\,
  -s+6}\theta_0 \big( (A,\nu)^\rho \big)
  \, \mathrm{d}\mu_p(A,\nu). 
\end{equation}

For convenience, we consider $\dot{\mathbf{H}}=\GL_2\times \GL_1$ as a
subgroup of $\GL_3$, embedded as block matrices via
$(A,\nu) \mapsto \operatorname{diag}(A,\nu)$.  In particular,
$T = \mathbf{T}(\Q_p) = \{
\operatorname{diag}(\lambda_1,\lambda_2,\nu) \mid \lambda_1,\lambda_2,
\nu\in \Q_p^\times \}$ is a maximal torus in~$\dot{H}$.

By Proposition~\ref{proposition:formula-integral} we obtain
\[
  \zeta^\wedge_{L, p}(s) = \sum_{w\in W} p^{-\lf(w)} \sum_{\xi\in
    w\Xi_w^+} \lvert \alpha(\xi(\pi)) \rvert_p^{\, -1} \,\, \lvert
  \det(\xi(\pi)^\rho) \rvert_p^{\, s} \,\, \theta(\xi(\pi)^\rho),
\]
where we choose
$\alpha \in \Hom(\mathbf{T},\Gm),
\alpha(\operatorname{diag}(\lambda_1, \lambda_2,
\nu))=\lambda_1\lambda_2^{-1}$ as the single positive root, and we
have
\[
  w\Xi_w^+=\{\xi\in \Xi^+ \mid \text{$\alpha(\xi(\pi))\in \Z_p$, and
    $\alpha(\xi(\pi))\in p\Z_p$ if $w=w_0$} \},
\]
where the Weyl group is $W=\{1, w_0\}$. In order to describe the set
$w\Xi_w^+$ we will need to consider dominant weights for the
contragredient representation, following \cite{dSLu96}. These are
given by
\[
  \omega_1^{-1}(h)=\lambda_2 \nu^{-1}, \quad
  \omega_2^{-1}(h)=\lambda_2, \quad \quad
  \omega_3^{-1}(h)=\lambda_2\nu, \quad
  \omega_4^{-1}(h)=\lambda_1\lambda_2\nu^{-1}, \quad
  \omega_5^{-1}(h)=\lambda_1\lambda_2
\]
for $h=\operatorname{diag}(\lambda_1,\lambda_2,\nu) \in T$. It follows
that $\alpha, \omega_1^{-1}, \omega_2^{-1}$ form a $\Z$-basis for
$\Hom(\mathbf{T},\Gm)$. Unlike the situation in
Section~\ref{section:x^2}, the $\N_0$-span of these three dominant
weights does not contain all the weights of~$\rho$. In the current
situation an element $h\in T$ is integral if and only if
$\alpha(h), \omega_1^{-1}(h), \omega_2^{-1}(h), \omega_3^{-1}(h)$ all
lie in $\Z_p$. Note that $\omega_3^{-1}=\omega_1\omega_2^{-2}$.  We
rewrite $\alpha_1=\alpha$, $\alpha_2=\omega_1^{-1}$,
$\alpha_3=\omega_2^{-1}$ and seek a dual basis, namely elements
$\xi_1, \xi_2, \xi_3\in\Xi$ such that
\[
\langle \alpha_i, \xi_j\rangle=
\begin{cases}
  1& \text{if $i=j$,} \\
  0& \text{if $i\neq j$.}
\end{cases}
\]
A routine calculation shows that the following elements suffice:
\[
  \xi_1(\tau)=(\tau,1,1), \quad \xi_2(\tau)=(1,1,\tau^{-1}), \quad
  \xi_3(\tau)=(\tau,\tau,\tau)\quad \text{for $\tau\in\Q_p^\times$}.
\]
A general element of $\Xi$ has the form
$\xi_\mathbf{e} = \xi_1^{\, e_1} \xi_2^{\, e_2} \xi_3^{\, e_3}$ with
$\mathbf{e} = (e_1, e_2, e_3) \in \Z^3$ and then
\begin{equation}\label{eqn-form-of-xi}
  \xi_\mathbf{e}(\pi)=\diag(\pi^{e_1+e_3}, \pi^{e_3},\pi^{e_3-e_2}).
\end{equation}
Hence
\[
  \xi_\mathbf{e}(\pi)^\rho= \operatorname{diag}( \pi^{e_1+e_2}, \pi^{e_2},
  \pi^{e_1+e_3}, \pi^{e_3}, \pi^{e_1-e_2+2e_3}, \pi^{-e_2+2e_3},
  \pi^{e_1+e_2+e_3}, \pi^{e_1+2e_3} )
\]
and we read off
\[
  \lvert \det \xi_\mathbf{e}(\pi)^\rho \rvert_p^{\, s} = p^{-(5e_1+e_2+9e_3)s},
  \qquad \theta_1 \big( \xi_\mathbf{e}(\pi)^\rho \big) = p^{12e_1+6e_2+18e_3}.
\] 
Note that
$\lvert \alpha(\xi_\mathbf{e}(\pi)) \rvert^{-1}=p^{\langle \alpha,
  \xi_\mathbf{e}\rangle}=p^{\langle \alpha_1,
  \xi_\mathbf{e}\rangle}=p^{e_1}$ and we can rewrite
\begin{align*}
  w\Xi^+_w&=\{ \xi \in \Xi \mid \text{$\langle \alpha_i,
            \xi \rangle\geq 0$ for  $i \in \{1,2,3\}$; $\langle
            \omega_3^{-1}, \xi\rangle \geq 0$, and $\langle \alpha_1, 
            \xi\rangle > 0$ if $w=w_0$} \} \\
          &=\{ \xi_\mathbf{e} \mid 
            \text{$e_i \geq 0$ for  $i \in \{1,2,3\}$; $2e_3\geq e_2$,
            and $e_1 > 0$ if $w=w_0$} 
            \},
\end{align*}
since $\alpha_1\in w(\Phi^-)$ if and only if $w\neq 1$ and
$\omega_3^{-1}=\omega_1\omega_2^{-2}=\alpha_2^{-1}\alpha_3^2$.
Writing
\begin{equation}\label{eq:definition-cone}
  \cone=\{\mathbf{e} \in \N_0^{\, 3}\mid 2e_3\geq e_2\}\end{equation}
we obtain
\begin{align}
  \mathcal{Z}_{\dot{\mathbf{H}},\rho, \theta, p}(s) = %
  & \sum_{w\in W} p^{-\lf(w)} \sum_{\xi\in
    w\Xi^+_w} p^{\langle \alpha, \xi \rangle} \,
    \pavbig{\det \xi(\pi)^\rho}^{\, s} \, \theta \big( \xi(\pi)^\rho
    \big) \nonumber \\
  = & \sum_{w\in W}p^{-\lf(w)} \sum_{\substack{\mathbf{e} \, \in\, \cone
      \text{ with}\\ e_1>0\ \text{if}\ w\neq 1}}
  p^{(13-5s)e_1+(6-s)e_2+(18-9s)e_3}\theta_0 \big( \xi_\mathbf{e}(\pi)^\rho
  \big). \label{eqn-t3-intermediate-expr} 
\end{align}


\subsection{Determining the function $\theta_0$}
In view of \eqref{equ:form-of-h-for-t3} and
\eqref{eqn-t3-intermediate-expr}, we need only compute $\theta_0$ for
elements of $H$ of a rather special form; for $n,m,k \in \Z$ we set
\begin{align*}
  \theta_0(\pi^{n},\pi^{m},\pi^{k})  %
  & =
    \theta_0 \big(
    \operatorname{diag}(\pi^{n},\pi^{m},\pi^{k})^\rho \big)
  \\
  & = \theta_0 \big(
    \operatorname{diag}(\pi^{n-k},\pi^{m-k}, \, \pi^{n},
    \pi^{m}, \, \pi^{n+k}, \pi^{m+k}, \, \pi^{m+n-k},
    \pi^{m+n}) \big),
\end{align*}
where the first expression is a mild, but convenient abuse of
notation.  Recall that we could choose $\pi = p$, but prefer to make
clear the different roles played by $\pi$ and~$p$.  This is beneficial
also with a view toward the more general situation considered in
Section~\ref{sec:base-ext}; we refrain from generalising all the
notation in the current section as we did in
Section~\ref{sec:points-on-surface}, but explain in
Remark~\ref{rem:char-2-extra} how one particular step carries over.
We assume throughout that $n \ge m$ since this is the only case of
interest to us; see~\eqref{eqn-form-of-xi}.  Write $l=n-m \in \N_0$,
and recall from Definition~\ref{def-f} with
${\scriptstyle \mathcal{O}} = \Z_p$ that
$f(\alpha, \beta, m) = \#\{(x, y, z) \in (\Z_p / \pi^m\Z_p)^3 \mid
\pi^\alpha x^2+\pi^\beta y z = 0 \}$ for $\alpha,\beta, m \in \N_0$.

\begin{lemma}\label{lemma:theta-tilde}
  For $n,m,k \in \Z$ with $l = n-m \ge 0$, we have
  \[
    \theta_0(\pi^{n},\pi^{m},\pi^{k})=p^{4k+3m+n} \,
    \tilde\theta(\pi^{n},\pi^{m},\pi^{k}),
  \]
  where $\tilde\theta(\pi^{n},\pi^{m},\pi^{k})$ equals
 
  \[
    \left\{\begin{array}{lll}
      p^{3k+l}f(l,0,m-k) & \text{if $k\geq 0$} &\textup{(Case 1)},\\
      p^{-k+l}f(2k+l,0,m+k) & \text{if
        $\max \{-m, -l\} \leq
        k<0$ and $2k+l\geq 0$} &\textup{(Case 2a)}, \\
      p^{5k+4l} f(0,-2k-l, m-k-l) & \text{if $\max
        \{-m, -l\}\leq k<0$ and $2k+l< 0$}; &\textup{(Case 2b)}, \\
      p^{-l} f(0,l,n+k) &
      \text{if $-m\leq k<-l$} &\textup{(Case 3).}
    \end{array}\right.
  \]
\end{lemma}

\begin{remark}\label{remark:theta-not-character}
  It follows from Lemma~\ref{lemma:theta-tilde} that
  $\theta \colon H \to \R_{>0}$ is not a character.  For instance,
  \[
    \tilde\theta(\pi^2,\pi^2,1)=f(0,0,2)=p^4+p^3-p^2\neq
    p^4=f(0,0,1)^2=\tilde\theta(\pi,\pi,1)^2.
  \]
  In fact, this calculation shows that the lifting condition
  \cite[Assumption 2.3]{dSLu96} fails for all primes $p$. Suppose that
  the lifting condition were to hold. By \cite[Lem.~3.12]{Be11}, it
  would follow that $\theta$ is a character on subsets of a maximal
  torus of $H$ with a designated ordering of valuations along the
  diagonal. It is readily seen that the elements
  $\operatorname{diag}(\pi^{2},\pi^{2},1)^\rho$,
  $\operatorname{diag}(\pi,\pi,1)^\rho$ belong to such a subset.
\end{remark}

\begin{proof}[Proof of Lemma~\ref{lemma:theta-tilde}]
  We consider the action of a diagonal element
  \[
    h = \operatorname{diag}(\pi^{n-k}, \pi^{m-k}, \,
    \pi^{n}, \pi^{m}, \, \pi^{n+k}, \pi^{m+k}, \,
    \pi^{m+n-k}, \pi^{m+n})
  \]
  on an element
  \begin{equation}\label{eq:n}
    u =
    \left(
      \begin{array}{ccc|cc}
	\mathrm{I}_2 & \left(\begin{smallmatrix} a&b\\ c&-a \end{smallmatrix}\right)&\left(\begin{smallmatrix} d&e\\ f&a^2+bc-d \end{smallmatrix}\right)&*&*\\
	0& \mathrm{I}_2 &\left(\begin{smallmatrix} \lambda+a&b\\ c&\lambda-a \end{smallmatrix}\right)&*&*\\
	0&0& \mathrm{I}_2 &*&*\\
	\hline
	0&0&0&1&\lambda\\
	0&0&0&0&1
      \end{array}
    \right),
  \end{equation}
  the latter being an explicit parametrisation of~\eqref{eqn:unipotent
    element}.  The situation of interest to us, i.e., when $h$ is
  integral, is equivalent to the conditions
  $n\geq m \geq \lvert k \rvert$.  We obtain the following
  necessary and sufficient conditions for $uh$ to be integral:
  \begin{align}
    \label{eq:5.a}   v_p(a)&\geq -m,\\
    \label{eq:5.b}   v_p(b)&\geq -m+\max\{0,-k\},\\
    \label{eq:5.c}   v_p(c)&\geq -n+\max\{0,-k\},\\
    \label{eq:5.d}    v_p(d)&\geq-n-k,\\
    \label{eq:5.e}   v_p(e)&\geq -m-k,\\
    \label{eq:5.f}    v_p(f)&\geq -n-k,\\
    \label{eq:5.g}    v_p(a^2+bc-d)&\geq -m-k,\\
    \label{eq:5.h}    v_p(\lambda+a)&\geq -n-k,\\
    \label{eq:5.i}    v_p(\lambda-a)&\geq -m-k,\\    
    \label{eq:5.j}    v_p(\lambda)&\geq -m-n.                               
  \end{align}
  Condition~\eqref{eq:5.j} is implied by conditions \eqref{eq:5.a} and
  \eqref{eq:5.i}; it is therefore redundant.  One readily sees the
  following equivalences:
  \begin{align*}
    \eqref{eq:5.d}: \; v_p(d) \geq-n-k %
    & \quad\iff\quad 	
      v_p(a^2+bc) \geq -n-k, && \text{if
                                          \eqref{eq:5.g}
                                          holds;} \\
    \eqref{eq:5.h}: \; v_p(\lambda+a) \geq -n-k %
    & \quad\iff\quad v_p(2a)\geq -n-k, && \text{if \eqref{eq:5.i}
                                                holds;} 
  \end{align*}
  so we may replace \eqref{eq:5.d} and \eqref{eq:5.h} respectively by
  \begin{align*}
    \tag*{$\mathrm{\eqref{eq:5.d}'}$}	v_p(a^2+bc) & \geq -n-k\\	
    \tag*{$\mathrm{\eqref{eq:5.h}'}$} v_p(a) & \geq  -n-k-\delta,
  \end{align*}
  where $\delta = v_p(2) \in \{0,1\}$ takes the value $1$ for $p=2$
  and the value $0$ otherwise.

  In our calculation we use the fact that the measure $\mu_{N/N_1}$
  may be treated as an additive measure on the parameter space
  $\Q_p^{\, 7}$ with $(N/N_1)(\Z_p)$ corresponding to
  $\Z_p^{\, 7}$. Indeed, using the notation introduced
  in~\eqref{eqn:unipotent element}, we see that the map
  $\Mat_2(\Q_p) \times \mathsf{sl}_2(\Q_p) \to N/N_1(\Q_p)$,
  $(X,Y) \mapsto u(X,Y)$ is a homeomorphism. The claim thus follows
  from~\cite[Thm.~8.32]{Kn02} and the fact that the groups involved
  are unimodular.
  
  For fixed parameters $(a, b, c) \in \Q_p^{\, 3}$, we obtain
  \[
    \mu_{\Q_p^{\, 4}} \big\{ (d,e,f,\lambda)\in \Q_p^{\, 4} \mid
    \text{\eqref{eq:5.e}, \eqref{eq:5.f}, \eqref{eq:5.g},
      \eqref{eq:5.i} hold} \big\} = p^{3m+n+4k}.
  \]
  It follows that
  $\theta_0(\pi^n,\pi^m,\pi^k) = p^{3m+n+4k} \,
  \tilde\theta(\pi^n,\pi^m,\pi^k)$, where
  \[
    \tilde\theta(\pi^n,\pi^m,\pi^k) = \mu_{\Q_p^{\, 3}} \big\{ (a,b,c)
    \in \Q_p^{\, 3} \mid \text{\eqref{eq:5.a}, \eqref{eq:5.b},
      \eqref{eq:5.c}, $\mathrm{\eqref{eq:5.d}'}$,
      $\mathrm{\eqref{eq:5.h}'}$ hold} \big\}.
  \]
  For convenience, we summarise the conditions \eqref{eq:5.a},
  \eqref{eq:5.b}, \eqref{eq:5.c}, $\mathrm{\eqref{eq:5.d}'}$,
  $\mathrm{\eqref{eq:5.h}'}$:
  \begin{equation}\label{equ:conds-with-delta} \tag{$\dagger$}
    \begin{aligned}
      v_p(a) & \geq \max \{ -m,-n-k-\delta\}, &
      \phantom{xxxxx} & v_p(b) \geq
      -m+\max\{0,-k\}, \\
      v_p(c) & \geq -n+\max\{0,-k\}, && v_p(a^2+bc) \geq
      -n-k.
    \end{aligned}
  \end{equation}
  The next step is to show that we can drop $\delta$, even for $p=2$.
  Suppose for a contradiction, that there are $a,b,c \in \Q_p$
  satisfying \eqref{equ:conds-with-delta} and such that
  $v_p(a) = -n-k-1 \ge -m$; in particular, $k<0$.  Then
  $v_p(a^2) = -2n-2k-2 < -n-k$ and we conclude from
  $\mathrm{\eqref{eq:5.d}'}$ that
  $v_p(b c) = v_p(a^2) = -2n-2k-2$.  On the other hand
  \eqref{eq:5.b} and \eqref{eq:5.c} yield
  $v_p(b c) \ge -n -m - 2 k$.  This gives
  $-2n-2k-2 \ge -n-m-2k$, hence $m-2 \ge n$, a
  contradiction.

  \begin{remark}\label{rem:char-2-extra}
    The last consideration carries through also in a more general
    setting, considered in Section~\ref{sec:base-ext}.  If we work
    over a compact discrete valuation ring
    ${\scriptstyle \mathcal{O}}$ with valuation $v_\wp$, replacing
    $\Z_p$ with valuation~$v_p$, then $\delta = v_\wp(2)$.  If
    ${\scriptstyle \mathcal{O}}$ has residue characteristic~$2$ this
    is the absolute ramification index
    of~${\scriptstyle \mathcal{O}}$, and the assumption
    $v_\wp(a) = -n-k-{\bar \delta} \ge -m$ with
    ${\bar \delta} \in \{1,\ldots,\delta\}$ leads again to a
    contradiction.
  \end{remark}
  
  Thus we can work with the simpler set of conditions
  \begin{equation}\label{equ:conds-without-delta} \tag{$\ddagger$}
    \begin{aligned}
      v_p(a) & \geq \max \{ -m,-n-k\}, &
      \phantom{xxxxx} & v_p(b) \geq
      -m+\max\{0,-k\}, \\
      v_p(c) & \geq -n+\max\{0,-k\}, && v_p(a^2+bc) \geq
      -n-k.
    \end{aligned}
  \end{equation}
  
  We perform a change of variables $\Q_p^{\ 3}\to \Q_p^{\, 3}$ by
  \[
    (a, b, c)\mapsto (x,y,z) = \left( a p^{\min\{m,
      n+k\}},b p^{m+\min\{0,k\}},c
    p^{n+\min\{0,k\}} \right).
  \]
  The new variables are all unconstrained
  elements of $\Z_p$, and the change of variables introduces a
  Jacobian equal to
  \[
    p^{\min\{m, n+k\}+m+n+\min\{0,2k\}}.
  \]
  It follows that
  \begin{multline*}
    \tilde\theta(\pi^n,\pi^m,\pi^k) = \mu_{\Q_p^{\, 3}} \big\{
    (a,b,c) \in \Q_p^{\, 3} \mid
    \text{\eqref{equ:conds-without-delta} holds} \big\} =
    p^{\min\{m,
      n+k\}+m+n+\min\{0,2k\}} \\
    \cdot \mu_{\Z_p^{\, 3}} \big\{ (x,y,z) \in \Z_p^{\, 3} \mid
    p^{-2\min\{m,n+k\}} x^2 + p^{-m-n-\min\{0,2k\}} yz
    \equiv 0 \mod p^{-n-k} \big\}.
  \end{multline*}
  \showproofs{Lemma~\ref{lemma:theta-tilde} now follows immediately
    by specialising to the four cases.}{\\ \color{purple}
    \textbf{Case 1: $k\geq 0$.} Then
    \begin{align*}
      \tilde\theta(\pi^n,\pi^m,\pi^k) %
      &= p^{2m+n} 
        \cdot \mu_{\Z_p^{\, 3}} \big\{ (x,y,z) \in \Z_p^{\, 3} \mid
        p^{-2m} x^2 +
        p^{-m-n} yz \equiv 0 \mod p^{-n-k} \big\}\\
      &=p^{2m+n} 
        f(l,0,m-k)\cdot \mu_{\Z_p^{\, 3}}\left(p^{m-k}\Z_p^3\right)\\
      &=p^{3k+l}f(l,0,m-k). 	
    \end{align*}
    \\
    \textbf{Case 2a: $\max\{-m,-l\}\leq k<0$ and
      $2k+l\geq 0$.} Then
    \begin{align*}
      \tilde\theta(\pi^n,\pi^m,\pi^k) %
      &= p^{2m+n+2k} 
	\cdot \mu_{\Z_p^{\, 3}} \big\{ (x,y,z) \in \Z_p^{\, 3} \mid
	p^{-2m} x^2 +
	p^{-m-n-2k} yz \equiv 0 \mod p^{-n-k} \big\}\\
      &= p^{2m+n+2k} 
        \cdot \mu_{\Z_p^{\, 3}} \big\{ (x,y,z) \in \Z_p^{\, 3} \mid
        p^{2k+l} x^2 +
        yz \equiv 0 \mod p^{m+k} \big\}\\	
      &=p^{2m+n+2k}f(2k+l,0,m+k) \mu_{\Z_p^{\,
        3}}\left(p^{m+k}\Z_p^3\right)\\ 
      &=p^{-k+l}f(2k+l,0,m+k).
    \end{align*}
    \\
    \textbf{Case 2b: $\max\{-m,-l\}\leq k<0$ and $2k+l< 0$.}
    Then
    \begin{align*}
      \tilde\theta(\pi^n,\pi^m,\pi^k) %
      &= p^{2m+n+2k} 
	\cdot \mu_{\Z_p^{\, 3}} \big\{ (x,y,z) \in \Z_p^{\, 3} \mid
	p^{-2m} x^2 +
	p^{-m-n-2k} yz \equiv 0 \mod p^{-n-k} \big\}\\
      &= p^{2m+n+2k} 
	\cdot \mu_{\Z_p^{\, 3}} \big\{ (x,y,z) \in \Z_p^{\, 3} \mid
        x^2 +
	p^{-2k-l} yz \equiv 0 \mod p^{m-k-l} \big\}\\	
      &=p^{2m+n+2k}f(0,-2k-l,m-k-l) \mu_{\Z_p^{\,
        3}}\left(p^{m-k-l}\Z_p^3\right)\\ 
      &=p^{5k+4l}f(0,-2k-l,m-k-l).
    \end{align*}
    \\
    \textbf{Case 3: $-m\leq k<-l$.} Then
    \begin{align*}
      \tilde\theta(\pi^n,\pi^m,\pi^k)  %
      &= p^{m+2n+3k} 
        \cdot \mu_{\Z_p^{\, 3}} \big\{ (x,y,z) \in \Z_p^{\, 3} \mid
        p^{-2(n+k)} x^2 +
        p^{-m-n-2k} yz \equiv 0 \mod p^{-n-k} \big\}\\
      &= p^{2l+3m+3k} 
        \cdot \mu_{\Z_p^{\, 3}} \big\{ (x,y,z) \in \Z_p^{\, 3} \mid
        x^2 +
        p^{l} yz \equiv 0 \mod p^{n+k} \big\}\\
      &= p^{2l+3m+3k} 
        f(0,l,n+k) \mu_{\Z_p^{\, 3}}\left(p^{n+k}\Z_p^3\right)\\
      &= p^{-l} 
        f(0,l,n+k).\qedhere
    \end{align*}
  }
\end{proof}

In order to continue the calculation paused at
\eqref{eqn-t3-intermediate-expr}, we recall that
$\xi_\mathbf{e} = \xi_1^{\, e_1}\xi_2^{\, e_2}\xi_3^{\, e_3}$ and, setting
\begin{equation}\label{eq:transcribing-to-e_i}
  n_\mathbf{e} =e_1+e_3, \quad
  m_\mathbf{e} = e_3, \quad k_\mathbf{e} =e_3-e_2, \qquad
  \text{thus} \quad l_\mathbf{e} 
  = n_\mathbf{e} - m_\mathbf{e} =e_1,
\end{equation}
we see from \eqref{eqn-form-of-xi} that
$\xi_\mathbf{e}(\pi) = \diag(\pi^{e_1+e_3}, \pi^{e_3}, \pi^{e_3-e_2})
= \diag(\pi^{n_\mathbf{e}}, \pi^{m_\mathbf{e}}, \pi^{k_\mathbf{e}})$.
Applying Lemma~\ref{lemma:theta-tilde} and using
\eqref{eq:transcribing-to-e_i} to resubstitute, we obtain
\begin{equation}\label{eqn:theta-zero-with-e_i}
  \theta_0 \big( \xi_\mathbf{e}(\pi)^\rho \big) = p^{e_1-4e_2+8e_3} \,
  \tilde\theta(\xi_\mathbf{e}(\pi)), 
\end{equation}
where 
\[
  \tilde\theta(\xi_\mathbf{e}(\pi)) = \left\{ \begin{array}{lll}
      p^{e_1-3e_2+3e_3} f(e_1,0,e_2) %
                                                & \text{if
                                                  $e_2\leq e_3$} &
                                                                   \textup{(Case 1)}, \\
                                                p^{e_1+e_2-e_3}f(e_1-2e_2+2e_3,0,-e_2+2e_3)
                                                & \text{if
                                                  $e_3<e_2\leq e_3+\min\{e_1,e_3\}$}& \\
                                                & \text{\; and
                                                  $2e_2\leq
                                                  e_1+2e_3$}&
                                                              \textup{(Case 2a)}, \\
                                                p^{4e_1-5e_2+5e_3}f(0,-e_1+2e_2-2e_3,-e_1+e_2)
                                                & \text{if
                                                  $e_3<e_2\leq e_3+\min\{e_1,e_3\}$}&\\
                                                & \text{\; and
                                                  $e_1+2e_3< 2e_2$}&
                                                                     \textup{(Case 2b)}, \\
                                                p^{-e_1}
                                                f(0,e_1,e_1-e_2+2e_3)
                                                & \text{if
                                                  $e_1+e_3<e_2\leq
                                                  2e_3$} &
                                                           \textup{(Case
                                                           3)}.
  \end{array}\right.
\]
Referring to \eqref{eqn-t3-intermediate-expr}, we obtain
\begin{equation}\label{eqn-t3-second-interm-expr} 
  \mathcal{Z}_{\dot{\mathbf{H}},\rho, \theta, p}(s)  
  = \sum_{w\in W}p^{-\lf(w)} \sum_{\substack{\mathbf{e} \, \in\, \cone
      \text{ with}\\ e_1>0\ \text{if}\ w\neq 1}}
  X_1^{e_1}X_2^{e_2}X_3^{e_3} \,
  \tilde\theta(\xi_\mathbf{e}(\pi)),
\end{equation}
where
\[
  X_1=p^{14-5s}, \quad X_2=p^{2-s},  \quad X_3=p^{26-9s}.
\]


\subsection{Decomposing the polyhedral cone}
In preparation of the final stage of the calculation, we consider the
following subsets of the `integral' cone $\cone$ introduced
in~\eqref{eq:definition-cone}; each subset is, in fact, a submonoid of
$\N_0^{\, 3}$.  Refer to Figure~\ref{figure1} for a pictorial illustration.

\begin{definition}\label{def:submonoids}
  Write
  \[
    \begin{array}{ccc}
      \bv_1=(1,0,0),&\bv_2=(0,2,1),&\bv_3=(0,0,1)\\
      \bv_4=(0,1,1),&\bv_5=(2,2,1),&\bv_6=(1,2,1)
    \end{array}
  \]
  and set
  \begin{align*}
    \cone_{ijk} %
      &=\Span_{\N_0}\{\bv_i, \bv_j, \bv_k\} && \text{for $1 \leq  i,j,k
                                               \leq 6$}, \\ 
    \cone_{ijk^+} %
      &=\Span_{\N_0}\{\bv_i, \bv_j\}+\N \bv_k &&
                                                 \text{for $1\leq i,j,k\leq 6$}, \\
    \cone_{ij} %
      &=\Span_{\N_0}\{\bv_i, \bv_j\} && \text{for $1\leq i,j\leq 6$}, \\
    \cone_{ij^+} %
      &=\N_0\bv_i+\N \bv_j && \text{for $1\leq i,j\leq 6$}, \\					
    \cone_*^{\shift} %
      & = \{(e_1,e_2,e_3)\in\cone_*\mid e_1>0\} && \text{for any
                                                   (possibly empty)
                                                   index $*$}. 
  \end{align*}	
\end{definition}

\begin{figure} 
\tdplotsetmaincoords{60}{135}

\begin{tikzpicture}[scale=2,tdplot_main_coords]
	
	\pgfmathsetmacro{\r}{2}
	\pgfmathsetmacro{\axisoffset}{0.3}
		
	\pgfmathsetmacro{\vax}{1}
	\pgfmathsetmacro{\vay}{0}
	\pgfmathsetmacro{\vaz}{0}
	\pgfmathsetmacro{\val}{1} 
	
	\pgfmathsetmacro{\vbx}{0}
	\pgfmathsetmacro{\vby}{2}
	\pgfmathsetmacro{\vbz}{1}
	\pgfmathsetmacro{\vbl}{2.236} 
	
	\pgfmathsetmacro{\vcx}{0}
	\pgfmathsetmacro{\vcy}{0}
	\pgfmathsetmacro{\vcz}{1}
	\pgfmathsetmacro{\vcl}{1}
	
	\pgfmathsetmacro{\vdx}{0}
	\pgfmathsetmacro{\vdy}{1}
	\pgfmathsetmacro{\vdz}{1}
	\pgfmathsetmacro{\vdl}{1.414}
	
	\pgfmathsetmacro{\vex}{2}
	\pgfmathsetmacro{\vey}{2}
	\pgfmathsetmacro{\vez}{1}
	\pgfmathsetmacro{\vel}{3}
	
	\pgfmathsetmacro{\vfx}{1}
	\pgfmathsetmacro{\vfy}{2}
	\pgfmathsetmacro{\vfz}{1}
	\pgfmathsetmacro{\vfl}{2.449} 

	
	\coordinate (a)		at	(\vax, \vay, \vaz);
	\coordinate (b)		at	(\vbx, \vby, \vbz);
	\coordinate (c)		at	(\vcx, \vcy, \vcz);
	\coordinate (d)		at	(\vdx, \vdy, \vdz);
	\coordinate (e)		at	(\vex, \vey, \vez);
	\coordinate (f)		at	(\vfx, \vfy, \vfz);
	
	\coordinate (aN)	at	(\r/\val*\vax, \r/\val*\vay, \r/\val*\vaz);	
	\coordinate (bN)	at	(\r/\vbl*\vbx, \r/\vbl*\vby, \r/\vbl*\vbz);
	\coordinate (cN)	at	(\r/\vcl*\vcx, \r/\vcl*\vcy, \r/\vcl*\vcz);
	\coordinate (dN)	at	(\r/\vdl*\vdx, \r/\vdl*\vdy, \r/\vdl*\vdz);
	\coordinate (eN)	at	(\r/\vel*\vex, \r/\vel*\vey, \r/\vel*\vez);
	\coordinate (fN)	at	(\r/\vfl*\vfx, \r/\vfl*\vfy, \r/\vfl*\vfz);
	
    \coordinate (O) at (0,0,0);
	
	
    \draw[thin,->,gray] (O) -- (\r + \axisoffset,0,0) node[anchor=north east]{$x$};
    \draw[thin,->,gray] (O) -- (0,\r + \axisoffset,0) node[anchor=north west]{$y$};
	\draw[thin,->,gray] (O) -- (0,0,\r + \axisoffset) node[anchor=south]{$z$};
	
	
	\draw[->, blue, thick]	(O) -- (a)	node[anchor=north] {$\bv_1$};
	\draw[->,blue,thick]	(O) -- (b)
        node[anchor=north west] {$\bv_2$};
	\draw[->,blue,thick]	(O) -- (c)	node[anchor=west] {$\bv_3$};
	\draw[->,blue,thick]	(O) -- (d)
        node[anchor=south east] {$\bv_4$};
	\draw[->,blue,thick]	(O) -- (e)	node[anchor=west] {$\bv_5$};
	\draw[->,blue,thick]	(O) -- (f)	node[anchor=west] {$\bv_6$};


	\draw[-,blue,dashed,thick]	(a) -- (aN)	node[anchor=north,purple] {$\bv_1'$};
	\draw[-,purple,dashed]	(bN) -- (bN)	node[anchor=south west] {$\bv_2'$};
	\draw[-,blue,dashed,thick]	(c) -- (cN)	node[anchor=south west,purple] {$\bv_3'$};
	\draw[-,blue,dashed,thick]	(d) -- (dN)	node[anchor=west,purple] {$\bv_4'$};
	\draw[-,purple,dashed]	(eN) -- (eN) node[anchor=south east] {$\bv_5'$};
	\draw[-,purple,dashed]	(fN) -- (fN)  node[anchor=north] {$\bv_6'$};

        
        \filldraw[blue] (O) circle (.7pt);
        \filldraw[red] (aN) circle (.7pt);
        \filldraw[red] (bN) circle (.7pt);
        \filldraw[red] (cN) circle (.7pt);
        \filldraw[red] (dN) circle (.7pt);
        \filldraw[red] (eN) circle (.7pt);
        \filldraw[red] (fN) circle (.7pt);
        
	
	\tdplotdefinepoints(0,0,0)(\r,0,0)(0,\r,0);
	    \tdplotdrawpolytopearc[thin,red,dashed]{\r}{}{};
	\tdplotdefinepoints(0,0,0)(0,0,\r)(0,\r,0);
	    \tdplotdrawpolytopearc[thin,red,dashed]{\r}{}{};
	\tdplotdefinepoints(0,0,0)(\r,0,0)(0,0,\r);
	    \tdplotdrawpolytopearc[thin,red,dashed]{\r}{}{};
	
	
	\tdplotdefinepoints(0,0,0)(\r/\val*\vax, \r/\val*\vay, \r/\val*\vaz)(\r/\vcl*\vcx, \r/\vcl*\vcy, \r/\vcl*\vcz);
	    \tdplotdrawpolytopearc[thick,red]{\r}{}{};
			
	\tdplotdefinepoints(0,0,0)(\r/\vdl*\vdx, \r/\vdl*\vdy, \r/\vdl*\vdz)(\r/\vcl*\vcx, \r/\vcl*\vcy, \r/\vcl*\vcz);
	    \tdplotdrawpolytopearc[thick,red]{\r}{}{};
		
	\tdplotdefinepoints(0,0,0)(\r/\val*\vax, \r/\val*\vay, \r/\val*\vaz)(\r/\vdl*\vdx, \r/\vdl*\vdy, \r/\vdl*\vdz);
	    \tdplotdrawpolytopearc[thick,red]{\r}{}{};

	
	\tdplotdefinepoints(0,0,0)(\r/\val*\vax, \r/\val*\vay, \r/\val*\vaz)(\r/\vel*\vex, \r/\vel*\vey, \r/\vel*\vez);
	    \tdplotdrawpolytopearc[thick,red]{\r}{}{};
	\tdplotdefinepoints(0,0,0)(\r/\vel*\vex, \r/\vel*\vey, \r/\vel*\vez)(\r/\vdl*\vdx, \r/\vdl*\vdy, \r/\vdl*\vdz);
	    \tdplotdrawpolytopearc[thick,red]{\r}{}{};
		
	
	\tdplotdefinepoints(0,0,0)(\r/\vfl*\vfx, \r/\vfl*\vfy, \r/\vfl*\vfz)(\r/\vel*\vex, \r/\vel*\vey, \r/\vel*\vez);
	    \tdplotdrawpolytopearc[thick,red]{\r}{}{};
	\tdplotdefinepoints(0,0,0)(\r/\vfl*\vfx, \r/\vfl*\vfy, \r/\vfl*\vfz)(\r/\vdl*\vdx, \r/\vdl*\vdy, \r/\vdl*\vdz);
	    \tdplotdrawpolytopearc[thick,red]{\r}{}{};

	
	\tdplotdefinepoints(0,0,0)(\r/\vbl*\vbx, \r/\vbl*\vby, \r/\vbl*\vbz)(\r/\vdl*\vdx, \r/\vdl*\vdy, \r/\vdl*\vdz);
	    \tdplotdrawpolytopearc[thick,red]{\r}{}{};
	\tdplotdefinepoints(0,0,0)(\r/\vfl*\vfx, \r/\vfl*\vfy, \r/\vfl*\vfz)(\r/\vbl*\vbx, \r/\vbl*\vby, \r/\vbl*\vbz);
	    \tdplotdrawpolytopearc[thick,red]{\r}{}{};
	
\end{tikzpicture}
\caption{Decomposition of the cone $\cone$.}  \label{figure1} 
\end{figure}

\begin{observation}\label{remark:cones-cases}
  The elements $\bv_1, \bv_2, \bv_3$ are the completely fundamental
  elements of $\cone$, while $\bv_4 = \frac{1}{2}(\bv_2 + \bv_3)$ is
  merely fundamental; compare with~\cite[Chap.~I]{St96}.  A
  routine verification shows that
  \begin{align*}
    \cone_{134}&=\{(e_1,e_2,e_3)\in\cone \mid e_2\leq e_3 \},\\ 
    \cone_{145^+}&=\{(e_1,e_2,e_3)\in\cone \mid e_3<e_2\leq
                   e_3+\min\{e_1,e_3\}, \, 2e_2\leq e_1+2e_3\},\\ 
    \cone_{456^+}&=\{(e_1,e_2,e_3)\in\cone \mid  e_3<e_2\leq
                   e_3+\min\{e_1,e_3\}, \, e_1+2e_3< 2e_2 \},\\
    \cone_{462^+}&=\{(e_1,e_2,e_3)\in\cone \mid  e_1+e_3<e_2\leq 2e_3
                   \};
  \end{align*}
  hence the sets $\cone_{134}$, $\cone_{145^+}$, $\cone_{456^+}$,
  $\cone_{462^+}$ correspond precisely to Cases 1, 2a, 2b and 3 in
  Lemma~\ref{lemma:theta-tilde};
  compare with~\eqref{eqn:theta-zero-with-e_i}.
\end{observation}

The following decompositions are easily
verified: 
\begin{equation}\label{eq:decompositions-of-cones}
  \begin{array}{lcl}
    \cone = \cone_{134} \cupdot \cone_{145^+}\cupdot
    \cone_{456^+}\cupdot \cone_{462^+}, %
    &&
       \cone^{\shift}=\cone_{134}^{\shift} \cupdot
       \cone_{145^+}^{\shift} \cupdot 
       \cone_{456^+}^{\shift} \cupdot
       \cone_{462^+}^{\shift},\\ 
    \cone_{234} = \cone_{34}\cupdot \cone_{42^+}, %
    &\phantom{xx}&
                   \cone = \cone^{\shift} \cupdot \cone_{234}.
  \end{array}
\end{equation}
For convenience, for a subset $\cone_{ijk}\subseteq \cone$ write
$Z_{ijk}(s)=\sum_{\mathbf{e} \,\in \cone_{ijk}} X_1^{e_1}X_2^{e_2}X_3^{e_3} \,
\tilde\theta(\xi_\mathbf{e}(\pi))$ and adopt a similar shorthand notation for
subsets of the form
$\cone_{ijk^+}, \cone_{ij}, \cone_{ij^+}$. From
\eqref{eqn-t3-second-interm-expr} and
\eqref{eq:decompositions-of-cones} we deduce that
\begin{equation}\label{eqn:decomp-into-cones-beginning}
  \begin{split}
    \mathcal{Z}_{\dot{\mathbf{H}},\rho, \theta, p}(s) %
    &= \sum_{\mathbf{e} \, \in\, \cone}X_1^{e_1}X_2^{e_2}X_3^{e_3} \,
    \tilde\theta(\xi_\mathbf{e}(\pi))+p^{-1} \sum_{\mathbf{e} \,
      \in\,  \cone^{\shift}}X_1^{e_1}X_2^{e_2}X_3^{e_3}\tilde\theta(\xi_\mathbf{e}(\pi))\\
    &=(1+p^{-1}) \sum_{\mathbf{e} \, \in\,
      \cone}X_1^{e_1}X_2^{e_2}X_3^{e_3} \,
    \tilde\theta(\xi_\mathbf{e}(\pi))-p^{-1}Z_{234}(s) \\
    & =(1+p^{-1}) \big(
    Z_{134}(s)+Z_{145^+}(s)+Z_{456^+}(s)+Z_{462^+}(s) \big) - p^{-1}
    \big( Z_{34}(s)+Z_{42^+}(s) \big).
  \end{split}
\end{equation}

\begin{lemma}\label{lemma-computations-subcones}
  Referring to Definition~\ref{def-f}, we have
  \begin{align*}
    Z_{134}(s) %
    & =\frac{1}{1-p^3X_3}\sum_{i=0}^\infty (pX_1)^{i}F_{i,0}(X_2X_3),\\
    Z_{145^+}(s) %
    & = \frac{p^3X_1^2X_2^2X_3}{1-p^3X_1^2X_2^2X_3}\sum_{i=0}^\infty
      (pX_1)^{i}F_{i,0}(X_2X_3), \\
    Z_{456^+}(s) %
    & = \frac{1}{1-p^3X_1^2X_2^2X_3}\sum_{i=1}^\infty
      (p^{-1}X_1X_2)^{i}F_{0,i}^*(X_2X_3), \\
    Z_{462^+}(s) %
    & = \frac{X_2^2X_3}{1-X_2^2X_3}\sum_{i=0}^\infty
      (p^{-1}X_1X_2)^{i} \, F^*_{0,i}(X_2X_3), \\ 
    Z_{34}(s) %
    & = \frac{1}{1-p^3X_3}F_{0,0}(X_2X_3),\\
    Z_{42^+}(s) %
    & = \frac{X_2^2X_3}{1-X_2^2X_3} F_{0,0}(X_2X_3).
  \end{align*}
\end{lemma}

\begin{proof}
  The description appearing immediately after
  \eqref{eqn:theta-zero-with-e_i} provides explicit formulae for
  $\tilde\theta(\xi_\mathbf{e}(\pi))$ in each of the Cases 1, 2a, 2b
  and 3 which, by Remark~\ref{remark:cones-cases}, correspond to the
  subcones $\cone_{134}, \cone_{145^+}, \cone_{456^+}$ and
  $\cone_{462^+}$ respectively.  The sets $\cone_{34}$ and
  $\cone_{42^+}$ correspond to parts of Cases 1 and 3 respectively.
  \showproofs{The calculations are all similar; we show one of
    them.}{\color{purple}}%
  \showproofs{}{\\ \textbf{Case 1: }}%
  ~Elements of $\mathbf{e} \in \cone_{134}$ can be expressed in the
  form $\mathbf{e} = r_1\bv_1+r_3\bv_3+r_4\bv_4$, where
  $r_1, r_3, r_4\in \N_0$, so that
  \[
    \mathbf{e} = (e_1,e_2,e_3) =
    r_1\bv_1+r_3\bv_3+r_4\bv_4=(r_1,r_4,r_3+r_4).
  \]
  From this we obtain
  \begin{align*}
    Z_{134}(s) %
    &=\sum_{\mathbf{e} \in \cone_{134}} X_1^{e_1} X_2^{e_2} X_3^{e_3} \,
      \tilde\theta(\xi_\mathbf{e}(\pi))\\ 
    &=\sum_{\mathbf{e} \in \cone_{134}} X_1^{e_1} X_2^{e_2} X_3^{e_3} \,
      p^{e_1-3e_2+3e_3} \, f(e_1,0,e_2)\\ 
    &=\sum_{r_1,r_3,r_4\geq 0} (pX_1)^{r_1} \, (p^{3}X_3)^{r_3} \,
      (X_2X_3)^{r_4} \, f(r_1,0,r_4)\\
    &=\frac{1}{1-p^3X_3} \sum_{i=0}^\infty (pX_1)^{i} \,
      F_{i,0}(X_2X_3). \showproofs{\qedhere}{}
  \end{align*}
  \showproofs{}{Setting $r_1 = 0$ in the above
    computation, we obtain 
    \begin{equation*}
      Z_{34}(s) 
       =\sum_{r_3,r_4\geq 0} (p^3X_3)^{r_3} \,
	(X_2X_3)^{r_4} \, f(0,0,r_4)
      =\frac{1}{1-p^3X_3} F_{0,0}(X_2X_3). 
    \end{equation*}
    
    \noindent\textbf{Case 2a:}
    Elements of $\e \in \cone_{145^+}$ can be expressed in the form
    $\e = r_1\bv_1+r_4\bv_4+r_5\bv_5$, where $r_1, r_4 \in \N_0$ and
    $r_5\in \N$, so that
    \[
      \e = (e_1,e_2,e_3) =
      r_1\bv_1+r_4\bv_4+r_5\bv_5=(r_1+2r_5,r_4+2r_5,r_4+r_5).
    \]
    From this we obtain
    \begin{align*}
      Z_{145^+}(s) %
      &=\sum_{\mathbf{e} \in \cone_{145^+}} X_1^{e_1} X_2^{e_2} X_3^{e_3} \,
        \tilde\theta(\xi_\mathbf{e}(\pi))\\ 
      &=\sum_{\mathbf{e} \in \cone_{145^+}} X_1^{e_1} X_2^{e_2} X_3^{e_3} \,
	p^{e_1+e_2-e_3} \, f(e_1-2e_2+2e_3,0,-e_2+2e_3)\\ 
      &=\sum_{r_1,r_4\geq 0;\ r_5>0} (pX_1)^{r_1} \, (X_2X_3)^{r_4} \,
	(p^3X_1^2X_2^2X_3)^{r_5} \, f(r_1,0,r_4)\\
      &=\frac{p^3X_1^2X_2^2X_3}{1-p^3X_1^2X_2^2X_3}\sum_{r_1,r_4\geq
        0} (pX_1)^{r_1} \, (X_2X_3)^{r_4} \,  \, f(r_1,0,r_4)\\
      &=\frac{p^3X_1^2X_2^2X_3}{1-p^3X_1^2X_2^2X_3}\sum_{i= 0}^\infty
        (pX_1)^{i}  \, F_{i, 0}(X_2X_3).
    \end{align*}

    \noindent \textbf{Case 2b:} Elements of $\e \in \cone_{456^+}$ can
    be expressed in the form $\e = r_4\bv_4+r_5\bv_5+r_6\bv_6$, where
    $r_4, r_5 \in \N_0$ and $r_6\in \N$, so that
    \[
      \e = (e_1,e_2,e_3) =
      r_4\bv_4+r_5\bv_5+r_6\bv_6=(2r_5+r_6,r_4+2r_5+2r_6,r_4+r_5+r_6).
    \]
    From this we obtain
    \begin{align*}
      Z_{456^+}(s) %
      &=\sum_{\mathbf{e} \in \cone_{456^+}} X_1^{e_1} X_2^{e_2} X_3^{e_3} \,
        \tilde\theta(\xi_\mathbf{e}(\pi))\\ 
      &=\sum_{\mathbf{e} \in \cone_{456^+}} X_1^{e_1} X_2^{e_2} X_3^{e_3} \,
	p^{4e_1-5e_2+5e_3} \, f(0,-e_1+2e_2-2e_3,-e_1+e_2)\\ 
      &=\sum_{r_4,r_5\geq 0;\ r_6>0} (X_2X_3)^{r_4} \, (p^3X_1^2X_2^2X_3)^{r_5} \,
	(p^{-1}X_1X_2^2X_3)^{r_6} \, f(0,r_6,r_4+r_6)\\
      &=\frac{1}{1-p^3X_1^2X_2^2X_3}\sum_{r_4\geq 0;\ r_6>0}   \,
	(p^{-1}X_1X_2)^{r_6}\, (X_2X_3)^{r_4+r_6} \, f(0,r_6,r_4+r_6)\\
      &=\frac{1}{1-p^3X_1^2X_2^2X_3}\sum_{i=1}^\infty   \,
        (p^{-1}X_1X_2)^{i}\, F^*_{0,i}(X_2X_3).
    \end{align*}

    \noindent\textbf{Case 3:}  
    Elements of $\e \in \cone_{462^+}$ can be expressed in the form
    $\e = r_4\bv_4+r_6\bv_6+r_2\bv_2$, where $r_4, r_6 \in \N_0$ and
    $r_2\in \N$, so that
    \[
      \e = (e_1,e_2,e_3) =
      r_4\bv_4+r_6\bv_6+r_2\bv_2=(r_6,r_4+2r_6+2r_2,r_4+r_6+r_2).
    \]
    From this we obtain
    \begin{align*}
      Z_{462^+}(s) %
      &=\sum_{\mathbf{e} \in \cone_{462^+}} X_1^{e_1} X_2^{e_2} X_3^{e_3} \,
        \tilde\theta(\xi_\mathbf{e}(\pi))\\ 
      &=\sum_{\mathbf{e} \in \cone_{462^+}} X_1^{e_1} X_2^{e_2} X_3^{e_3} \,
	p^{-e_1} \, f(0,e_1,e_1-e_2+2e_3)\\ 
      &=\sum_{r_4,r_6\geq 0;\ r_2>0} (X_2X_3)^{r_4} \, (p^{-1}X_1X_2^2X_3)^{r_6} \,
	(X_2^2X_3)^{r_2} \, f(0,r_6,r_4+r_6)\\
      &=\frac{X_2^2X_3}{1-X_2^2X_3}\sum_{r_4,r_6\geq 0} (X_2X_3)^{r_4} \, (p^{-1}X_1X_2^2X_3)^{r_6} \, f(0,r_6,r_4+r_6)\\
      &=\frac{X_2^2X_3}{1-X_2^2X_3}\sum_{r_4,r_6\geq 0}
        (p^{-1}X_1X_2)^{r_6} \,(X_2X_3)^{r_4+r_6}\, f(0,r_6,r_4+r_6)\\ 
      &=\frac{X_2^2X_3}{1-X_2^2X_3} \sum_{i=0}^\infty
        (p^{-1}X_1X_2)^{i} \, F^*_{0,i}(X_2X_3).	 
    \end{align*}
    Setting $r_6=0$ in the above computation, we obtain
    \begin{equation*}
      Z_{42^+}(s) =\sum_{r_4\geq 0;\ r_2>0} (X_2X_3)^{r_4} \, 
      (X_2^2X_3)^{r_2} \, f(0,0,r_4)
      =\frac{X_2^2X_3}{1-X_2^2X_3}F_{0,0}(X_2X_3).
      \showproofs{}{\qedhere}			
    \end{equation*}}
\end{proof}

\showproofs{Explicit formulae for the expressions in
  Lemma~\ref{lemma-computations-subcones} can now be obtained via
  Proposition~\ref{prop:F-fns} and Lemma~\ref{lemma:Fstar-fns-new}.}{
  \color{purple} By Lemma~\ref{lemma-computations-subcones} and
  Proposition~\ref{prop:F-fns}, we have
  \begin{align*}&Z_{134}(s)+Z_{145^+}(s)\\
                &=\frac{1-p^6X_1^2X_2^2X_3^2}{(1-p^3X_3)(1-p^3X_1^2X_2^2X_3)}\sum_{i=0}^\infty
                  (pX_1)^{i}F_{i,0}(X_2X_3)\\
                &=\frac{1-p^6X_1^2X_2^2X_3^2}{(1-p^3X_3)(1-p^3X_1^2X_2^2X_3)}\sum_{i=0}^\infty
                  (pX_1)^i \left[p^{2i}(X_2X_3)^i
                  F_{0,0}(X_2X_3)+\frac{(1-pX_2X_3)(1-p^{2i}(X_2X_3)^i)}{(1-p^2X_2X_3)^2}\right]\\
                &=\frac{1-p^6X_1^2X_2^2X_3^2}{(1-p^3X_3)(1-p^3X_1^2X_2^2X_3)} \\
                &\quad\quad \cdot\left[
                  \frac{(1-pX_2X_3)(1-p^3X_1X_2X_3)+(1-pX_1)\left(-1+pX_2X_3+F_{0,0}(X_2X_3)(1-p^2X_2X_3)^2\right)
                  }{ (1-pX_1)(1-p^2X_2X_3)^2(1-p^3X_1X_2X_3) }
                  \right].
  \end{align*}
  By Lemmata~\ref{lemma-computations-subcones} and
  \ref{lemma:Fstar-fns-new}, we have
  \begin{align*}
    &Z_{456^+}+Z_{462^+}\\
    &=\frac{1 - X_2^2X_3 + X_2^2X_3 - p^3X_1^2X_2^4X_3^2	
      }{(1-X_2^2X_3)(1-p^3X_1^2X_2^2X_3)}
      \sum_{i=0}^\infty 
      (p^{-1}X_1X_2)^i F^*_{0,i}(X_2X_3)\\
    &\quad\quad - 
      \frac{1}{1-p^3X_1^2X_2^2X_3}F_{0,0}(X_2X_3)\\
    &=\left(1 -p^3X_1^2X_2^4X_3^2\right)\\
    &\quad \quad \cdot \frac{(1-pX_2X_3)(1+pX_2X_3+pX_1X_2^2X_3(1+p^2X_2X_3))}{(1-X_2^2X_3)(1-p^3X_1^2X_2^2X_3)(1-p^3X_1^2X_2^4X_3^2)(1-p^2X_2X_3)(1-p^3X_2^2X_3^2)}\\
    &\quad\quad - 
      \frac{1}{1-p^3X_1^2X_2^2X_3}F_{0,0}(X_2X_3)
  \end{align*}
  By Lemma~\ref{lemma-computations-subcones}, we have
  \begin{align*}
Z_{34}(s)+Z_{42^+}(s)=F_{0,0}(X_2X_3)\left(
\frac{1}{1-p^3X_3}+\frac{X_2^2X_3}{1-X_2^2X_3}
\right)=
\frac{(1-p^3X_2^2X_3^2)F_{0,0}(X_2X_3)}{(1-p^3X_3)(1-X_2^2X_3)}.
\end{align*}
Note that by Proposition~\ref{prop:F-fns} and Lemma~\ref{lemma:Fstar-fns-new},
\begin{align*}
	F_{0,0}(X_2X_3)&=\frac{1-p^2X_2^2X_3^2}{(1-p^2X_2X_3)(1-p^3X_2^2X_3^2)},\\	
	F^*_{0,1}(X_2X_3)&=p^2X_2X_3\frac{(1-pX_2X_3)(1+p^2X_2X_3)}{(1-p^2X_2X_3)(1-p^3X_2^2X_3^2)}.
\end{align*}
Substituting the above expressions into
\eqref{eqn:decomp-into-cones-beginning}, we obtain }
\showproofs{Substituting these into \eqref{eqn:decomp-into-cones-beginning} yields}{}
\begin{multline} \label{equ:X1X2X3-expression}
  \mathcal{Z}_{\dot{\mathbf{H}},\rho, \theta, p}(s) = \\
  \frac{(1-pX_2X_3)(-p^4X_1^3X_2^3X_3^2-p^3X_1^3X_2^2X_3-p^4X_1^2X_2^3X_3^2-pX_1^2X_2^2X_3+p^3X_1X_2X_3+X_1+pX_2X_3+1)}{(1-pX_1)(1-p^3X_3)(1-X_2^2X_3)(1-p^2X_2X_3)(1-p^3X_1^2X_2^2X_3)}.
\end{multline}
Recalling that $X_1=p^{14-5s}, X_2=p^{2-s}$, $X_3=p^{26-9s}$ we obtain
\begin{align*}
  \zeta^\wedge_{\Gamma_{t^3},p}(s) %
  &=\frac{(1-p^{29-10s})(-p^{104-36s}-p^{90-31s}-p^{75-26s}-p^{59-21s}+p^{45-15s}+p^{29-10s}+p^{14-5s}+1)}{(1-p^{15-5s})(1-p^{29-9s})(1-p^{30-11s})(1-p^{30-10s})(1-p^{61-21s})}\\ 
  &=
    \frac{(1-p^{29-10s})(-p^{89-31s}-p^{75-26s}+p^{74-26s}-p^{59-21s}+p^{30-10s}-p^{15-5s}+p^{14-5s}+1)}{(1-p^{15-5s})^2(1-p^{29-9s})(1-p^{30-11s})(1-p^{61-21s})}, 
\end{align*}
proving Theorem~\ref{x3example}.


\section{Meromorphic continuation for the pro-isomorphic zeta function
  of $\Gamma_{t^3}$} \label{sec:mero-cont} In this section we consider
the pro-isomorphic zeta function of the $D^*$-group
$\Gamma = \Gamma_{t^3}$ of Hirsch length~$8$, defined
in~\eqref{equ:def-Gamma-tm}.  Our task is to deduce the assertions
about $\zeta^{\wedge}_\Gamma(s)$ in Corollary~\ref{cor:t3} from the
Euler product decomposition~\eqref{equ:euler-decomp} and the explicit
description of the local zeta functions in Theorem~\ref{x3example}.
The main step is to establish that the line
$\{ s \in \C \mid \operatorname{Re}(s)=3 \}$ is a natural
boundary for the meromorphic continuation of
$\zeta^{\wedge}_\Gamma(s)$.  We follow the strategy laid out
in~\cite[Chap.~5]{dSWo08} and use a compatible notation; in the
terminology of~\cite{dSWo08}, we are dealing with a Type~II situation,
which requires approximations up to terms of degree~$3$, as we shall
see.

Theorem~\ref{x3example} shows that
\begin{equation} \label{equ:psi-appears} \zeta^{\wedge}_\Gamma(s) =
  \frac{\zeta(5s-15)^2 \zeta(9s-29) \zeta(11s-30) \zeta(21s-61)
  }{\zeta(10s-29)} \; \psi(s),
\end{equation}
where $\zeta(s)$ denotes the Riemann zeta function and
\begin{equation} \label{equ:psi-prod}
   \psi(s) =
   \prod_p \widetilde{W}(p,p^{-s})
 \end{equation}
 for
 \[
   \widetilde{W}(X,Y) = 1 + X^{14}Y^5 - X^{15}Y^5 + X^{30}Y^{10} -
   X^{59}Y^{21} + X^{74}Y^{26} - X^{75}Y^{26} -
   X^{89}Y^{31},
 \]
 as in the statement of Corollary~\ref{cor:t3}.  It is routine to
 check that the infinite product in~\eqref{equ:psi-prod} converges
 absolutely for all $s \in \C$ with
\[
  \operatorname{Re}(s) > \max \big\{ 15/5, 16/5, 31/10, 60/21, 75/26,
  76/26, 90/31 \big\} = 16/5
\]
and yields a holomorphic function
on~$\{ s \in \C \mid \operatorname{Re}(s) > 16/5 \}$.  In
passing, we observe that the abscissa of convergence of the Dirichlet
generating series $\zeta^{\wedge}_\Gamma(s)$, which has non-negative
coefficients, can be detected by looking for the right-most
singularity on the real line; from \eqref{equ:psi-appears} we see that
this singularity lies at $s = 30/9 = 10/3$ and yields a simple pole.

Next we show that the function $\psi(s)$, and thus
$\zeta^{\wedge}_\Gamma(s)$, can be meromorphically continued further
to the right-half plane
$\mathcal{H} = \{ s \in \C \mid \operatorname{Re}(s) > 3
\}$.  Indeed, the cyclotomic polynomial $1-t+t^2$ does not vanish at
$t = p^{15-5s}$ for $s \in \mathcal{H}$, because
$\lvert p^{15-5s} \rvert < 1$.  We consider
\begin{equation} \label{equ:psi-tilde}
  \widetilde{\psi}(s) =  \prod_p
  \frac{\widetilde{W}(p,p^{-s})}{1-p^{15-5s}+p^{30-10s}} = \prod_p
  \left( 1 + \frac{p^{14-5s} - p^{59-21s} + p^{74-26s} - p^{75-26s} -
      p^{89 -31s}}{1-p^{15-5s}+p^{30-10s}} \right);
\end{equation}
this infinite product converges absolutely and yields a holomorphic
function for $s \in \mathcal{H}$, because
\[
  \max \big\{ 15/5, 60/21, 75/26, 76/26, 90/31 \big\} = 3.
\]
As $1-t+t^2 = (1-t^6)(1-t)(1-t^2)^{-1}(1-t^3)^{-1}$, we see that
\[
  \psi(s) = \frac{\zeta(10s-30) \zeta(15s-45)}{\zeta(30s-90)
    \zeta(5s-15)} \; \widetilde{\psi}(s), \qquad \text{for
    $s \in \mathcal{H}$,}
\]
yields the desired meromorphic continuation.  Furthermore, using a
Tauberian theorem~\cite[Thm.~4.20]{dSGr00} as in the proof of
Corollary~\ref{cor:t2}, we obtain the description of the asymptotic
growth of pro-isomorphic subgroups in $\Gamma_{t^3}$ as recorded in
Remark~\ref{rem:asymp-growth-t3}.

It remains to show that the line
$\mathcal{L} = \{ s \in \C \mid \operatorname{Re}(s) = 3 \}$
is a natural boundary for $\psi(s)$; in view
of~\eqref{equ:psi-appears}, this implies that $\mathcal{L}$ is
also a natural boundary for $\zeta^{\wedge}_\Gamma(s)$ and
Corollary~\ref{cor:t3} follows.  The strategy is to show that each
point $s \in \mathcal{L}$ is a limit point of zeros of the
meromorphic function $\psi(s)$, defined on~$\mathcal{H}$; since
poles and zeros of the Riemann zeta function are isolated, it suffices
to show that each $s \in \mathcal{L}$ is a limit point of zeros
of the holomorphic function $\widetilde{\psi}(s)$, defined
on~$\mathcal{H}$.  Recall from~\eqref{equ:psi-tilde} that
$\widetilde{\psi}(s)$ is given as an infinite product, indexed by~$p$;
thus $\widetilde{\psi}(s)$ vanishes, for any given
$s \in \mathcal{H}$, if and only if $\widetilde{W}(p,p^{-s})$
vanishes for at least one prime~$p$.

This leads us to study the zeros of the polynomial
\[
  F(V,U) = 1 + (V-1)U^5 + U^{10} - V^4 U^{21} + \left( V^4-V^3 \right)
  U^{26} - V^4 U^{31} \in \Z[V][U].
\]
Observe that $F(X^{-1},X^3Y) = \widetilde{W}(X,Y)$; we will be
interested in evaluating~$F$ at $V = p^{-1} \to 0$, as the prime $p$
tends to infinity, and $U = p^{3-s}$, for suitable
$s \in \mathcal{H}$ depending on~$p$.  We see that
\[
  F(0,U) = 1 - U^5 + U^{10}
\]
is a product of the $6$th and the $30$th cyclotomic polynomial.  We
fix the primitive $6$th root of unity
$\lambda = \exp(\pi i / 3) = (1+\sqrt{3}\,i)/2$ so that $\lambda$ is a
root of $1-t+t^2$, and we fix the primitive $30$th root of unity
$\omega = \exp(\pi i / 15)$ so that $\omega$ is a root of $F(0,U)$.
By the Holomorphic Implicit Function Theorem, there is a holomorphic
function $u = u(v)$, defined in a small complex neighbourhood of
$v=0$, such that $u(0) = \omega$ and $F\big(v,u(v)\big)=0$;
furthermore, being analytic, this function admits a local
representation as a power series
\[
  u(v) = \omega \; \left( 1 + a_1 v + a_2 v^2 + a_3 v^3 + \ldots \right)
\]
in $v$ with uniquely determined complex coefficients.  A routine power
series calculation and comparison of coefficients yield
\[
  a_1 = \tfrac{1}{15} (2\lambda -1), \qquad
  a_2  = \tfrac{1}{15^2} (1-5\lambda), \qquad
  a_3  = \tfrac{1}{15^3} \left(-17-450 \omega + 49 \omega^5 + 225
    \omega^6 \right).
\]
Writing $u(v) = p^{3-s}$ and $v = p^{-1}$, for sufficiently large $p$,
we solve for $s \in \C$ to obtain a set 
\[
  \mathcal{N}_p = \mathcal{H} \cap \bigg\{ 3 -
  \underbrace{\frac{\log(\omega)}{\log(p)}}_{\in \R i} -
  \underbrace{\frac{\log(1 + a_1 p^{-1} + a_2 p^{-2} + a_3 p^{-3} +
      \ldots)}{\log(p)}}_{(\ast)} - \underbrace{\frac{2\pi k}{\log(p)}
    i}_{\in \R i} \mid k \in \Z \bigg\}
\]
of zeros of $\psi(s)$, where $k$ is a parameter that we can use, for
increasing~$p$, to approximate any given point on the
line~$\mathcal{L}$ to any required degree.  However, we still need to
verify that, for sufficiently large~$p$, the real part of the
numerator in~$(\ast)$ is negative so that the resulting candidate zero
lies in~$\mathcal{H}$, as required.

Using the logarithm series $\log(1 + t) = t - \frac{1}{2} t^2 +
\frac{1}{3} t^3 - \ldots$ for small $t = a_1 p^{-1} + a_2 p^{-2} + a_3
p^{-3} + \ldots$, we see that the relevant numerator in~$(\ast)$ is
\begin{equation} \label{equ:approx-mod-Omega}
  \underbrace{\phantom{\big(}\!\!a_1\phantom{\big)}\!\!\!}_{\in
    \R i} p^{-1} + \underbrace{\left( a_2 - \tfrac{1}{2}
      a_1^{\, 2} \right)}_{\in \R i} p^{-2} +
  \underbrace{\left( a_3 - a_1 a_2 + \tfrac{1}{3} a_1^{\, 3}
    \right)}_{\text{has negative real part}} p^{-3} + \Omega(p^{-1}),
\end{equation}
where $\Omega(v)$ is a complex power series in $v$ starting with $v^4$
or some higher term.  Indeed, short calculations yield
\[
    a_1 = \tfrac{1}{15} (2\lambda -1) = \tfrac{1}{15} \sqrt{3}\, i \in
        \R i \qquad \text{and} \qquad
  a_2 - \tfrac{1}{2}  a_1^{\, 2} = \tfrac{1}{15^2} \left(
                                   \tfrac{5}{2} - 5\lambda \right) =
                                   \tfrac{-1}{90} \sqrt{3}\, i \in
                                   \R i.
\]
Furthermore, a slightly longer, but routine calculation gives
\[
  a_3 - a_1 a_2 + \tfrac{1}{3} a_1^{\, 3} = \tfrac{1}{15^3} \left( -25
    + 50 \lambda + 15^2 (\lambda-2) \omega \right) = \tfrac{1}{135}
  \left( -1 - 18\omega + 2\omega^5 + 9\omega^6 \right)
\]
and, since
$\operatorname{Re}(\omega) > \operatorname{Re}(\omega^5) >
\operatorname{Re}(\omega^6)$, we deduce that
\[
  \operatorname{Re} (a_3 - a_1 a_2 + \tfrac{1}{3} a_1^{\, 3} ) <
  0
\]
as asserted.  For sufficiently large~$p$, the contribution of
$\Omega(p^{-1})$ in~\eqref{equ:approx-mod-Omega} is much smaller than
the $p^{-3}$-term; hence $\mathcal{N}_p$ supplies the required zeros
of~$\psi(s)$.
  
  
\section{Base extensions}\label{sec:base-ext}

Following a suggestion of the referee, we extend in this section our
results for the $\Q$-indecomposable $D^*$-groups $\Gamma_{t^2}$ and
$\Gamma_{t^3}$ to two infinite families of class-two nilpotent groups
that result naturally from the initial groups via `base extensions' of
the corresponding Lie lattices; for completeness, we
  also discuss what happens if we start with the decomposable
  $D^*$-group $\Gamma_t$.  The outcome illustrates that the
investigation in~\cite{BGS22}, which was carried out partly after,
partly in parallel to our original work, has an impact in the
situation that we consider in this paper.  We exercise some care not
to exclude any primes; this allows us to get explicit results in the
global setting.  In a nutshell we will see that the calculations
carried out in Sections~\ref{section:x^2} and~\ref{section:x^3}
require only mild modifications, once the relevant algebraic
automorphism groups are understood.  In particular, we establish
Theorem~\ref{thm:t2-base-change}.

We briefly set up the scene.  Let $L$ be a nilpotent $\Z$-Lie
lattice; our main interest will be in $L = L_{t^m}$, the Lie lattice
associated to the nilpotent group $\Gamma_{t^m}$ with
presentation~\eqref{equ:def-Gamma-tm}, with an extra focus on
$m \in \{2,3\}$.  We consider a number field $k$ of absolute degree
$d = [k:\Q]$, with ring of integers~$\mathfrak{o}$.  By
extension of scalars from $\Z$ to $\mathfrak{o}$ and
restriction of scalars back to $\Z$, we obtain a
$\Z$-Lie lattice
$\widetilde{L} = {_{\Z,\mathfrak{o}} L}$ of $\Z$-rank
$\dim_{\Z}(\widetilde{L}) = d \dim_{\Z}(L)$.  Clearly,
$\widetilde{L}$ is nilpotent, of the same class as~$L$.

Automorphisms of $L$ induce in a natural way automorphisms of
$\widetilde{L}$, but, in general, the automorphism group of
$\widetilde{L}$ may turn out considerably more `complex' than that
of~$L$.  Consequently, the pro-isomorphic zeta functions of
$\widetilde{L}$ and of $L$ may bear little resemblance to one another.
Our aim in this section is to show that Lie lattices of the form
$L = L_{t^m}$, for $m \in \N_{\ge 2}$, are sufficiently
`rigid' so that $\mathbf{Aut}(\widetilde{L})$ is strongly linked to
$\mathbf{Aut}(L)$, in an appropriate local sense.  For
$m \in \{2,3\}$, this allows us to determine the local pro-isomorphic
zeta
functions~$\zeta^\wedge_{\widetilde{L},p}(s) =
\zeta_{\widetilde{L}_p}^{\mathrm{iso}}(s)$ for all primes~$p$ and, via
the Euler product~\eqref{equ:euler-decomp}, we deduce analytic
properties of the pro-isomorphic zeta
function~$\zeta^\wedge_{\widetilde{\Gamma}}(s)$ of the class-two
nilpotent group $\widetilde{\Gamma}$ associated to~$\widetilde{L}$;
compare with
Section~\ref{section:class-2-correspondence}.  The Lie
  lattice $L_t$ is not quite `rigid', but a slight modification of the
  approach in \cite{BGS22} allows us to bypass the problem and we
  obtain the local and global pro-isomorphic zeta functions also in
  this basic case.


\subsection{Local rigidity of the Lie lattices $L_{t^m}$ for
  $m\geq 2$.} \label{sec:local-rgidity}

As above, let $\widetilde{L} = {_{\Z,\mathfrak{o}} L}$ denote
the $\Z$-Lie lattice associated, via `base extension', to a
$\Z$-Lie lattice $L$ and a number field $k$ with ring of
integers~$\mathfrak{o}$.  Fix a rational prime~$p$, and recall that
there are finitely many non-archimedean primes
$\mathfrak{p} \in \operatorname{Spec}(\mathfrak{o})$ dividing~$p$.  It
is well known that there is a natural ring isomorphism
${_{\Z_p} \mathfrak{o}} = \Z_p \otimes_\Z
\mathfrak{o} \cong \prod_{\mathfrak{p} \mid p}
\mathfrak{o}_\mathfrak{p}$, where $\mathfrak{o}_\mathfrak{p}$ denotes
the completion of $\mathfrak{o}$ at the prime~$\mathfrak{p}$.  From
this one sees that the $\Z_p$-Lie lattice
$\widetilde{L}_p = \Z_p \otimes_\Z \widetilde{L}$,
relevant to our investigation, is isomorphic to
$\bigoplus_{\mathfrak{p} \mid p} \widetilde{L}_\mathfrak{p}$, where
$\widetilde{L}_\mathfrak{p} = {_{\Z_p,
    \mathfrak{o}_\mathfrak{p}} L_p}$ denotes the $\Z_p$-Lie
lattice that results from $L_p = \Z_p \otimes_\Z L$
via extension of scalars to the complete valuation ring
$\mathfrak{o}_\mathfrak{p}$ and restriction back to~$\Z_p$.

This prompts us to consider the $\Z_p$-Lie lattice
$\widetilde{L}_\wp = {_{\Z_p, {\scriptscriptstyle
      \mathcal{O}}} L_p}$, for any given finite extension
$\mathcal{F}$ of~$\Q_p$, with valuation ring
${\scriptstyle \mathcal{O}}$ and valuation ideal~$\wp$.  Write
$\widetilde{\G}_\wp = \mathbf{Aut}(\widetilde{L}_\wp)$ and
$\G_p = \mathbf{Aut}(L_p)$ for the algebraic automorphism groups of
the $\Z_p$-Lie lattices $\widetilde{L}_\wp$ and~$L_p$.  Here
$\G_p$ is simply the $\Z_p$-group scheme that results from the
algebraic automorphism group $\G = \mathbf{Aut}(L)$ of the original
$\Z$-Lie lattice via base change: any $\Z$-basis
$\mathcal{S}$ of~$L$ naturally identifies a $\Z_p$-basis
of~$L_p$, and via $\mathcal{S}$ we realise $\G_p \le \GL_n$ as an
affine $\Z_p$-group scheme, for
$n = \dim_\Z(L) = \dim_{\Z_p}(L_p)$.  In the following
we write $\G$ in place of $\G_p$, when the base ring is insignificant.
Moreover, tensoring $\mathcal{S}$ with a $\Z_p$-basis of
${\scriptstyle \mathcal{O}}$, we obtain a $\Z_p$-basis
$\widetilde{\mathcal{S}}$ of $\widetilde{L}_\wp$, which allows us to
realise $\widetilde{\G}_\wp \le \GL_{nd}$ as an affine
$\Z_p$-group scheme, where
$d = \dim_{\Z_p}({\scriptstyle \mathcal{O}}) = [\mathcal{F} :
\Q_p]$.  Our explicit construction yields, in particular,
\[
  \G({\scriptstyle \mathcal{O}}) \cong \Aut({_{\scriptscriptstyle
      \mathcal{O}} L_p}) \le \Aut({_{\Z_p,{\scriptscriptstyle
        \mathcal{O}}} L_p}) \cong \widetilde{\G}_\wp(\Z_p),
  \qquad \G(\mathcal{F}) \cong \Aut({_\mathcal{F}
    L_p}) \le \Aut({_{\Q_p,\mathcal{F}} L_p}) \cong
  \widetilde{\G}_\wp(\Q_p).
\]

Typically, these embeddings are proper, because $\Z_p$-linear
automorphisms are not necessarily ${\scriptstyle \mathcal{O}}$-linear.
Suppose that $L$, hence also $L_p$, is nilpotent of class~$2$.  In
this situation we can easily make out two types of automorphisms,
which could be used to fill this gap: central automorphisms and field
automorphisms.  More precisely, we set
\[
  \mathbf{J}_\wp = \mathsf{C}_{\widetilde{\G}_\wp}
  \big(\widetilde{L}_\wp / \mathrm{Z}(\widetilde{L}_\wp)\big) = \ker
  \big( \mathbf{Aut}(\widetilde{L}_\wp) \to \mathbf{Aut}
  \big(\widetilde{L}_\wp / \mathrm{Z}(\widetilde{L}_\wp) \big) \big)
  \trianglelefteq \widetilde{\G}_\wp,
\]
the affine $\Z_p$-group scheme which is the algebraic
centraliser of the $\Z_p$-module
$\widetilde{L}_\wp / \mathrm{Z}(\widetilde{L}_\wp)$.  For the concrete
realisation as a subgroup scheme in~$\GL_{nd}$, it is convenient to
choose the underlying $\Z_p$-basis $\mathcal{S}$ in such a way
that it includes a $\Z_p$-basis for $\mathrm{Z}(L_p)$; then
$\widetilde{\mathcal{S}}$ includes a $\Z_p$-basis
for~$\mathrm{Z}(\widetilde{L}_\wp)$ and $\mathbf{J}_\wp$ can be
defined rather directly.  In addition, we consider the algebraic
automorphism group of the extension $\mathcal{F} \,|\, \Q_p$
as a subgroup scheme of $\mathbf{Aut}(\widetilde{L}_\wp)$, that is the finite
group scheme
\[
  \mathbf{F}_\wp \le \widetilde{\G}_\wp \qquad \text{with} \qquad
  \mathbf{F}_\wp \cong \mathbf{Aut}({\scriptstyle \mathcal{O}} \,|\,
  \Z_p) \cong \mathbf{Aut}(\mathcal{F} \,|\, \Q_p)
\]
such that, in particular,
$\mathbf{F}_\wp(\Z_p) \cong \Aut({\scriptstyle \mathcal{O}}
\,|\, \Z_p) \cong \Aut(\mathcal{F} \,|\, \Q_p)$ acts
naturally via field automorphisms on the Lie
lattice~$\widetilde{L}_\wp$.  Furthermore, we observe that $\G$, now
regarded as an affine ${\scriptstyle \mathcal{O}}$-group scheme,
serves as the algebraic automorphism group
$\mathbf{Aut}({_{\scriptscriptstyle \mathcal{O}} L_p})$ of the
${\scriptstyle \mathcal{O}}$-Lie lattice
${_{\scriptscriptstyle \mathcal{O}} L_p}$; accordingly, the affine
$\Z_p$-group scheme
$\mathsf{Res}_{{\scriptscriptstyle \mathcal{O}} \,|\,
  \Z_p}(\G)$ which results via restriction of scalars can be
realised as a subgroup scheme of~$\widetilde{\G}_\wp$.  We are
interested in situations where the following `rigidity'
holds:
\begin{equation} \label{equ:JRF-equality} \big( \mathbf{J}_\wp \cdot
  \mathsf{Res}_{{\scriptscriptstyle \mathcal{O}} \,|\,
    \Z_p}(\G) \big) \rtimes \mathbf{F}_\wp =
  \widetilde{\G}_\wp \qquad \text{as $\Q_p$-defined algebraic
    subgroups of $\GL_{nd}$}.
\end{equation}
Actually, for us it suffices that the two group schemes yield the same
groups of $\Q_p$-rational points; this condition is slightly
weaker, but implies, for instance, that the two
$\Q_p$-algebraic groups have the same connected component.  In
down-to-earth terms we require that the $\mathcal{F}$- and thus also
$\Q_p$-Lie algebra
$\widetilde{\mathcal{L}} = {_\mathcal{F} L_p} = \mathcal{F}
\otimes_{\Z_p} L_p$ satisfies
\begin{equation} \label{equ:down-to-earth}
  \left( \mathrm{C}_{\Aut_{\Q_p}(\widetilde{\mathcal{L}})}
  \big( \widetilde{\mathcal{L}} / \mathrm{Z}(\widetilde{\mathcal{L}}) \big) \;
  \Aut_\mathcal{F}(\widetilde{\mathcal{L}}) \right) \rtimes
  \Aut(\mathcal{F} \,|\, \Q_p) =
  \Aut_{\Q_p}(\widetilde{\mathcal{L}}). 
\end{equation}

In~\cite{BGS22}, Berman, Glazer and Schein extend results of
Segal~\cite{Se89} for algebraic automorphism groups of certain Lie
algebras, with a view toward studying pro-isomorphic zeta functions
under `base extensions'.  In particular, they formulate sufficient
conditions under which \eqref{equ:JRF-equality} holds true; see
\cite[Thm.~3.9]{BGS22}.  For the discussion at hand, a special and
thus simpler version of their criterion is sufficient.  We say that
the $\Q_p$-Lie
algebra~${_{\Q_p} L} \cong {_{\Q_p} L_p}$ is
\emph{absolutely indecomposable} if, for every finite extension
$\mathcal{F}$ of~$\Q_p$, the $\mathcal{F}$-Lie algebra
${_\mathcal{F} L} \cong {_\mathcal{F} L_p}$ is indecomposable.  We
make use of the following special instance of~\cite[Thm.~3.9]{BGS22}.

\begin{lemma} \label{lem:BGS-lemma}
  Let $L$ be a class-two nilpotent $\Z$-Lie lattice, as above,
  and such that $[L,L] = \mathrm{Z}(L)$.  Let $p$ be a prime such
  that the $\Q_p$-Lie algebra
  $\mathcal{L} = {_{\Q_p} L_p}$ is absolutely indecomposable and
  generated by
  \[
    \mathcal{Y} = \big\{ w \in \mathcal{L} \smallsetminus
    \mathrm{Z}(\mathcal{L}) \mid \mathrm{C}_\mathcal{L} \big(
    \mathrm{C}_\mathcal{L}(w) \big) = \Q_p w +
    \mathrm{Z}(\mathcal{L}) \big\}.
  \]
  Then \eqref{equ:JRF-equality} holds, for every finite extension
  $\mathcal{F}$ of~$\Q_p$, with valuation
  ring~${\scriptstyle \mathcal{O}}$ and valuation ideal~$\wp$.
\end{lemma}

Next we consider the $\Z$-Lie lattices $L_{t^m}$,
$m \in \N$, associated to the $D^*$-groups $\Gamma_{t^m}$
with presentation~\eqref{equ:def-Gamma-tm}.  This means that
$L_{t^m}$ has $\Z$-rank $2m+2$ and admits the presentation
\begin{multline} \label{equ:presentation-Lm}
  L_{t^m} = \langle x_1,\dots, x_m, y_1, \dots, y_m, z_1, z_2 \, \mid
  \, [x_i,y_j] = \delta_{i,j} z_1 + \delta_{i+1,j} z_2
  \quad \text{and} \\
  [x_i,x_j] = [y_i,y_j] = [x_i,z_1]= [x_i,z_2] = [y_i,z_1] = [y_i,z_2]
  = 0 \text{ for $1 \le i,j \le m$} \rangle,
\end{multline}
a special instance of~\eqref{eq:structure.of.L}.  Furthermore,
$\mathrm{Z}(L) = \Z z_1 + \Z z_2$, and
$[L,L] = \mathrm{Z}(L)$ for $m \ge 2$. 

\begin{lemma} \label{lem:indecomposable} Let $L = L_{t^m}$ with
  $m \ge 2$, and let $\mathcal{F}$ be any field.  Then the
  $\mathcal{F}$-Lie algebra $\mathcal{L} = {_\mathcal{F} L}$ is
  indecomposable.
\end{lemma}

\begin{proof}
  Put
  $\mathcal{Z} = \mathrm{Z}(\mathcal{L}) = \Span_\mathcal{F} \{ z_1,
  z_2 \}$.  A routine check shows that
  \begin{multline} \label{equ:Aut-invariant-W} \mathcal{W} = \big\{ w
    \in \mathcal{L} \mid \dim_\mathcal{F} (\Span_\mathcal{F} \{ [w,v]
    \mid v \in \mathcal{L} \}) \leq 1 \big\} \\ = \big\{ w \in
    \mathcal{L} \mid \Span_\mathcal{F} \{ [w,v] \mid v \in \mathcal{L}
    \} \subseteq \mathcal{F} z_1 \big\} = \Span_\mathcal{F} \{x_m ,
    y_1 \} + \mathcal{Z}
  \end{multline}
  so that $\mathcal{W}$ is a vector subspace and
  $\dim_\mathcal{F}(\mathcal{W}) = 4$.  For a contradiction, suppose
  that $\mathcal{L} = \mathcal{A} \oplus \mathcal{B}$ for non-zero Lie
  ideals $\mathcal{A}, \mathcal{B} \trianglelefteq \mathcal{L}$.
  Since $\mathcal{A}$ is nilpotent, $\mathcal{A}$ has non-zero centre
  $\mathrm{Z}(\mathcal{A}) \ne \{ 0 \}$, and likewise
  $\mathrm{Z}(\mathcal{B}) \ne \{ 0 \}$.  Thus
  $\mathcal{Z} = \mathrm{Z}(\mathcal{A}) \oplus
  \mathrm{Z}(\mathcal{B})$ implies
  $\dim_\mathcal{F}(\mathrm{Z}(\mathcal{A})) =
  \dim_\mathcal{F}(\mathrm{Z}(\mathcal{B})) = 1$.  We deduce that
  $\mathcal{A} \cup \mathcal{B} \subseteq \mathcal{W}$ and hence
  $\mathcal{W} = \mathcal{L}$, in contradiction to
  $\dim_\mathcal{F}(\mathcal{L}) = 2m + 2 > 4$.
\end{proof}

We remark that, in contrast to the situation treated in
Lemma~\ref{lem:indecomposable}, the Lie lattice $L_t$ is already
decomposable over~$\Z$: clearly,
$L_t = (\Z x_1 + \Z y_1 + \Z z_1) \oplus
\Z z_2$ decomposes as a direct sum of two non-zero Lie ideals.

\begin{lemma} \label{lem:Y-generates} Let $L = L_{t^m}$ with
  $m \ge 2$, and let $\mathcal{F}$ be any field.  Then the
  $\mathcal{F}$-Lie algebra $\mathcal{L} = {_\mathcal{F} L}$ is
  generated by
  \[
    \mathcal{Y} = \big\{ w \in \mathcal{L} \smallsetminus
    \mathrm{Z}(\mathcal{L}) \mid \mathrm{C}_\mathcal{L} \big(
    \mathrm{C}_\mathcal{L}(w) \big) = \mathcal{F} w +
    \mathrm{Z}(\mathcal{L}) \big\}
  \]
 if and only if $m \ne 2$.
\end{lemma}

\begin{proof}
  For short we put
  $\mathcal{Z} = \mathrm{Z}(\mathcal{L}) = [\mathcal{L},\mathcal{L}] =
  \Span_\mathcal{F} \{ z_1, z_2 \}$.

  First consider the special case $m=2$.  We claim that $\mathcal{Y}$
  is contained in the proper Lie subalgebra
  $\mathcal{W} = \Span_\mathcal{F} \{ x_2,y_1\} + \mathcal{Z}$; thus
  $\mathcal{Y}$ fails to generate~$\mathcal{L}$.  Indeed, from the
  description~\eqref{equ:Aut-invariant-W} and the definition of
  $\mathcal{Y}$ we see that both $\mathcal{W}$ and $\mathcal{Y}$ are
  $\Aut(\mathcal{L})$-invariant.  Thus it suffices to check that
  $x_1\notin \mathcal{Y}$ and that for every
  $w \in \mathcal{L} \smallsetminus \mathcal{W}$ there exists
  $g \in \Aut(\mathcal{L})$ such that $w g =x_1$.  From
  $\mathrm{C}_\mathcal{L}(x_1) = \Span_\mathcal{F} \{ x_1, x_2 \} +
  \mathcal{Z}$ we deduce that
  $\mathrm{C}_\mathcal{L} ( \mathrm{C}_\mathcal{L}(x_1)) =
  \Span_\mathcal{F} \{ x_1, x_2 \} + \mathcal{Z}$, and this gives
  $x_1\notin \mathcal{Y}$.  Now let
  $w \in \mathcal{L} \smallsetminus \mathcal{W}$.
  Corollary~\ref{cor:corollary-reductive-part} describes the reductive
  part of $\Aut(\mathcal{L})$; compare with~\eqref{equ:form-of-h}.
  From this description we see that there exists
  $g_1 \in \Aut(\mathcal{L})$ such that $w g_1 \in x_1 + \mathcal{W}$.
  Finally, the description of the unipotent radical of
  $\Aut(\mathcal{L})$ in Example~\ref{exa:unipotent-radical-m23} shows
  that there exists $g_2 \in \Aut(\mathcal{L})$ such that
  $w g_1 g_2 = x_1$.

  Now suppose that $m\geq 3$.  We claim that $\mathcal{Y}$ contains
  the generating set
  \[
    x_1,\, x_2,\, \ldots,\, x_{m-2},\, x_m,\quad y_1,\, y_3,\, y_4,\,
    \ldots,\, y_m,\quad \sum\nolimits_{i=1}^m x_i,\,
    \sum\nolimits_{i=1}^m y_i
  \]
  for~$\mathcal{L}$.  Indeed, for $i \in \{1,\ldots, m\}$ it is easily
  checked that
  \begin{align*}
    \mathrm{C}_\mathcal{L}(x_i) %
    & = \Span_\mathcal{F} \{ x_1, \ldots, x_m,
      y_1, \ldots, y_{i-1}, y_{i+2}, \ldots, y_m \} +\mathcal{Z}, \\
    \mathrm{C}_\mathcal{L}(y_i) %
    & = \Span_\mathcal{F} \{ x_1, \ldots, x_{i-2}, x_{i+1}, \ldots, x_m,
      y_1, \ldots, y_m \} +\mathcal{Z}.
  \end{align*}
  For $i \ne m-1$ this implies
  $\mathrm{C}_\mathcal{L} \big( \mathrm{C}_\mathcal{L}(x_i) \big) =
  \mathcal{F} x_i + \mathcal{Z}$, hence $x_i \in \mathcal{Y}$.
  Likewise $y_i \in \mathcal{Y}$ for $i \ne 2$, but it can be seen
  that $x_{m-1}, y_2$ do not belong to~$\mathcal{Y}$.  In order to
  bypass these exceptions, it suffices to show that $\sum_{i=1}^m x_i$
  and $\sum_{i=1}^m y_i$ lie in~$\mathcal{Y}$.  We deduce from
  \[
    \mathrm{C}_\mathcal{L} \left( \sum\nolimits_{i=1}^m x_i \right) =
    \Span_\mathcal{F} \{ x_1,\ldots, x_m, y_2-y_3, y_3-y_4, \ldots,
    y_{m-1}-y_m \}
    + \mathcal{Z}
  \]
 that
  $\mathrm{C}_\mathcal{L} \big( \mathrm{C}_\mathcal{L} ( \sum_{i=1}^m
  x_i ) \big) = \mathcal{F} (\sum_{i=1}^m x_i) + \mathcal{Z}$.  This
  gives $\sum_{i=1}^m x_i \in \mathcal{Y}$ and similarly
  $\sum_{i=1}^m y_i \in \mathcal{Y}$.
\end{proof}

We remark that, for $m=1$, the set
$\mathcal{Y} \subseteq \mathcal{L} = {_\mathcal{F} L_t}$ defined in
Lemma~\ref{lem:Y-generates} coincides with
$\mathcal{L} \smallsetminus \mathrm{Z}(\mathcal{L})$ and thus
generates $\mathcal{L}$ for trivial reasons.

\begin{proposition}
  Let $L = L_{t^m}$ with $m \ge 2$, and let $p$ be a prime.  Then
  \eqref{equ:down-to-earth} holds for every finite extension
  $\mathcal{F}$ of $\Q_p$.
\end{proposition}

\begin{proof}
  For $m > 2$ we can use the criterion established
  in~\cite[Thm.~3.9]{BGS22}: the stronger `rigidity condition'
  \eqref{equ:JRF-equality} follows, for every finite extension
  $\mathcal{F}$ of $\Q_p$, with valuation ring
  ${\scriptstyle \mathcal{O}}$ and valuation ideal $\wp$, from
  Lemmata~\ref{lem:BGS-lemma}, \ref{lem:indecomposable}
  and~\ref{lem:Y-generates}.  For $m=2$ we give a direct proof
  of~\eqref{equ:down-to-earth}, as follows.

  Fix a finite extension $\mathcal{F}$ of $\Q_p$
  of degree $d = [\mathcal{F}:\Q_p]$ and pick a primitive
  element $\alpha$ for the extension
  so that
  \[
  \mathcal{F} = \Q_p(\alpha) = \Q_p \, 1 +
  \Q_p \, \alpha + \ldots + \Q_p \, \alpha^{d-1}. 
  \]
  The $\Q_p$-Lie algebra
  $\widetilde{\mathcal{L}} = {_{\Q_p, \mathcal{F}}
    \mathcal{L}}$ results from the $6$-dimensional $\Q_p$-Lie
  algebra $\mathcal{L} = \Q_p \otimes_\Z L$ with basis
  $x_1,x_2,y_1,y_2,z_1,z_2$, subject to the relations indicated
  in~\eqref{equ:presentation-Lm}, via extension and restriction of
  scalars; we have $\dim_{\Q_p}(\widetilde{\mathcal{L}}) = 6d$
  and $\widetilde{\mathcal{L}}$ admits a $\Q_p$-basis
  consisting of the elementary tensors
  \[
    x_i \alpha^j = \alpha^j \otimes x_i, \quad y_i \alpha^j = \alpha^j
    \otimes y_i, \quad z_i \alpha^j = \alpha^j \otimes z_i, \qquad
    \text{for $i \in \{1,2\}$, $ j \in \{0,\ldots,d-1\}$,}
  \]
  where we write the powers of $\alpha$ on the right so that they are
  visibly separated from scalars coming from~$\Q_p$.  Likewise
  we find it convenient in the calculations below to treat
  $\widetilde{\mathcal{L}}$ formally as a
  $(\Q_p,\mathcal{F})$-bimodule.  We put
  $\widetilde{\mathcal{Z}} = \mathrm{Z}(\widetilde{\mathcal{L}})$ and
  recall that
  \[
    [ \widetilde{\mathcal{L}}, \widetilde{\mathcal{L}} ] =
    \widetilde{\mathcal{Z}} = \Span_\mathcal{F} \{ z_1,
    z_2 \} = \Span_{\Q_p} \big\{ z_i \alpha^j \mid i \in
    \{1,2\}, \; j \in \{0,\ldots,d-1\} \big\}.
  \]
  Furthermore, we observe that with
  \begin{equation} \label{equ:W-tilde-basis}
    \begin{split}
      \widetilde{\mathcal{W}} %
      & = \Span_\mathcal{F} \{ x_2, y_1 \} +
      \mathrm{Z}(\widetilde{\mathcal{L}}) = \big\{ w \in
      \widetilde{\mathcal{L}} \mid \dim_\mathcal{F} [w,
      \widetilde{\mathcal{L}}] \le 1 \big\} \\
      & = \Span_{\Q_p} \{ x_2 \alpha^j \mid 0 \le j < d \}
      \cup \{ y_1 \alpha^j \mid 0 \le j < d \} +
      \mathrm{Z}(\widetilde{\mathcal{L}}) = \big\{ w \in
      \widetilde{\mathcal{L}} \mid \dim_{\Q_p} [w,
      \widetilde{\mathcal{L}}] \le d \big\}
    \end{split}
  \end{equation}
  we obtain a chain of
  $\Aut_{\Q_p}(\widetilde{\mathcal{L}})$-invariant
  $\mathcal{F}$- and hence $\Q_p$-subspaces
  \begin{equation} \label{equ:chain-invertible-subspaces}
    \{ 0 \} \; \subseteq \; \underbrace{[\widetilde{\mathcal{W}},
      \widetilde{\mathcal{L}}]}_{= z_1 \mathcal{F}} \; \subseteq \;
    \widetilde{\mathcal{Z}} \; \subseteq \; \widetilde{\mathcal{W} }
    \; \subseteq \; \widetilde{\mathcal{L}}
  \end{equation}
  with
  $\dim_{\Q_p} [\widetilde{\mathcal{W}},
  \widetilde{\mathcal{L}}] = d$ and
  $\dim_{\Q_p} \widetilde{\mathcal{W}}= 4d$; compare
  with~\eqref{equ:Aut-invariant-W}.

  \smallskip
  
  Now consider an arbitrary automorphism
  $\varphi \in \Aut_{\Q_p}(\widetilde{\mathcal{L}})$.  By
  means of a finite number of basic reductions, we show that $\varphi$
  is contained in the subgroup that appears on the left-hand side
  of~\eqref{equ:down-to-earth}.

  \smallskip

  \noindent \emph{Step $1$.}  By
  Proposition~\ref{pro:short-split-sequ}, the group
  $\Aut_\mathcal{F}(\widetilde{L})$ induces on
  $\widetilde{\mathcal{Z}} = z_2 \mathcal{F} + z_1 \mathcal{F}$ the
  group of all invertible upper triangular matrices; in particular, it
  acts transitively on
  $(z_1 \mathcal{F} \smallsetminus \{0\}) \times ((z_2 \mathcal{F} +
  z_1 \mathcal{F}) \smallsetminus z_1 \mathcal{F})$.  In view of the
  $\varphi$-invariance of $z_1 \mathcal{F}$ and
  $z_2 \mathcal{F} + z_1 \mathcal{F}$
  in~\eqref{equ:chain-invertible-subspaces} we may thus suppose
  without loss of generality that $\varphi$ fixes $z_1$ and $z_2$:
  \[
    z_1 \varphi = z_1 \qquad \text{and} \qquad z_2 \varphi = z_2.
  \]
  
  \smallskip

  \noindent \emph{Step $2$.}  Next we focus on
  $[\widetilde{W},\widetilde{L}] = z_1\mathcal{F} =
  \Span_{\Q_p} \{ z_1 \alpha^j \mid 0 \le j < d \}$, with the
  aim to reduce to the situation where $\varphi$ induces the identity
  on this subspace.  In view of~\eqref{equ:chain-invertible-subspaces}
  we may write
  \[
    (z_1 \alpha^j) \varphi = z_1 \lambda_j \quad \text{for suitable
      $\lambda_j \in \mathcal{F}$,} \qquad \text{for $0 \le j \le d$.}
  \]
  Due to the reduction in Step~1 we have $\lambda_0 = 1$, and
  $\lambda_d$ is actually determined by
  $\lambda_0, \ldots, \lambda_{d-1}$, because $\alpha^d$ can be
  expressed as a $\Q_p$-linear combination of
  $\alpha^0, \ldots, \alpha^{d-1}$; in~\eqref{equ:why-d-mattes} below
  it becomes clear why our analysis includes $\lambda_d$.
  Furthermore, for $0 \le j < d$, the images of $x_1 \alpha^j$ and
  $y_1 \alpha^j \in \widetilde{\mathcal{W}}$ under $\varphi$ can be
  written, modulo $\widetilde{\mathcal{Z}}$, as $\mathcal{F}$-linear
  combinations
  \[
    (x_1 \alpha^j) \varphi \equiv_{\widetilde{\mathcal{Z}}} x_1 a_j +
    y_2 b_j + x_2 a_j' + y_1 b_j' \qquad \text{and} \qquad (y_1
    \alpha^j) \varphi \equiv_{\widetilde{\mathcal{Z}}} x_2 c_j  + y_1 d_j.
  \]
  For $0 \le j <d$ we deduce that
  \[
    0 = [x_1, x_1 \alpha^j] \varphi = [x_1 a_0 + y_2 b_0 + \ldots, x_1
    a_j + y_2 b_j + \ldots] \equiv_{z_1 \mathcal{F}} z_2 (a_0 b_j -
    b_0 a_j)
  \]
  so that $a_0 b_j = a_j b_0$.  In a similar way, for $0 \le j \le d$
  and $0 \le i \le \min\{1,j\}$ we see that
  \[
    z_1 \lambda_j = (z_1 \alpha^j) \varphi = [x_1 \alpha^{j-i}, y_1
    \alpha^i] \varphi = [ x_1 a_{j-i} + y_2 b_{j-i} + \ldots, x_2 c_i
    + y_1 d_i] = z_1 (a_{j-i} d_i - b_{j-i} c_i)
  \]
  so that $\lambda_j = a_{j-i} d_i - b_{j-i} c_i$.  Using
  $b_0 a_{j-1} = a_0 b_{j-1}$ to modify the underlined terms and
  $\lambda_0 =1$ for the final simplification, we deduce that for
  $1 \le j \le d$,
  \begin{multline*}
    \lambda_1 \lambda_{j-1} %
    = (a_0 d_1 - b_0 c_1) (a_{j-1} d_0 - b_{j-1} c_0) = a_0 d_1
    a_{j-1} d_0 - \uwave{b_0 c_1 a_{j-1} d_0} - \uwave{a_0
      d_1 b_{j-1} c_0} + b_0 c_1 b_{j-1} c_0 \\
    = a_0 d_0 (a_{j-1} d_1 - b_{j-1} c_1 ) - b_0 c_0 (a_{j-1} d_1 -
    b_{j-1} c_1 ) = (a_0 d_0 - b_0 c_0) \lambda_j = \lambda_0
    \lambda_j = \lambda_j.
  \end{multline*}
  By induction, we obtain
  $\lambda_j = \lambda_1^{\, j}$ for $0 \le j \le d$.  Let
  $f = \sum_{j=0}^d f_j t^j \in \Q_p[t]$ denote the minimal
  polynomial of $\alpha$ over $\Q_p$.  Then
  \begin{equation} \label{equ:why-d-mattes}
    0 = \big( z_1 f(\alpha) \big) \varphi = \sum\nolimits_{j=0}^d f_j
    \, (z_1 \alpha^j) \varphi = \sum\nolimits_{j=0}^d f_j \, (z_1
    \lambda_j) = z_1 \big( \sum\nolimits_{j=0}^d f_j \lambda_j \big) =
    z_1 f(\lambda_1)
  \end{equation}
  implies $f(\lambda_1) = 0$.  Hence $\alpha$ and $\lambda_1$ are
  Galois conjugates in
  $\mathcal{F} = \Q_p(\alpha) = \Q_p(\lambda_1)$.
  Modifying $\varphi$ by a field automorphism, i.e.\ an element of
  $\Aut(\mathcal{F} \,|\, \Q_p)$, we may suppose without loss of
  generality that
  \[
    \big( z_1 \alpha^j \big) \varphi = z_1 \alpha^j \qquad \text{for
      $0 \le j < d$.}
  \]

  \smallskip

  \noindent \emph{Step $3$.}  Next we focus on the action of $\varphi$
  on $\widetilde{\mathcal{Z}} = z_1 \mathcal{F} + z_2 \mathcal{F}$
  modulo $[\widetilde{W},\widetilde{L}] = z_1 \mathcal{F}$; this
  factor space admits $z_2 \alpha^j$, $0 \le j <d$, as a $\Q_p$-basis.
  In view of~\eqref{equ:chain-invertible-subspaces} we may write
  \[
    (z_2 \alpha^j) \varphi \equiv_{z_1 \mathcal{F}} z_2 \mu_j \quad
    \text{for suitable $\mu_j \in \mathcal{F}$,} \qquad \text{for
      $0 \le j \le d$;}
  \]
  our aim is to show that $\mu_j = \beta^j$, with $\beta = \mu_1$
  Galois conjugate to $\alpha$.

  Due to the reduction in Step~1 we have $\mu_0 = 1$, and $\mu_d$ is
  actually determined by $\mu_0, \ldots, \mu_{d-1}$; compare with
  Step~2.  For $0 \le j < d$, the images of $x_1 \alpha^j$ and
  $y_2 \alpha^j$ under $\varphi$ can be written,
  modulo~$\widetilde{\mathcal{Z}}$, as $\mathcal{F}$-linear
  combinations
  \[
    (x_1 \alpha^j) \varphi \equiv_{\widetilde{\mathcal{Z}}} x_1 a_j +
    y_2 b_j + x_2 a_j' + y_1 b_j' \qquad \text{and} \qquad (y_2
    \alpha^j) \varphi \equiv_{\widetilde{\mathcal{Z}}} x_1 c_j + y_2
    d_j + x_2 c_j' + y_1 d_j'.
  \]
  In Step~2 we saw that $a_0 b_j = a_j b_0$ for $0 \le j <d$.
  Furthermore, for $0 \le j \le d$ and $0 \le i \le \min\{1,j\}$ we
  get, modulo $z_1 \mathcal{F}$,
  \[
    z_2 \mu_j \equiv_{z_1 \mathcal{F}} (z_2 \alpha^j) \varphi = [x_1
    \alpha^{j-i}, y_2 \alpha^i] \varphi = [x_1 a_{j-i} + y_2 b_{j-i} +
    \ldots, x_1 c_i + y_2 d_i + \ldots] \equiv_{z_1 \mathcal{F}} z_2
    (a_{j-i} d_i - b_{j-i} c_i)
  \]
  so that $\mu_j = a_{j-i} d_i - b_{j-i} c_i$.  A similar argument as
  in Step~2 shows that $\mu_j = \mu_1^{\, j}$ for $0 \le j \le d$ and
  that $\alpha$ and $\beta = \mu_1$ are Galois conjugates in
  $\mathcal{F} = \Q_p(\alpha) = \Q_p(\beta)$.

  \smallskip

  \noindent \emph{Step $4$.} We analyse further the action of
  $\varphi$ on
  $\widetilde{\mathcal{Z}} = z_1 \mathcal{F} + z_2 \mathcal{F}$.  So
  far we have reduced to the situation in which $\varphi$ acts as the
  identity on $z_1 \mathcal{F}$ and
  \[
    (z_2 \alpha^j) \varphi = z_2 \beta^j + z_1 \nu_j \quad \text{for
      $1 \le j \le d$,}
  \]
  where $\beta$ denotes a Galois conjugate of $\alpha$ and
  $0 = \nu_0, \nu_1, \ldots, \nu_d \in \mathcal{F}$ are suitable
  coefficients.  As before, $\nu_d$ is actually determined by the
  previous parameters.  Proposition~\ref{prop:Grp-new-coordinates}
  describes the pointwise stabiliser of
  $\mathrm{Z}(\widetilde{\mathcal{L}})$ inside
  $\Aut_\mathcal{F}(\widetilde{\mathcal{L}})$; a short reflection
  reveals that this stabiliser acts transitively on
  $\widetilde{\mathcal{L}} \smallsetminus \widetilde{\mathcal{W}}$ and
  consequently we may suppose without loss of generality that
  \[
    x_1 \varphi = x_1;
  \]
  in particular, the abelian Lie subalgebra
  \[
    \widetilde{\mathcal{X}} =
    \mathrm{C}_{\widetilde{\mathcal{L}}}(x_1) =
    \mathrm{C}_{\widetilde{\mathcal{L}}}(x_1 \mathcal{F} + x_2
    \mathcal{F}) = \Span_\mathcal{F} \{ x_1, x_2, z_1, z_2\}
  \]
  is $\varphi$-invariant.  For $0 \le j \le d$ we deduce from
  \[
    [x_1, (y_2 \alpha^j)\varphi] = [x_1, y_2 \alpha^j] \varphi = (z_2
    \alpha^j) \varphi = z_2 \beta^j + z_1 \nu_j
  \]
  that, modulo $\widetilde{\mathcal{X}}$,
  \[
    (y_2 \alpha^j) \varphi \equiv_{\widetilde{\mathcal{X}}} y_1 \nu_j
    + y_2 \beta^j;
  \]
  in particular, $y_2 \varphi \equiv y_2$ modulo
  $\widetilde{\mathcal{X}}$.  Furthermore,
  $(x_1 \alpha^j)\varphi \in \widetilde{\mathcal{X}}$ and
  \[
    [(x_1 \alpha^j)\varphi, y_2 ] = [(x_1 \alpha^j)\varphi, y_2\varphi
    ] = [x_1 \alpha^j, y_2] \varphi = (z_2 \alpha^j) \varphi = z_2
    \beta^j + z_1 \nu_j
  \]
  yield, modulo $\widetilde{\mathcal{Z}}$,
  \[
    (x_1 \alpha^j) \varphi \equiv_{\widetilde{\mathcal{Z}}} x_1
    \beta^j + x_2 \nu_j.
  \]
  For $0 \le j < d$ we deduce from
  \[
    z_2 \beta^{j+1} + z_1 (\beta^j \nu_1 + \beta \nu_j) = [ x_1
    \beta^j + x_2 \nu_j, y_1 \nu_1 + y_2 \beta ] = [x_1 \alpha^j, y_2
    \alpha] \varphi = (z_2 \alpha^{j+1} ) \varphi = z_2 \beta^{j+1} +
    z_1 \nu_{j+1}
  \]
  that $\nu_{j+1} = \beta \nu_j + \beta^j \nu_1$.  By recursion, this
  gives
  \[
    \nu_j = j \beta^{j-1} \nu_1 \quad \text{for $0 \le j \le d$.}
  \]
  Let $f = \sum_{j=0}^d f_j t^j \in \Q_p[t]$ denote the
  minimal polynomial of $\alpha$ and of its conjugate $\beta$ over
  $\Q_p$.  Then
  \begin{multline*}
    0 = \big( z_2 \underbrace{f(\alpha)}_{=0} \big) \varphi = \big(
    \sum\nolimits_{j=0}^d f_j (z_2 \alpha^j) \big) \varphi =
    \sum\nolimits_{j=0}^d f_j \big( (z_2 \alpha^j) \varphi \big) =
    \sum\nolimits_{j=0}^d f_j \, (z_2 \beta^j + z_1 \nu_j) \\
    = \sum\nolimits_{j=0}^d f_j \, \big( z_2 \beta^j + z_1 (j
    \beta^{j-1} \nu_1) \big) = z_2 \underbrace{f(\beta)}_{=0} + z_1
    (\nu_1 f'(\beta)) = z_1 (\nu_1 \underbrace{f'(\beta)}_{\ne 0})
  \end{multline*}
  implies $\nu_1 = 0$, hence $\nu_j = 0$ for $0 \le j \le d$ and
  \[
    (z_2 \alpha^j) \varphi = z_2 \beta^j \quad \text{for $0 \le j <
      d$.}
  \]

  \smallskip

  \noindent \emph{Step $5$.}  Finally, let us see how $\varphi$ acts
  modulo the centre.  In Step~4 we saw that $y_2 \varphi \equiv y_2$
  modulo~$\widetilde{\mathcal{X}}$.
  Proposition~\ref{prop:Grp-new-coordinates} describes the pointwise
  stabiliser of $\widetilde{\mathcal{Z}}$ inside
  $\Aut_\mathcal{F}(\widetilde{\mathcal{L}})$; in particular, this
  stabiliser acts transitively on $y_2 + \widetilde{\mathcal{X}}$,
  even if we add the condition that $x_1$ is to remain fixed: in the
  notation of the proposition, we can take
  \[
    X_1 =
    \begin{pmatrix}
      1 & 0 \\ c_1 & 1
    \end{pmatrix}
    \quad \text{and} \quad X_2 =
    \begin{pmatrix}
      0 & 0 \\ c_2 & 0
    \end{pmatrix},
    \qquad \text{where $c_1, c_2 \in \mathcal{F}$ are
      free parameters.}
  \]
  Thus we may suppose, without interfering with the previous
  reductions, that $\varphi$ fixes $y_2$, i.e.\
  \[
    y_2 \varphi = y_2.
  \]
  From $x_2 \varphi \in \widetilde{\mathcal{X}}$ and
  $[x_2 \varphi, y_2 ] = [x_2,y_2] \varphi = z_1 \varphi = z_1$ we
  deduce that $x_2 \varphi \equiv x_2$ modulo
  $\widetilde{\mathcal{Z}}$.  Recall that
  $y_1 \in \widetilde{\mathcal{W}}$ implies
  $y_1 \varphi \in \widetilde{\mathcal{W}}$; moreover, $\varphi$ fixes
  $x_1$ and~$z_1$.  Hence
  $[x_1, y_1 \varphi] = [x_1, y_1] \varphi = z_1 \varphi = z_1$ gives
  $y_1 \varphi \equiv y_1$ modulo
  $x_2 \mathcal{F} + \widetilde{\mathcal{Z}}$.  From
  $[y_1 \varphi, y_2] = [y_1,y_2]\varphi = 0$ we conclude that
  $y_1 \varphi \equiv y_1$ modulo~$\widetilde{\mathcal{Z}}$.  We have
  gained
  \[
    x_2 \varphi \equiv_{\widetilde{\mathcal{Z}}} x_2 \qquad \text{and}
    \qquad y_1 \varphi \equiv_{\widetilde{\mathcal{Z}}} y_1
  \]
  Now let
  $0 \le j < d$.  From
  \[
    [x_1, (y_1 \alpha^j)\varphi] = [x_1, y_1 \alpha^j]\varphi = (z_1
    \alpha^j) \varphi = z_1 \alpha^j \quad \text{and}
    \quad [y_1,(y_1 \alpha^j)\varphi] = [y_2, (y_1 \alpha^j) \varphi] = 0
  \]
  we see that
    $(y_1 \alpha^j)\varphi \equiv_{\widetilde{\mathcal{Z}}}
    y_1 \alpha^j$.
  Similarly,
  $[x_1, (y_2 \alpha^j)\varphi] = (z_2 \alpha^j)\varphi = z_2 \beta^j$
  and $[y_1,(y_2 \alpha^j)\varphi] = [y_2,(y_2 \alpha^j)\varphi] = 0$ imply
  $(y_2 \alpha^j)\varphi \equiv_{\widetilde{\mathcal{Z}}}
  y_2 \beta^j$.  Moreover
  \[
    z_1 \alpha = (z_1 \alpha)\varphi = [x_2, y_2 \alpha]\varphi
    = [x_2,(y_2 \alpha)\varphi] = [x_2, y_2 \beta ] = z_1 \beta
  \]
  implies $\alpha = \beta$.
   
  In summary, this shows that $\varphi$ fixes pointwise the centre
  $\widetilde{\mathcal{Z}}$, and that, modulo
  $\widetilde{\mathcal{Z}}$,
  \[
    (y_1 \alpha^j)\varphi \equiv_{\widetilde{\mathcal{Z}}} y_1
    \alpha^j \quad \text{and} \quad (y_2 \alpha^j)\varphi
    \equiv_{\widetilde{\mathcal{Z}}} y_2 \alpha^j, \qquad \text{for
      $0 \le j < d$.}
  \]
  Finally, we observe that
  $(x_1 \alpha^j) \varphi, (x_2 \alpha^j) \varphi \in
  \widetilde{\mathcal{X}}$ satisfy
  $[(x_1 \alpha^j)\varphi, y_2] = [ x_1 \alpha^j, y_2] \varphi = (z_2
  \alpha^j)\varphi = z_2 \alpha^j$ and, by similar considerations,
  $[(x_2 \alpha^j)\varphi, y_1] = 0$ and
  $[ (x_2 \alpha^j) \varphi, y_2] \varphi = z_1 \alpha^j$.  From this
  we conclude that, modulo $\widetilde{\mathcal{Z}}$,
  \[
    (x_1 \alpha^j)\varphi \equiv_{\widetilde{\mathcal{Z}}} x_1
    \alpha^j \quad \text{and} \quad (x_2 \alpha^j)\varphi
    \equiv_{\widetilde{\mathcal{Z}}} x_2 \alpha^j, \qquad \text{for
      $0 \le j < d$.}
  \]
  As $\widetilde{\mathcal{Z}} = \mathrm{Z}(\widetilde{\mathcal{L}})$
  it follows that
  $\varphi \in
  \mathrm{C}_{\Aut_{\Q_p}(\widetilde{\mathcal{L}})} \big(
  \widetilde{\mathcal{L}} / \mathrm{Z}(\widetilde{\mathcal{L}}) \big)$
  is contained in the subgroup on the left-hand side
  of~\eqref{equ:down-to-earth}.
\end{proof}


\subsection{The local pro-isomorphic zeta functions of groups
  $\widetilde{\Gamma}$ associated to $L_{t^2}$ and
  $L_{t^3}$} \label{sec:local-Lt2-Lt3} We return to the setting
described at the beginning of the section.  Let $k$ be a number field
of absolute degree $d = [k : \Q]$, with ring of
integers~$\mathfrak{o}$.  Let $\widetilde{L} = {_{\Z,\mathfrak{o}} L}$
be the nilpotent $\Z$-Lie lattice associated to $L_{t^m}$ for
$m \in \{2,3\}$ via `base extension', with algebraic automorphism
group~$\widetilde{\G} = \mathbf{Aut}(\widetilde{L}) \le \GL_{nd}$,
where $n = \dim_\Z L = 2m+2$, and let $p$ be a prime.  The basic
ingredients for the `fine' Euler decomposition established in
\cite[Prop.~3.14]{BGS22} are the natural isomorphisms
${_{\Z_p} \mathfrak{o}} = \Z_p \otimes_\Z \mathfrak{o} \cong
\prod_{\mathfrak{p} \mid p} \mathfrak{o}_\mathfrak{p}$ and
${_{\Q_p} k} = \Q_p \otimes_\Q k \cong \prod_{\mathfrak{p} \mid p}
k_\mathfrak{p}$; we summarise the technical steps and implications in
our setting.  We write $\widetilde{\mathbf{H}}$ for the reductive part
of the $1$-component~$\widetilde{\G}^\circ$.  As described in
Section~\ref{section:reduction-integral}, the local zeta function
associated to the $\Z_p$-Lie lattice
$\widetilde{L}_p = {_{\Z_p} \widetilde{L}}$ can be expressed as a
$p$-adic integral
\begin{equation} \label{equ:integral-Lp-tilde}
  \zeta_{\widetilde{L}_p}^\mathrm{iso}(s) = \int_{\widetilde{H}_p^+}
  \lvert \det h \rvert_p^{\, s} \, \theta_0(h) \, \theta_1(h) \,
  \mathrm{d}\mu_{\widetilde{H}_p}(h),
\end{equation}
where $\widetilde{H}_p = \widetilde{\mathbf{H}}(\Q_p)$,
$\widetilde{H}_p^+ = \widetilde{H}_p \cap \Mat_{nd}(\Z_p)$ and
$\theta_0, \theta_1 \colon \widetilde{H}_p \to \R_{\ge 0}$ are
suitable volume functions, modulo a small technical issue to be taken
care of: while $\G = \mathbf{Aut}(L)$ is connected (as we proved), the
group $\widetilde{\G}$ is typically not connected.  But in the
presence of~\eqref{equ:JRF-equality} or the somewhat weaker
condition~\eqref{equ:down-to-earth}, which we established in
Section~\ref{sec:local-rgidity}, the finite group scheme
$\mathbf{F} \cong \mathbf{Aut}(\mathfrak{o} \,|\, \Z) \cong
\mathbf{Aut}(k \,|\, \Q)$, which potentially renders the group
$\widetilde{\G}$ non-connected, has the feature that
$\mathbf{F}(\Z_p) = \mathbf{F}(\Q_p)$ and can thus be safely ignored,
by using the same argument as in the proof
of~\cite[Prop.~2.1]{dSLu96}.  Moreover, the group
$\widetilde{G}_p = \widetilde{\G}(\Q_p)$ almost, but not quite
decomposes as a direct product indexed by the primes
$\mathfrak{p} \mid p$.  In the reductive part $\widetilde{H}_p$ the
troublesome central automorphisms disappear and we have
\[
  \widetilde{H}_p \cong \prod_{\mathfrak{p} \mid p} H_\mathfrak{p}
  \qquad \text{with
    $H_\mathfrak{p} = \mathbf{H}(k_\mathfrak{p}) \cong
    \widetilde{\mathbf{H}}_\mathfrak{p}(\Q_p)$ for each
    $\mathfrak{p} \mid p$,}
\]
where $\mathbf{H}$ is the reductive part of of the original group $\G$
and, setting
$d_\mathfrak{p} = \dim_{\Z_p}(\mathfrak{o}_\mathfrak{p}) =
[k_\mathfrak{p} : \Q_p]$, we denote by
$\widetilde{\mathbf{H}}_\mathfrak{p}$ the reductive part of the
algebraic automorphism group
$\widetilde{\G}_\mathfrak{p} =
\mathbf{Aut}(\widetilde{L}_\mathfrak{p}) \le \GL_{n d_\mathfrak{p}}$
of the $\Z_p$-Lie lattice
$\widetilde{L}_\mathfrak{p} = {_{\Z_p,\mathfrak{o}_\mathfrak{p}} L}$,
which we analysed in Section~\ref{sec:local-rgidity}.  The next step
is to transform the integral in~\eqref{equ:integral-Lp-tilde} over
$\widetilde{H}_p^+$ into a product of integrals over
$H_\mathfrak{p}^+ = H_\mathfrak{p} \cap
\Mat_n(\mathfrak{o}_\mathfrak{p})$ for $\mathfrak{p} \mid p$; this
essentially uses the natural isomorphism between locally compact
groups
$\big( \mathsf{Res}_{\mathfrak{o}_\mathfrak{p} \,|\, \Z_p}(\mathbf{H})
\big) (\Q_p) \cong \mathbf{H}(k_\mathfrak{p})$, but one also needs to
pay attention to the accommodation of central automorphisms.

Modulo this small wrinkle, it is not difficult to carry out the
analysis in Section~\ref{section:reduction-integral} for the local
field $k_\mathfrak{p}$ in place of $\Q_p$ to obtain
\begin{equation} \label{equ:integral-p-pfrak}
  \zeta_{\widetilde{L}_p}^\mathrm{iso}(s) = \prod_{\mathfrak{p} \mid
    p} \int_{H_\mathfrak{p}^+} \lvert \det h \rvert_\mathfrak{p}^{\,
    s} \, \theta_0(h) \, \theta_1(h)^d \,
  \mathrm{d}\mu_{H_\mathfrak{p}}(h),
\end{equation}
where, for each $\mathfrak{p} \mid p$, the volume functions
$\theta_0, \theta_1$ are defined in analogy to
Section~\ref{section:reduction-integral} (we refrain from adding the
decoration `$\mathfrak{p}$') and $\mu_{H_\mathfrak{p}}$ denotes the
right Haar measure on $H_\mathfrak{p}$ with the normalisation
$\mu_{H_\mathfrak{p}}(H_\mathfrak{p}(\mathfrak{o}_\mathfrak{p})) = 1$;
compare with the discussion in~\cite[\S3]{BGS22}, in particular
with~\cite[Prop.~3.14]{BGS22}.  It is worth pointing out that on the
right-hand side of~\eqref{equ:integral-p-pfrak} the exponent of
$\theta_1(h)$ is $d = [k:\Q]$ and not the corresponding local
parameter~$d_\mathfrak{p}$; this feature results from the treatment of
central automorphisms and justifies that we consider the finite
product of integrals as one `package'.

\smallskip

It remains to carry out the explicit calculation of the integrals
in~\eqref{equ:integral-p-pfrak} arising from the concrete cases
$L = L_{t^2}$ and $L = L_{t^3}$.  Consider first $L = L_{t^2}$.  The
calculation of the integral in Section~\ref{section:x^2} carries over
with little change.  The only material difference is that
$\theta_1(h)$ in the integrand is replaced by $\theta_1(h)^d$.  The
intermediate integral~\eqref{equ:intermediate-integral-t2} now takes
the form
\begin{align*}
  \int_{\substack{(A,\nu)\in \dot{H}_\mathfrak{p} \text{ with}\\
  v_\mathfrak{p}(A) \ge 0 \text{ and} \\
  v_\mathfrak{p}(A) + v_\mathfrak{p}(\nu) \geq 0}} \,
  \lvert \det A \rvert_\mathfrak{p}^{\, 4s-8 d -2} \, \lvert \nu
  \rvert_\mathfrak{p}^{\, 5s-12 d} 
  \, \mathrm{d}\mu_\mathfrak{p}(A,\nu) 
\end{align*}
where
$\dot{H}_\mathfrak{p} = \dot{\mathbf{H}}(k_\mathfrak{p}) =
\GL_2(k_\mathfrak{p}) \times \GL_1(k_\mathfrak{p})$ and the valuation
map $v_\mathfrak{p}$ on $k_\mathfrak{p}$ and on
$\Mat_2(k_\mathfrak{p})$ replaces the $p$-adic valuation $v_p$ used
previously.  Due to the dependence on~$d$, we then obtain
\[
  \theta_0 \big( \xi_\mathbf{e}(\pi)^\rho \big) \, \theta_1 \big(
  \xi_\mathbf{e}(\pi)^\rho \big)^d = q_\mathfrak{p}^{\,
    (8d+2)e_1+12de_2+(4d+4)e_3},
\]
where $\pi$ now denotes a uniformising element for $k_\mathfrak{p}$,
that is $v_\mathfrak{p}(\pi) = 1$, and where $q_\mathfrak{p}$ denotes
the residue field size of~$k_\mathfrak{p}$, that is
$q_\mathfrak{p} = \lvert \mathfrak{o} / \mathfrak{p} \rvert = \lvert
\pi \rvert_\mathfrak{p}^{\, -1}$.  We obtain the formula
\[
  \zeta_{\widetilde{L}_p}^\mathrm{iso}(s) = \prod_{\mathfrak{p} \mid
    p} \frac{ 1 + q_\mathfrak{p}^{\, 8d+2-4s}}{ (1 -
    q_\mathfrak{p}^{\, 4 d +4-3s}) (1 -  q_\mathfrak{p}^{\,
      8d +3-4s} ) (1 - q_\mathfrak{p}^{\, 12d-5s})}.
\]
This is the local pro-isomorphic zeta function, at the prime~$p$, of
the class-two nilpotent group
$\widetilde{\Gamma} = \widetilde{\Gamma}_{t^2}$ associated
to~$\widetilde{L} = \widetilde{L}_{t^2}$, as described in
Section~\ref{section:class-2-correspondence}.  It is straightforward
to deduce Theorem~\ref{thm:t2-base-change} and the assertions in
Remark~\ref{rem:t2-base-change}.

Finally we consider $L=L_{t^3}$. The calculation of the integral in
Section~\ref{section:x^3} carries over with little change.  Indeed,
the treatment there was already performed so that it applies equally
well to the more general situation.  Again, the only material
difference is that $\theta_1(h)$ in the integrand is replaced by
$\theta_1(h)^d$.  The intermediate
integral~\eqref{equ:intermediate-integral-t3} now takes the form
\begin{align*}
  \int_{\substack{(A,\nu)\in \dot{H} \text{ with}\\
  v_\mathfrak{p}(A) \ge 0 \text{ and} \\
  v_\mathfrak{p}(A) + v_\mathfrak{p}(\nu) \geq 0}}
  \lvert \det A \rvert_\mathfrak{p}^{\, 5s-12d} \, \lvert \nu
  \rvert_\mathfrak{p}^{\, -s+6d}\theta_0 \big( (A,\nu)^\rho \big)
  \, \mathrm{d}\mu_\mathfrak{p}(A,\nu) 
\end{align*}
where
$\dot{H}_\mathfrak{p} = \dot{\mathbf{H}}(k_\mathfrak{p}) =
\GL_2(k_\mathfrak{p}) \times \GL_1(k_\mathfrak{p})$ and the valuation
map $v_\mathfrak{p}$ on $k_\mathfrak{p}$ and on
$\Mat_2(k_\mathfrak{p})$ replaces the $p$-adic valuation $v_p$ used
previously.  Due to the dependence on~$d$, we have
$\theta_1 \big( \xi_\mathbf{e}(\pi)^\rho \big)^d = q_\mathfrak{p}^{\,
  12 d e_1+ 5d e_2+ 18d e_3}$, where $\pi$ is a uniformising element
for $k_\mathfrak{p}$ and where $q_\mathfrak{p}$ denotes the residue
field size of~$k_\mathfrak{p}$.
Equation~\ref{eqn-t3-second-interm-expr} becomes
\[
  \mathcal{Z}_{\mathfrak{p}}(s) = \mathcal{Z}_{\dot{\mathbf{H}},\rho,
    \theta, \mathfrak{p}}(s) = \sum_{w\in W} q_\mathfrak{p}^{\,
    -\lf(w)} \sum_{\substack{\mathbf{e} \, \in\, \cone \text{ with}\\
      e_1>0\ \text{if}\ w\neq 1}} X_1^{e_1}X_2^{e_2}X_3^{e_3} \,
  \tilde\theta(\xi_\mathbf{e}(\pi)),
\]
where the numerical data is now
\[
  X_1 = q_\mathfrak{p}^{\, 12d + 2 - 5s}, \quad X_2 =
  q_\mathfrak{p}^{\, 6d-4-s}, \quad X_3 = q_\mathfrak{p}^{\,
    18d+8-9s}.
\]
The remaining calculations go through unchanged: the analogue of
equation~\eqref{equ:X1X2X3-expression} yields the formula
\begin{equation} \label{equ:L-iso-widetilde-t3}
  \zeta_{\widetilde{L}_p}^\mathrm{iso}(s) = \prod_{\mathfrak{p} \mid
    p} \frac{(1 - q_\mathfrak{p}^{\, 24d + 5 -10s}) \,
    V_{\mathfrak{p}}(s)}{ (1 - q_\mathfrak{p}^{\, 12d+3-5s})^2 \, (1 -
    q_\mathfrak{p}^{\, 18d+11-9s}) \, (1-q_\mathfrak{p}^{\, 30d-11s})
    \, (1-q_\mathfrak{p}^{54d+7-21s})},
\end{equation}
where
\begin{equation} \label{equ:V-p}
  \begin{split}
    V_{\mathfrak{p}}(s) & = \frac{-q_\mathfrak{p}^{\, 90d+14-36s} -
      q_\mathfrak{p}^{\, 78d +12-31s} - q_\mathfrak{p}^{\, 66d+9-26s}
      - q_\mathfrak{p}^{54d +5-21s} + q_\mathfrak{p}^{36d+9-15s} +
      q_\mathfrak{p}^{24d+5-10s} + q_\mathfrak{p}^{12d+2-5s}+1}{1 +
      q_\mathfrak{p}^{\, 12d+3-5s}}
    \\
    & = 1 + q_\mathfrak{p}^{\, 12d + 2 -5s} - q_\mathfrak{p}^{\, 12d +
      3 -5s} + q_\mathfrak{p}^{\, 24d + 6 - 10s} - q_\mathfrak{p}^{\,
      54d + 5 - 21s} + q_\mathfrak{p}^{\, 66d + 8 - 26s} -
    q_\mathfrak{p}^{\, 66d + 9 - 26s} - q_\mathfrak{p}^{\, 78d + 11 -
      31s}.
  \end{split}
\end{equation}
This is the local pro-isomorphic zeta function, at the prime~$p$, of
the class-two nilpotent group
$\widetilde{\Gamma} = \widetilde{\Gamma}_{t^2}$ associated
to~$\widetilde{L} = \widetilde{L}_{t^3}$, as described in
Section~\ref{section:class-2-correspondence}.  We formulate, in
analogy to Theorem~\ref{thm:t2-base-change}, a partial generalisation of
Corollary~\ref{cor:t3}.

\begin{theorem} \label{thm:t3-base-change} Let $k$ be a number field
  of absolute degree $d = [k:\Q]$, with ring of
  integers~$\mathfrak{o}$.  Let
  $\widetilde{\Gamma} = \widetilde{\Gamma}_{t^3,k}$ be the class-two
  nilpotent group of Hirsch length $8d$ and with rank-$2d$ centre,
  corresponding to the class-two nilpotent $\Z$-Lie lattice
  $\widetilde{L} = \widetilde{L}_{t^3,k}$ which results from the Lie
  lattice $L = L_{t^3}$ via `base extension' as defined above.

  Then the pro-isomorphic zeta function of the group
  $\widetilde{\Gamma}$ is
  \[
    \zeta^{\wedge}_{\widetilde{\Gamma}}(s) =
    \frac{\zeta_k(5s-(12d+3))^2 \, \zeta_k(9s-(18d+11)) \,
      \zeta_k(11s-30d) \, \zeta_k(21s-(54d+7))}{\zeta_k(10s-(24d+5))}
    \, \omega(s),
  \]
  where $\zeta_k(s)$ denotes the Dedekind zeta function of~$k$ and
  \begin{equation} \label{equ:omega-prod} \omega(s) =
    \prod_\mathfrak{p} V_\mathfrak{p}(s)
  \end{equation}
  with the product running over all non-archimedean primes
  $\mathfrak{p}$ of $k$ and $V_\mathfrak{p}(s)$ defined as in
  \eqref{equ:V-p}.
\end{theorem}

\begin{remark} \label{rem:t3-base-change} For $k = \Q$, i.e.\ $d=1$,
  the description is in agreement with Corollary~\ref{cor:t3};
  compare~\eqref{equ:psi-appears}.  Similar to the special situation
  covered in Section~\ref{sec:mero-cont}, it is routine to check that
  the infinite product in~\eqref{equ:omega-prod} converges absolutely
  and yields a holomorphic function on the half-plane consisting of
  all $s \in \C$ with
  \[
    \operatorname{Re}(s) > \max \big\{ \tfrac{12d+4}{5},
    \tfrac{24d+7}{10}, \tfrac{54d+6}{21}, \tfrac{66d+10}{26},
    \tfrac{78d+12}{31} \big\} =
    \begin{cases}
      \frac{12d+4}{5} & \text{if $d \in \{1,2\}$,} \\
      \frac{18d+2}{7} & \text{if $d \ge 3$.}
    \end{cases}
  \]
  Consequently, for number fields $k$ of absolute degree $d \ge 3$,
  the pro-isomorphic zeta function of
  $\widetilde{\Gamma} = \widetilde{\Gamma}_{t^3,k}$ has abscissa of
  convergence~$(30d+1)/11$ and can be meromorphically continued at
  least to $\{ s \in \C \mid \operatorname{Re}(s) > (18d+2)/7 \}$ with
  a simple pole at $s = (30d+1)/11$.  For quadratic fields $k$, i.e.\
  $d=2$, there is an extra twist, but a routine analysis shows that
  the pro-isomorphic zeta function has abscissa of convergence~$28/5$
  and can be meromorphically continued at least to
  $\{ s \in \C \mid \operatorname{Re}(s) > 11/2 \}$ with a simple pole
  at $s=28/5$.  Similar to Remark~\ref{rem:asymp-growth-t3}, the
  asymptotic growth of pro-isomorphic subgroups in
  $\widetilde{\Gamma}$ can be described by means of a suitable
  Tauberian theorem.  Via the Euler product, the formula for
  $\zeta^{\wedge}_{\widetilde{\Gamma}}(s)$ incorporates the
  description \eqref{equ:L-iso-widetilde-t3} of the local
  pro-isomorphic zeta functions
  $\zeta^{\wedge}_{\widetilde{\Gamma},p}(s) =
  \zeta_{\widetilde{L}_p}^\mathrm{iso}(s)$ for all primes~$p$ and thus
  also yields a generalisation of Theorem~\ref{x3example}.  Indeed,
  for $d = 2$ the zeta function
  $\zeta^{\wedge}_{\widetilde{\Gamma},p}(s)$ has abscissa of
  convergence $115/21$ and for $d \ge 3$ it has abscissa of
  convergence $30d/11$.  Whenever $p$ is unramified in $k$, the local
  zeta function satisfies the functional equation
  \[
    \zeta^{\wedge}_{\widetilde{\Gamma},p}(s) \vert_{p\to p^{-1}} = \pm
    p^{24d^2+8d-10ds}\, \zeta^{\wedge}_{\widetilde{\Gamma},p}(s).
  \]
\end{remark}


\subsection{The local pro-isomorphic zeta functions of groups
  $\widetilde{\Gamma}$ associated to $L_t$}
The Lie lattice $L = L_t$ is decomposable and does not quite fit into
the same drawer as the lattices $L_{t^m}$, $m \ge 2$.  For
completeness we indicate how the approach in \cite{BGS22} can be
adapted in this and similar situtations to obtain the local
pro-isomorphic zeta functions of class-two nilpotent groups
$\widetilde{\Gamma}$ associated to Lie lattices $\widetilde{L}$
obtained from $L$ via `base extension'.

We start our discussion more generally.  Let $L$ be any class-two
nilpotent $\Z$-Lie lattice, and throughout let $p$ denote a rational
prime.  Then $[L,L] \subseteq \mathrm{Z}(L)$ and we can decompose $L$
as a direct sum $L = L^\circ \oplus M$ of Lie sublattices, where
$L^\circ$ satisfies $[L^\circ,L^\circ] = \mathrm{Z}(L^\circ) = [L,L]$
and $M \subseteq \mathrm{Z}(L)$ is abelian.  Typically there are many
choices for $L^\circ$ and $M$, but both are uniquely determined up to
isomorphism.  We set $l = \dim_\Z(L^\circ/[L^\circ,L^\circ])$,
$m = \dim_\Z(M)$ and $n = \dim_\Z([L^\circ,L^\circ])$.  The algebraic
automorphism group $\mathbf{Aut}(L^\circ)$ can be realised as a
subgroup scheme $\G \le \GL_{l+n}$ via a $\Z$-basis
$\mathcal{S}^\circ = (x_1, \ldots, x_l, z_1, \ldots, z_n)$ for
$L^\circ$ consisting of a basis $x_1, \ldots, x_l$ for a complement of
$[L^\circ,L^\circ]$ in $L^\circ$ and a basis $z_1, \ldots, z_n$ for
$[L^\circ,L^\circ]$.  Similarly we view $\mathbf{Aut}(L)$ as a
subgroup scheme of $\GL_{l+m+n}$ via the extended $\Z$-basis
\[
  \mathcal{S} = (x_1, \ldots, x_l, y_1, \ldots, y_m, z_1, \ldots,
  z_n),
\]
where $y_1, \ldots, y_m$ form a basis for $M$.

There are polynomial conditions, which we will denote by $(\dagger)$,
and a polynomial map $f$ from $\GL_l$ to $\GL_n$, which can be made
explicit in terms of the structure constants of the Lie lattice
$L^\circ$, such that automorphisms of the $\Q_p$-Lie algebra
${_{\Q_p} L^\circ} = \Q_p \otimes_\Z L^\circ$, viewed as elements of
the group $G_p = \G(\Q_p) \le \GL_{l+n}(\Q_p)$, take the form
\begin{equation}\label{equ:shape-L-L0-a}
  \left(
    \begin{array}{c|c}
      A & * \\
      \hline
        & f(A)
    \end{array}
  \right), \quad
  \begin{array}{l} \text{where
    $A \in \GL_l(\Q_p)$ satisfies $(\dagger)$} \\
    \text{ and $*$ is a placeholder for arbitrary entries.}
  \end{array}
\end{equation}
Moreover, automorphisms of ${_{\Q_p} L} = \Q_p \otimes_\Z L$, viewed
as elements of $\GL_{l+m+n}(\Q_p)$, take the form
\begin{equation}\label{equ:shape-L-L0-b}
  \left(
    \begin{array}{c|c|c}
      A &* & *\\
      \hline
        & B & * \\
      \hline
        && f(A)
    \end{array}
  \right), \quad
  \begin{array}{l} \text{where
    $A \in \GL_l(\Q_p)$ satisfies $(\dagger)$, $B \in \GL_m(\Q_p)$ has no
    particular} \\
    \text{ restrictions and $*$ is a placeholder for arbitrary entries.}
  \end{array}
\end{equation}

We set ${_p L^\circ} = \Z_p \otimes_\Z L^\circ$ and
${_p L} = \Z_p \otimes_\Z L$. As described in
Section~\ref{section:reduction-integral}, the local zeta function
$\zeta_{L^\circ,p}^\wedge(s) = \zeta_{L^\circ_p}^{\mathrm{iso}}(s)$ is
under suitable assumptions given by an integral
\begin{equation} \label{equ:integral-L0}
  \zeta_{L^\circ_p}^\mathrm{iso}(s) = \int_{H_p^+} \lvert \det h
  \rvert_p^{\, s} \, \theta_0(h) \, \theta_1(h) \,
  \mathrm{d}\mu_{H_p}(h)
\end{equation}
over $H_p^+ = H_p \cap \Mat_{l+n}(\Z_p)$, where $H_p$ denotes the
reductive part of~$G_p$ and
$\theta_0, \theta_1 \colon H_p \to \R_{\ge 0}$ are suitable volume
functions.  The descriptions in \eqref{equ:shape-L-L0-a} and
\eqref{equ:shape-L-L0-b} provide a close link between the groups of
$\Q_p$-points of the algebraic automorphism groups of $L^\circ$
and~$L$; for instance, the reductive parts are $H_p$ and
$H_p \times \GL_m(\Q_p)$, up to isomorphism.  Utilizing this
connection, we obtain for the local zeta function
$\zeta_{L,p}^\wedge(s) = \zeta_{L_p}^{\mathrm{iso}}(s)$ the integral
formula
\begin{equation} \label{equ:integral-L-from-L0-GLm}
  \zeta_{L_p}^\mathrm{iso}(s) %
  = \int_{H_p^+} \lvert \det h \rvert_p^{\, s} \, \theta_0(h) \,
  \theta_1(h)^{(l+m)/l} \, \mathrm{d}\mu_{H_p}(h) \;\cdot\; %
  \int_{\GL_m(\Q_p)^+} \lvert \det g \rvert_p^{\, s-l} \,
  \mathrm{d}\mu_{\GL_m(\Q_p)}(g),
\end{equation}
in analogy to~\eqref{equ:integral-p-pfrak}; the first factor is a mild
modification of the integral in~\eqref{equ:integral-L0} and
accomodates for the extra middle block in the third column of
\eqref{equ:shape-L-L0-b}, the second factor accommodates for the extra
blocks in the middle column in~\eqref{equ:shape-L-L0-b}.  The second
factor in \eqref{equ:integral-L-from-L0-GLm} is well-known and easy to
compute; one gets
\[
  \int_{\GL_m(\Q_p)^+} \lvert \det g \rvert_p^{\, s-l} \,
  \mathrm{d}\mu_{\GL_m(\Q_p)}(g) = \prod_{j=1}^m \big( 1-p^{(l+j-1)-s}
  \big)^{-1}.
\]

\begin{example}
  Let $L = L^\circ \oplus M$, where
  $L^\circ = \Z x_1 \oplus \Z x_2 \oplus \Z z_1$ with
  $[x_1,x_2] = z_1$ denotes the Heisenberg Lie lattice and
  $M = \bigoplus_{i=1}^m \Z y_i \cong \Z^m$ is abelian.  Then
  $(l,n) = (2,1)$ and
  $H_p = \{ \diag(A,\det A) \mid A \in \GL_2(\Q_p) \} \le
  \GL_3(\Q_p)$, furthermore $\theta_0(h) = 1$ and
  $\theta_1(h) = \lvert \det A \rvert_p^{\, -2}$ for
  $h = \diag(A,\det A) \in H_p$ in the above approach; consequently,
  we obtain
  \begin{align*}
    \zeta_{L_p}^\mathrm{iso}(s) %
    & = \int_{\GL_2(\Q_p)^+} \lvert \det A \rvert_p^{\, 2s-m-2} \,
      \mathrm{d}\mu_{\GL_2(\Q_p)}(h) \;\cdot%
      \; \prod_{j=1}^m \big( 1-p^{(j+1)-s} \big)^{-1} \\
    & = \big(1 - p^{(m+2)-2s} \big)^{-1} \big(1 - p^{(m+3)-2s}
      \big)^{-1} \; \prod_{j=1}^m \big( 1-p^{(j+1)-s} \big)^{-1},
  \end{align*}
  in agreement with the calculation in~\cite[\S 3.3.4]{Be05}.  The
  pro-isomorphic zeta function of the corresponding group
  $\Gamma \cong \mathrm{Heis}(\Z) \times C_\infty^{\, m}$ is a product
  of shifted Riemann zeta functions.
\end{example}
  
In a situation, where $L^\circ$ satisfies rigidity conditions of the
kind described in Section~\ref{sec:local-rgidity} the approach
discussed in Section~\ref{sec:local-Lt2-Lt3} can easily be adapted to
yield a formula for the local zeta functions
$\zeta_{\widetilde{L},p}^\wedge(s) =
\zeta_{\widetilde{L}_p}^\mathrm{iso}(s)$ of Lie lattices
$\widetilde{L}$ obtained from $L$ via `base extension'.  As before let
$k$ be a number field of absolute degree $d = [k : \Q]$, with ring of
integers~$\mathfrak{o}$.  We apply the method and notation from
Section~\ref{sec:local-Lt2-Lt3} to $L^\circ$ in place of~$L$.  Let
$\widetilde{L} = {_{\Z,\mathfrak{o}} L}$ be the class-two nilpotent
$\Z$-Lie lattice obtained by extending and restricting scalars.  We
observe that
$\widetilde{L} = \widetilde{L^\circ} \oplus \widetilde{M}$ is a
suitable decomposition of $\widetilde{L}$ in the sense given at the
beginning of this section, with
$\widetilde{L}^\circ = \widetilde{L^\circ}$ satisfying
$[\widetilde{L}^\circ,\widetilde{L}^\circ] =
[\widetilde{L},\widetilde{L}]$ and
$\widetilde{M} \subseteq \mathrm{Z}(\widetilde{L})$ abelian;
furthermore, $\dim_\Z(\widetilde{L}^\circ) = (l+n)d$ and
$\dim_\Z(\widetilde{M}) = md$.  Combining the approaches taken in this
and the previous section, we arrive at the integral formula
\begin{align*}
  \zeta_{\widetilde{L}_p}^\mathrm{iso}(s) %
  & = \left( \prod_{\mathfrak{p}
    \mid p} \int_{H_\mathfrak{p}^+} \lvert \det h
    \rvert_\mathfrak{p}^{\, s} \, \theta_0(h) \,
    \theta_1(h)^{(l+m)d/l} \, \mathrm{d}\mu_{H_\mathfrak{p}}(h)
    \right) \;\cdot %
    \;\int_{\GL_{md}(\Q_p)^+} \lvert \det g \rvert_p^{\, s-ld} \,
    \mathrm{d}\mu_{\GL_{md}(\Q_p)}(g) \\
  & = \left( \prod_{\mathfrak{p}
    \mid p} \int_{H_\mathfrak{p}^+} \lvert \det h
    \rvert_\mathfrak{p}^{\, s} \, \theta_0(h) \,
    \theta_1(h)^{(l+m)d/l} \, \mathrm{d}\mu_{H_\mathfrak{p}}(h) \right)
    \;\cdot %
    \;\prod_{j=1}^{md} \big( 1-p^{(ld+j-1)-s} \big)^{-1},
\end{align*}
where $H_\mathfrak{p}$, $H_\mathfrak{p}^+$, $\theta_0$ and $\theta_1$
(again without the decoration `$\mathfrak{p}$') are defined as in
preparation for~\eqref{equ:integral-p-pfrak}, but applied to $L^\circ$
instead of~$L$.

\smallskip

Finally we apply the above considerations to
$L = L_t = L^\circ \oplus M$, where
$L^\circ = \Z x_1 \oplus \Z x_2 \oplus \Z z_1$ with $[x_1,x_2] = z_1$
is the Heisenberg Lie lattice and $M = \Z y_1$ is abelian of rank~$1$.
In this case $(l,m,n) = (2,1,1)$ and our general formula specialises
to
\begin{align*}
  \zeta_{\widetilde{L}_p}^\mathrm{iso}(s) %
  & = \left( \prod_{\mathfrak{p}
    \mid p} \int_{\GL_2(k_\mathfrak{p})^+} \lvert \det A
    \rvert_\mathfrak{p}^{\, 2s-3d} \, 
    \mathrm{d}\mu_{\GL_2(k_\mathfrak{p})}(h) \right) \;\cdot%
    \; \prod_{j=1}^d \big( 1-p^{(2d+j-1)-s} \big)^{-1} \\
  & = \left( \prod_{\mathfrak{p}
    \mid p} \big(1 - q_\mathfrak{p}^{\, 3d-2s} \big)^{-1}
    \big(1 - q_\mathfrak{p}^{\, (3d+1) -2s}
    \big)^{-1} \right) \;\cdot\; \prod_{j=1}^d \big( 1-p^{(2d+j-1)-s} \big)^{-1},
\end{align*}
where $q_\mathfrak{p}$ denotes the residue field size of
$k_\mathfrak{p}$.  The pro-isomorphic zeta function of the
corresponding class-two nilpotent group $\widetilde{\Gamma}$ is a
product of two shifts of the Dedekind zeta function of $k$ and $d$
shifted Riemann zeta functions.

\begin{theorem} \label{thm:t-base-change} Let $k$ be a number field of
  absolute degree $d = [k:\Q]$, with ring of integers~$\mathfrak{o}$.
  Let $\widetilde{\Gamma} = \widetilde{\Gamma}_{t,k}$ be the class-two
  nilpotent group of Hirsch length $4d$ and with rank-$2d$ centre,
  corresponding to the class-two nilpotent $\Z$-Lie lattice
  $\widetilde{L} = \widetilde{L}_{t,k}$ which results from the Lie
  lattice $L = L_{t}$ via `base extension' as defined above.

  Then the pro-isomorphic zeta function of the group
  $\widetilde{\Gamma}$ is
  \begin{equation*} 
    \zeta^{\wedge}_{\widetilde{\Gamma}}(s) = \zeta_k(2s-3d)
    \, \zeta_k(2s-(3d+1)) \cdot \prod_{j=1}^d \zeta_\Q(s-(2d+j-1)),
  \end{equation*}
  where $\zeta_k(s)$ denotes the Dedekind zeta function of~$k$ and
  $\zeta_\Q(s)$ denotes the Riemann zeta function; in particular, it
  admits meromorphic continuation to the entire complex plane.
\end{theorem}

\begin{remark} The abscissa of convergence of
  $\zeta^{\wedge}_{\widetilde{\Gamma}}(s)$ is~$3d$, with a single pole
  at $s=3d$.  The asymptotic growth of pro-isomorphic subgroups in
  $\widetilde{\Gamma}$ can be described by means of a suitable
  Tauberian theorem.  The local zeta function
  $\zeta^{\wedge}_{\widetilde{\Gamma},p}(s)$ has abscissa of
  convergence $3d-1$ and, if $p$ is unramified in $k$, it satisfies
  the functional equation
  \[
    \zeta^{\wedge}_{\widetilde{\Gamma},p}(s) \vert_{p\to p^{-1}} = \pm
    p^{8d^2+\frac{1}{2}d(d+1)-5ds}\, \zeta^{\wedge}_{\widetilde{\Gamma},p}(s).
  \]
\end{remark}
  

\end{document}